\documentclass[final,1p,times]{elsarticle}

\journal{Computer Physics Communications}

\usepackage{upgreek}

\usepackage{graphicx}
\usepackage{amsmath,amsfonts,amssymb,bm,times,mathtools}
\usepackage{pmbb-sym}
\usepackage{listings}

\usepackage{algorithm}
\usepackage{algorithmic}

\usepackage{amsthm}

\newtheorem{remark}{Remark}

\newtheorem{formula}{Formula}

\usepackage{subfigure}
\usepackage{booktabs}
\usepackage[colorlinks,bookmarksopen,bookmarksnumbered,citecolor=red,urlcolor=red,breaklinks=true]{hyperref}

\usepackage{tikz}
\usetikzlibrary{positioning,shapes,arrows,arrows.meta}

\usepackage{here}

\usepackage{pgfplots}

\def\pp#1#2{\frac{\partial #1}{\partial #2}}
\def\abs#1{\left|{#1}\right|}
\def\bm#1{\boldsymbol{#1}}

\def\ksnrsoeji#1#2#3#4{\setbox0\hbox{$\displaystyle #1$}\mathop{\hbox to
   \wd0{\hss$ \displaystyle #2 $\hss}\kern-\wd0 \hbox{$\displaystyle
     #1_{#3}^{#4}$}}} 
\def\pfint_#1^#2{\ksnrsoeji{\int}={#1}{#2}}
\def\vpint_#1^#2{\ksnrsoeji{\int}-{#1}{#2}}

\def\Order{\mathcal{O}}

\def\XXint#1#2#3{{\setbox0=\hbox{$#1{#2#3}{\int}$}
\vcenter{\hbox{$#2#3$}}\kern-.5\wd0}}

\def\imath{\textrm{i}}

\def\diff{\mathrm{d}}

\def\uin{u^{\rm in}} 

\iffalse
\usepackage{numprint}
\npthousandsep{} 
\npdecimalsign{.}
\npproductsign{\times}
\npaddplusexponent
\npaddmissingzero
\else
\usepackage{siunitx}
\NewDocumentCommand\numprint{m}{\num[round-mode = places]{#1}}
\NewDocumentCommand\nprounddigits{m}{\sisetup{round-precision = #1}}
\def\npproductsign#1{} 
\fi

\def\Ns{N_{\rm s}}
\def\Nt{N_{\rm t}}
\def\Ds{\Delta_{\rm s}}
\def\Dt{\Delta_{\rm t}}

\def\Hs{h_{\rm s}}
\def\Ht{h_{\rm t}}

\def\Ps{{p_{\rm s}}}
\def\Pt{{p_{\rm t}}}

\def\IL{{\mathcal I}} 

\def\diff{\mathrm{d}}

\def\Order{\mathcal{O}}

\def\red#1{\textcolor{red}{#1}}

\def\blue#1{\textcolor{blue}{#1}}
\def\magenta#1{\textcolor{magenta}{#1}}
\def\cyan#1{\textcolor{cyan}{#1}}
\def\darkgreen#1{\textcolor[rgb]{0.0,0.5,0.0}{#1}}

\def\binom#1#2{\genfrac{(}{)}{0pt}{}{#1}{#2}}

\def\matL{\mathbf{L}}

\def\matU{\mathbf{U}}

\def\mat#1{\mathbf{#1}}

\usepackage{framed} 

\iftrue 
\usepackage{ulem}

\def\del#1{\sout{#1}}
\usepackage{soul,xcolor}
\setstcolor{cyan}
\else

\def\del#1{}

\fi

\begin{document}

\begin{frontmatter}
  
  \title{%
    An enhancement of the fast time-domain boundary element method for the three-dimensional wave equation
  }%
  
  
  \author[NU]{Toru Takahashi\corref{cor}}
  \ead{toru.takahashi@mae.nagoya-u.ac.jp}
  \cortext[cor]{Corresponding author}
  \author[SMZ]{Masaki Tanigawa}
  \author[NU]{Naoya Miyazawa}
  
  \address[NU]{Department of Mechanical Systems Engineering, Nagoya University, Furo-cho, Chikusa-ku, Nagoya city, Aichi, 464-8603 Japan}
  \address[SMZ]{Institute of Technology, Shimizu Corporation, 3-4-17 Etchujima, Koto-ku, Tokyo 135-8530 Japan}
  
  \begin{abstract}
    Our objective is to stabilise and accelerate the time-domain boundary element method (TDBEM) for the three-dimensional wave equation. To overcome the potential time instability, we considered using the Burton--Miller-type boundary integral equation (BMBIE) instead of the ordinary boundary integral equation (OBIE), which consists of the single- and double-layer potentials. In addition, we introduced a smooth temporal basis, i.e. the B-spline temporal basis of order $d$, whereas $d=1$ was used together with the OBIE in a previous study~\cite{takahashi2014}. Corresponding to these new techniques, we generalised the interpolation-based fast multipole method that was developed in \cite{takahashi2014}. In particular, we constructed the multipole-to-local formula (M2L) so that even for $d\ge 2$ we can maintain the computational complexity of the entire algorithm, i.e. $O(\Ns^{1+\delta}\Nt)$, where $\Ns$ and $\Nt$ denote the number of boundary elements and the number of time steps, respectively, and $\delta$ is theoretically estimated as $1/3$ or $1/2$. The numerical examples indicated that the BMBIE is indispensable for solving the homogeneous Dirichlet problem, but the order $d$ cannot exceed 1 owing to the doubtful cancellation of significant digits when calculating the corresponding layer potentials. In regard to the homogeneous Neumann problem, the previous TDBEM based on the OBIE with $d=1$ can be unstable, whereas it was found that the BMBIE with $d=2$ can be stable and accurate. The present study will enhance the usefulness of the TDBEM for 3D scalar wave problems.
  \end{abstract}
  
  
\end{frontmatter}

\section{Introduction}\label{s:intro}

It is often necessary to analyse or simulate a wave phenomenon in an open space or external domain rather than in a closed space in physical applications. Further, in comparison with steady-state or frequency-domain approach, transient or time-domain approach is relatively useful because it can obtain even the frequency response through a Fourier analysis, although the time insatiability must be dealt with in the time domain analysis. So, developing a fast and stable computational method that is applicable to external problems in time domain is important and remains a challenge to be surmounted. This study proposes a noble method based on the boundary element method (BEM) regarding the 3D wave equation.

The BEM or boundary integral equation method is invaluable in the numerical analysis of exterior boundary value problems of (linear) partial differential equations in classical physics, in particular, acoustics and electromagnetics as well as elastodynamics.
This is because BEM can handle the (semi-)infinite domain without any approximation, which is the so-called absorbing boundary condition (ABC) in terms of domain-type solvers such as the finite element method.  However, the major drawback of the BEM is its very high computational cost. Nevertheless, since the emergence of the fast multipole method (FMM) proposed by Greengard and Rokhlin in 1987~\cite{greengard1987}, similar fast algorithms have been developed to accelerate the BEM for various types of problems, especially for steady-state (or frequency-domain) wave problems~\cite{nishimura2002,liu2009book,liu2011review}.

On the other hand, the acceleration of unsteady-state or time-domain BEM (TDBEM) has not been widely investigated thus far because of the additional efforts needed to take account of the time axis. The pioneering studies were performed by Michielssen's group around 2000~\cite{chew2001book}. In 2D/3D acoustics~\cite{ergin1998,ergin1999a,ergin1999c,ergin2000,lu2004,lu2004b} and 3D electromagnetics~\cite{shanker2003,aygun2004}, they developed the planewave time-domain (PWTD) algorithm as a time-domain version of the FMM~\cite{greengard1987}. In the 3D problems we are interested in, the PWTD algorithm can reduce the computational complexity from $O(\Ns^2\Nt)$ to $O(\Ns\log^n\Ns\cdot\Nt)$, where $\Ns$ and $\Nt$ denote the spatial and temporal degrees of freedom, respectively, and $n$ depends on the details of the planewave expansion. (In this study, $\Ns$ and $\Nt$ represent the number of boundary elements and the number of time steps, respectively.) Afterwards, the PWTD algorithm was applied to 2D/3D elastodynamic problems~\cite{takahashi2001,takahashi2003} and enhanced by using the wavelet~\cite{liu2014parallel_wavelet_PWTD}.

As a variant of the PWTD algorithm, an interpolation-based FMM for the 3D wave equation was proposed in the previous study~\cite{takahashi2014}. Because the interpolation can readily realise the separation of variables of the retarded layer potentials, the formulation and implementation of the interpolation-based FMM are simpler than those of the PWTD. The trade-off for the simplification is that the computational complexity of the fast TDBEM based on the interpolation-based FMM is relatively high, i.e. $O(\Ns^{1+\delta}\Nt)$, where $\delta$ is theoretically estimated as $1/3$ when boundary elements are distributed uniformly in 3D space and $1/2$ when they are on a plane in 3D. Although the interpolation-based FMM was not compared with the PWTD algorithm in \cite{takahashi2014}, the fast TDBEM using the interpolation-based FMM outperformed the conventional TDBEM, whose complexity is $O(\Ns^2\Nt)$, in the numerical test.

Whether the algorithm is fast or not, late-time instability is an important problem in the TDBEM. In general, there are several approaches to addressing this problem~\cite{hargreaves2007}. Probably the most implementation-friendly approach would be the $\alpha$-$\delta$ method~\cite{soares2007}. In this method, weights are introduced in the discretised boundary data over some successive time steps so that the possible oscillation of the boundary data can be averaged in those time steps~\cite{soares2007}. Okamura et al.~\cite{okamura2016} applied the $\alpha$-$\delta$ method to the fast TDBEM ~\cite{takahashi2014} and observed that the instability could be suppressed in their numerical analysis, although the numerical accuracy was not discussed quantitatively.

A more fundamental approach is to use the Burton--Miller-type boundary integral equation (BMBIE) instead of the ordinary BIE (OBIE), which consists of the single- and double-layer potentials. Ergin et al.~\cite{ergin1999b} found that the BMBIE, which is a linear combination of the normal and temporal derivatives of the OBIE, can be more stable than the OBIE. This can be regarded as the time-domain counterpart of removing the interior resonance in frequency domain~\cite{burton_miller1971}. The same authors accelerated the BMBIE using the two-level and multi-level PWTD algorithms~\cite{ergin1999a,ergin1999c}. Recently, Fukuhara et al.~\cite{fukuhara2019} analysed the stability of the TDBEM based on various types of BIE, including the BMBIE, for certain initial boundary value problems involving the 2D wave equation. They demonstrated that the distribution of the eigenvalues, which were obtained by means of the Sakurai--Sugiura method~\cite{asakura2009}, can differ according to the type of BIE, and that if the imaginary parts of all the eigenvalues are less than $0$, the corresponding BIE is stable. Their results suggest that the OBIE is unstable but the BMBIE is stable for the homogeneous Dirichlet problem in 2D. It is also implied that adding the single-layer potential (multiplied by a real constant) to the BMBIE can increase its accuracy. Similarly to Fukuhara et al.~\cite{fukuhara2019}, the 3D case was examined by Chiyoda et al.~\cite{chiyoda2019}.

The present study aims to stabilise the fast TDBEM~\cite{takahashi2014} for the 3D wave equation, maintaining the benefit obtained from the acceleration by the interpolation-based FMM. Following the former studies~\cite{ergin1999b,fukuhara2019,chiyoda2019}, we adopt the BMBIE instead of the OBIE, which was used in \cite{takahashi2014}, and incorporate the BMBIE into the FMM. In addition, we consider high-order discretisation with respect to time. That is, we adopt the B-spline temporal basis of order $d$ ($\ge 2$) instead of the piecewise-linear basis, which corresponds to $d=1$. We can evaluate the space-time integrals associated with both the discretised OBIE and the BMBIE for $d\ge 2$ similarly to how it was evaluated for $d=1$. Moreover, we need to generalise the interpolation-based FMM from $d=1$ to higher $d$'s, as the original FMM for $d=1$ fails because a certain linearity resulting from $d=1$ is no longer available for $d\ge 2$ in the multipole-to-local translation (M2L) of the FMM.

The rest of this paper is organised as follows: Section~\ref{s:biem} shows the formulation of the TDBEM based on the BMBIE and B-spline temporal basis. In Section~\ref{s:fmm2}, the interpolation-based FMM is constructed for the TDBEM formulated in the previous section. The details of constructing an efficient M2L are described in Section~\ref{s:m2l}. Sections~\ref{s:num} and \ref{s:app} numerically assess the computational stability and efficiency of the proposed fast TDBEM and demonstrate the applicability to parameter optimisation, respectively. Section~\ref{s:conclusion} concludes the paper.

\section{A TDBEM regarding BMBIE and B-spline temporal basis}\label{s:biem}

We formulate a TDBEM regarding the BMBIE as well as the OBIE when the B-spline basis of order $d$ is used as the temporal basis instead of the piecewise-linear basis, which exactly corresponds to the case $d=1$ and was used in the previous work~\cite{takahashi2014}.

\subsection{Problem statement}

Let $V$ be a finite domain in $\bbbr^3$ with the piecewise-smooth boundary $S:=\partial V$. We consider the following exterior problems of the 3D wave equation regarding the wave field or sound pressure $u$:
\begin{subequations}
  \begin{eqnarray}
    \triangle u(\bm{x},t)=\frac{1}{c^2}\frac{\partial^2 u}{\partial t^2}(\bm{x},t)&&\text{for $\bm{x}\in\bbbr^3\setminus\bar{V}$, $t>0$},\label{eq:wave3d}\\
    u(\bm{x},0)=\frac{\partial u}{\partial t}(\bm{x},0)=0 &&\text{for $\bm{x}\in\bbbr^3\setminus\bar{V}$},\\
    u(\bm{x},t)=\bar{u}(\bm{x},t)&&\text{for $\bm{x}\in S_u$, $t>0$},\\
    q(\bm{x},t):=\frac{\partial u}{\partial n_y}(\bm{x},t)=\bar{q}(\bm{x},t)&&\text{for $\bm{x}\in S_q$, $t>0$},\\
    u(\bm{x},t)\rightarrow 0 && \text{as $\abs{\bm{x}}\rightarrow\infty$, $t>0$},
  \end{eqnarray}%
  \label{eq:ibvp}%
\end{subequations}
where $\bm{n}$ is the unit outward normal, $S \equiv S_u \cup S_q$ and $S_u \cap S_q \equiv \emptyset$. Further, $c>0$ is the wave velocity. In the numerical examples in Sections~\ref{s:num} and \ref{s:app}, we will consider the incident wave, which will be denoted by $\uin$, but we omit it from the formulation in this and the following sections for the sake of simplicity.

\subsection{BIEs}

To solve the initial-boundary value problem in (\ref{eq:ibvp}), the first --- and simplest --- choice is to use the following OBIE, which consists of the single- and double-layer potentials:
\begin{eqnarray}
  \frac{1}{2}u(\bm{x},t)=\int_0^t\int_S\left(
  \Gamma(\bm{x}-\bm{y},t-\tau)q(\bm{y},\tau)
  -\frac{\partial\Gamma}{\partial n_y}(\bm{x},\bm{y},t-\tau)u(\bm{y},\tau)
  \right)\diff\tau\diff S_y\quad\text{for $\bm{x}\in S$ and $t>0$},
  \label{eq:bie}
\end{eqnarray}
where, with Dirac's delta function $\delta$, $\Gamma$ is the causal fundamental solution of (\ref{eq:wave3d}) given by
\begin{eqnarray*}
  \Gamma(\bm{x},t):=\frac{\delta(t-|\bm{x}|/c)}{4\pi |\bm{x}|}.
\end{eqnarray*}
However, the OBIE can suffer from the interior resonance problem~\cite{ergin1999b,fukuhara2019}. To avoid it, we consider the following BMBIE, which is obtained by applying the normal and temporal derivatives, i.e. $\frac{\partial}{\partial n_x}-\frac{1}{c}\frac{\partial}{\partial t}$, to the OBIE:
\begin{eqnarray}
  &&\frac{1}{2}\left(\pp{u(\bm{x},t)}{n_x}-\frac{1}{c}\pp{u(\bm{x},t)}{t}\right)
  =\int_0^t\int_S\left(\pp{\Gamma(\bm{x}-\bm{y},t-\tau)}{n_x}q(\bm{y},\tau)-\frac{\partial^2\Gamma(\bm{x}-\bm{y},t-\tau)}{\partial n_x\partial n_y}u(\bm{y},\tau)\right)\diff\tau\diff S_y\nonumber\\
  &&-\frac{1}{c}\frac{\partial}{\partial t}\int_0^t\int_S\left(
    \Gamma(\bm{x}-\bm{y},t-\tau)q(\bm{y},\tau)
    -\frac{\partial\Gamma}{\partial n_y}(\bm{x},\bm{y},t-\tau)u(\bm{y},\tau)
    \right)\diff\tau\diff S_y\quad\text{for $\bm{x}\in S$ and $t>0$},
    \label{eq:bmbie}
\end{eqnarray}

We use the collocation method to solve both BIEs. The alternative is the Galerkin method~\cite{aimi2009,joly2017,gimperlein2018}, but we do not consider it in this study.

\subsection{Discretisation of the BIEs}

To discretise $u$ and $q$ in the BIEs, we use the piecewise-constant basis with respect to space. The boundary $S$ is discretised with $\Ns$ triangular boundary elements, denoted by $E_i$ where $i\in[1,\Ns]$.  Then, $u(\bm{x},t)$ and $q(\bm{x},t)$ on $E_i$ are denoted by $u_i(t)$ and $q_i(t)$, respectively. We let the centre of $E_i$ be the collocation point $\bm{x}_i$, where $i=1,\ldots,\Ns$.

As the temporal basis of $u_i$ and $q_i$, we utilise the B-spline bases of order $d$ ($\ge 1$). That is, we interpolate $u_i$ and $q_i$ as
\begin{eqnarray}
  \{u_i,q_i\}(t)\approx\sum_{\beta=0}N^{\beta,d}(t)\{u_i^\beta,q_i^\beta\},
  \label{eq:u_bspline}
\end{eqnarray}
where $N^{\beta,d}$ represents the $\beta$th B-spline basis of order $d$, and $u_i^\beta$ is its coefficient. As usual, in the TDBEM, we assume that the knots of the B-spline basis or temporal nodes, denoted by $t_0,t_1,\ldots$, are equidistant, i.e. $t_{i+1}-t_i\equiv\Dt$ or $t_i\equiv i \Dt$, where $\Dt$ denotes the time-step length. Then, we let the $\alpha$th temporal collocation point be $t_\alpha$ ($\equiv\alpha\Dt$) for $\alpha=1,2,\ldots,\Nt-1$, where $\Nt$ denotes the prescribed number of time steps.


\begin{remark}[Support]\label{remark:support}
  The support of the basis $N^{\beta,d}$ is $[t_\beta,t_{\beta+d+1}]$.
\end{remark}

\begin{remark}[Translational invariance]
  For any time $s\in\bbbr$, it holds that $N^{\beta,d}(t_\alpha+s)\equiv N^{\beta+\gamma,d}(t_\alpha+t_\gamma+s)$ if the knots $\{t_i\}$ are uniform. This means that the value of a basis (i.e. $N^{\beta,d}$) is determined by the difference between an evaluation time $t_\alpha$ and a source time $t_\beta$. 
\end{remark}

As with the piecewise-linear temporal basis or the case of $d=1$, it is desirable to analytically integrate the space-time integrals in the resulting discretised layer potentials. To simplify the integration, we decompose $N^{\beta,d}$ into $d+2$ truncated power functions of order $d$, i.e. $(\cdot)_+^d$, as follows:
\begin{eqnarray}
  N^{\beta,d}(x)=\sum_{\kappa=0}^{d+1} w^{\kappa,d} \left(\frac{x-t_{\kappa+\beta}}{\Dt}\right)^d_+,
  \label{eq:decompose}
\end{eqnarray}
where $w^{\kappa,d}$ is defined by
\begin{eqnarray}
  w^{\kappa,d}:=\frac{d+1}{\prod_{k=0;k\ne \kappa}^{d+1}(k-\kappa)}
  \label{eq:w}
\end{eqnarray}
and have the values shown in Table~\ref{tab:w}. Hence, we can express $u_i$ and $q_i$ in (\ref{eq:u_bspline}) as follows:
\begin{eqnarray}
  \{u_i(t),q_i(t)\}\approx\sum_{\beta=0} \sum_{\kappa=0}^{d+1} w^{\kappa,d} \left(\frac{t-t_{\kappa+\beta}}{\Dt}\right)^d_+ \{u^\beta_i,q^\beta_i\}.
  \label{eq:u}
\end{eqnarray}

\begin{table}[H]
  \centering
  \caption{Values of Coefficients $w^{\kappa,d}$ in (\ref{eq:w}) for $d\le 4$.}
  \label{tab:w}
  \renewcommand{\arraystretch}{1.2}
  \begin{tabular}{ccccccc}
    \toprule
    $d$ $\backslash$ $\kappa$ & $0$ & $1$ & $2$ & $3$ & $4$ & $5$\\
    \midrule
    $1$ & $1$ & $-2$ & $1$ & na & na & na\\
    $2$ & $\frac{1}{2}$ & $-\frac{3}{2}$ & $\frac{3}{2}$ & $-\frac{1}{2}$ & na & na\\
    $3$ & $\frac{1}{6}$ & $-\frac{2}{3}$ & 1 & $-\frac{2}{3}$ & $\frac{1}{6}$ & na\\
    $4$ & $\frac{1}{24}$ & $-\frac{5}{24}$ & $\frac{5}{12}$ & $-\frac{5}{12}$ & $\frac{5}{24}$ & $-\frac{1}{24}$\\
    \bottomrule
  \end{tabular}
\end{table}

\begin{figure}[H]
  \centering
  \begin{tabular}{ccc}
    \includegraphics[width=.3\textwidth]{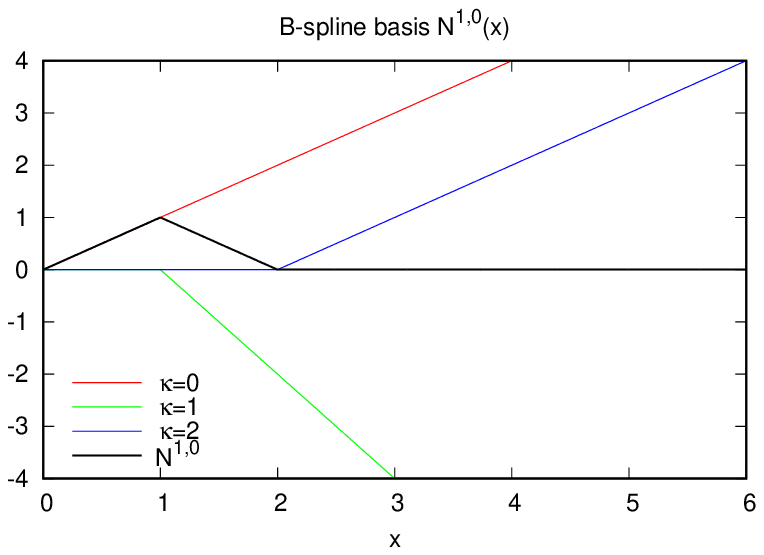}
    & \includegraphics[width=.3\textwidth]{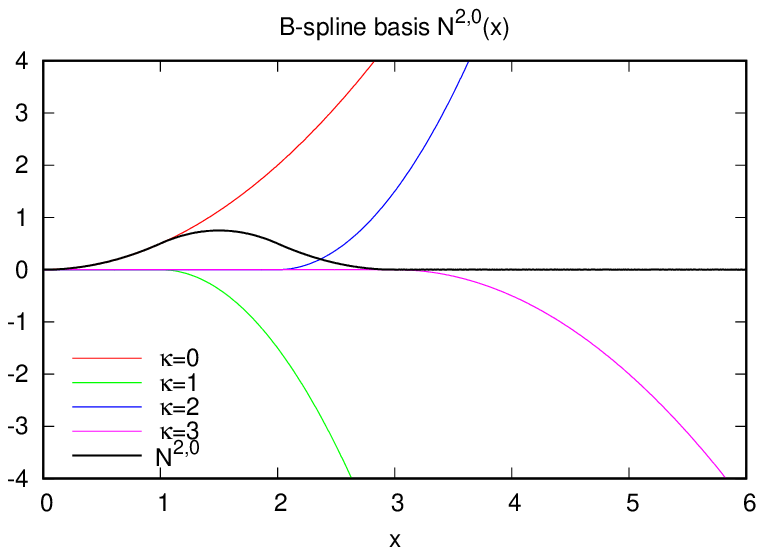}
    & \includegraphics[width=.3\textwidth]{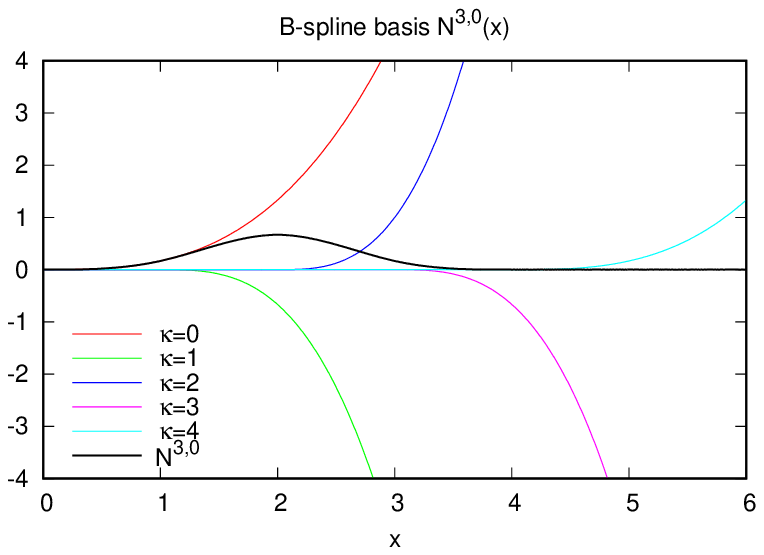}\\
    $d=1$ & $d=2$ & $d=3$
  \end{tabular}
  \caption{Decomposition of the B-spline basis $N^{d,0}$, whose support is $[t_0,t_{d+1}]$ from Remark~\ref{remark:support}, into $d+1$ truncated power functions $w^{\kappa,d}\left(\frac{t-t_{\kappa}}{\Dt}\right)^d_+$ of order $d$ according to (\ref{eq:decompose}). Here, $\Dt$ is given as $1$.}.
  \label{fig:decompose}
\end{figure}

Figure~\ref{fig:decompose} shows the decomposition of the B-spline basis $N^{d,0}$ into $d+1$ truncated power functions of order $d$, where $d$ is selected as 1, 2, or 3.

\subsection{Discretisation of the BIEs}\label{s:discretize OBIE}

With (\ref{eq:u}), we discretise both OBIE in (\ref{eq:bie}) and BMBIE in (\ref{eq:bmbie}). First, we consider the former and discretise the single-layer potential, denoted by $\Phi_{\rm single}$, as follows:
\begin{eqnarray*}
  \Phi_{\rm single}(\bm{x},t)
  &&:=\int_0^{t}\int_S \Gamma(\bm{x}-\bm{y},t-\tau)q(\bm{y},\tau) \diff\tau\diff S_y\nonumber\\
  &&\approx\sum_{j=1}^{\Ns}\int_{E_j}\frac{1}{4\pi\abs{\bm{x}-\bm{y}}}\sum_{\beta=0}\sum_{\kappa=0}^{d+1}w^{\kappa,d}\left(\frac{t-\frac{\abs{\bm{x}-\bm{y}}}{c}-t_{\kappa+\beta}}{\Delta_t}\right)^d_+ q_j^\beta\diff S_y
\end{eqnarray*}
For the collocation point $(\bm{x}_i,t_\alpha)$, we can obtain the following expression:
\begin{eqnarray}
  \Phi_{\rm single}(\bm{x}_i,t_\alpha)
  &=&\sum_{j=1}^{\Ns}\sum_{\beta=0}^{\alpha-1}\sum_{\kappa=0}^{d+1}w^{\kappa,d}\underbrace{\int_{E_j}\frac{(ct_{\alpha-\beta-\kappa}-\abs{\bm{x}_i-\bm{y}})_+^d}{4\pi(c\Dt)^d\abs{\bm{x}_i-\bm{y}}} \diff S_y}_{\displaystyle U_{ij}^{(\alpha-\beta-\kappa)}} q_j^\beta\nonumber\\
  &=&\sum_{j=1}^{\Ns}\sum_{\beta=0}^{\alpha-1}\sum_{\kappa=0}^{d+1}w^{\kappa,d} U_{ij}^{(\alpha-\beta-\kappa)} q_j^\beta,
     \label{eq:discretised_single_layer_potential}
\end{eqnarray}
where the index $\gamma$ in $U_{ij}^{(\gamma)}$ stands for an index of the time difference rather than an index of the time step itself.

Similarly, we can discretise the double-layer potential $\Phi_{\rm double}$ in (\ref{eq:bie}) as follows:
\begin{eqnarray}
  &&\Phi_{\rm double}(\bm{x}_i,t_\alpha)\nonumber\\
  &&:=\int_0^{t}\int_S\frac{\partial\Gamma}{\partial n_y}(\bm{x}_i,\bm{y},t_\alpha-\tau)u(\bm{y},\tau)\diff\tau \diff S_y\nonumber\\
  &&=\sum_{j=1}^{\Ns}\sum_{\beta=0}^{\alpha-1}\sum_{\kappa=0}^{d+1}w^{\kappa,d}
   \underbrace{\int_{E_j} \frac{\bm{n}\cdot(\bm{x}_i-\bm{y})}{4\pi (c\Dt)^d\abs{\bm{x}_i-\bm{y}}^3}\left[ d\left(ct_{\alpha-\beta-\kappa}-\abs{\bm{x}_i-\bm{y}}\right)^{d-1}_+ct_{\alpha-\beta-\kappa}-(d-1)\left(ct_{\alpha-\beta-\kappa} -\abs{\bm{x}_i-\bm{y}}\right)^d_+\right]\diff S_y}_{\displaystyle W_{ij}^{(\alpha-\beta-\kappa)}} u_j^\beta\nonumber\\
  &&=\sum_{j=1}^{\Ns}\sum_{\beta=0}^{\alpha-1}\sum_{\kappa=0}^{d+1}w^{\kappa,d} W_{ij}^{(\alpha-\beta-\kappa)} u_j^\beta,
 \label{eq:discretised_double_layer_potential}
\end{eqnarray}
where, denoting $c(t-s)$ by $T$, we consider the following normal derivative with respect to $\bm{y}$:
\begin{eqnarray*}
  \partial_p^y \left\{\frac{(T-\abs{\bm{x}-\bm{y}})^d_+}{\abs{\bm{x}-\bm{y}}}\right\}
  &=&\ \partial_p^y \left\{\left(\frac{cT}{\abs{\bm{x}-\bm{y}}}-1\right)^d_+ \abs{\bm{x}-\bm{y}}^{d-1}\right\}\\
  &=&\ d\left(T-\abs{\bm{x}-\bm{y}}\right)^{d-1}_+\frac{T(x_p-y_p)}{\abs{\bm{x}-\bm{y}}^3}-(d-1)\left(T -\abs{\bm{x}-\bm{y}}\right)^d_+ \frac{x_p-y_p}{\abs{\bm{x}-\bm{y}}^3}\quad\text{for $p=1,2,3$}.
\end{eqnarray*}

\begin{remark}[Non-positive index]
  $U_{ij}^{(\gamma)}=W_{ij}^{(\gamma)}=0$ for any $\gamma\le 0$ because of the property of the truncated power functions.
  \label{remark:non_positive_index}
\end{remark}

For later convenience, we introduce the following functions $U$ and $\bm{W}$ as the kernel functions of the integrals $U_{ij}^{(\alpha-\beta-\kappa)}$ in (\ref{eq:discretised_single_layer_potential}) and $W_{ij}^{(\alpha-\beta-\kappa)}$ in (\ref{eq:discretised_double_layer_potential}):
\begin{eqnarray}
  U(\bm{x},\bm{y},t,s):=&\frac{(c(t-s)-\abs{\bm{x}-\bm{y}})_+^d}{\abs{\bm{x}-\bm{y}}},\label{eq:U}\\
  \bm{W}(\bm{x},\bm{y},t,s):=&\nabla_y U(\bm{x},\bm{y},t,s) = \left(d\left(c(t-s)-\abs{\bm{x}-\bm{y}}\right)^{d-1}_+c(t-s)-(d-1)\left(c(t-s) -\abs{\bm{x}-\bm{y}}\right)^d_+\right)\frac{\bm{x}-\bm{y}}{\abs{\bm{x}-\bm{y}}^3}.\label{eq:W}
\end{eqnarray}

From (\ref{eq:discretised_single_layer_potential}) and (\ref{eq:discretised_double_layer_potential}), we can discretise the ordinary BIE in (\ref{eq:bie}) as follows:
\begin{eqnarray}
  \frac{1}{2}\sum_{\beta=0}^{\alpha-1} N^{\beta,d}(t_\alpha)u_i^\beta = \sum_{j=1}^{\Ns}\sum_{\beta=0}^{\alpha-1}\sum_{\kappa=0}^{d+1} w^{\kappa,d} \left(U_{ij}^{(\alpha-\beta-\kappa)} q_j^\beta - W_{ij}^{(\alpha-\beta-\kappa)} u_j^\beta \right)\quad\text{for $\alpha=1,2,\ldots$}
  \label{eq:bie_disc'}
\end{eqnarray}

We note that the free term in the LHS of (\ref{eq:bie_disc'}) can be included in the double-layer potential by considering that the collocation point $\bm{x}_i$ lies outside the domain $V$ but infinitesimally close to the boundary $S$.  Then, we can express (\ref{eq:bie_disc'}) in matrix form as follows:
\begin{eqnarray}
  \bm{0} = \sum_{\beta=0}^{\alpha-1}\sum_{\kappa=0}^{d+1} w^{\kappa,d} \left(\mat{U}^{(\alpha-\beta-\kappa)} \mat{q}^\beta - \mat{W}^{(\alpha-\beta-\kappa)} \mat{u}^\beta \right)\quad\text{for $\alpha=1,2,\ldots$},
  \label{eq:bie_disc}
\end{eqnarray}
where $\mat{U}^{(\gamma)},\mat{W}^{(\gamma)}\in\bbbr^{\Ns\times\Ns}$ and $\mat{u}^\beta,\mat{q}^\beta\in\bbbr^{\Ns}$.

It should be noted that the lower bound (starting index) `$0$' of the summation over $\beta$ in (\ref{eq:bie_disc}) can be replaced with a positive index $\beta^*$, which means that we can discard the information from all the passed time steps before $\beta^*$ at the current time step $\alpha$. To prove this, we first write down the coefficient of $q_j^\beta$ in (\ref{eq:discretised_single_layer_potential}), i.e. $\sum_{\kappa=0}^{d+1}w^{\kappa,d}U_{ij}^{(\alpha-\beta-\kappa)}$ (as well as that of $u_j^\beta$ in (\ref{eq:discretised_double_layer_potential}), i.e. $\sum_{\kappa=0}^{d+1}w^{\kappa,d}W_{ij}^{(\alpha-\beta-\kappa)}$) as
\begin{eqnarray*}
  \sum_{\kappa=0}^{d+1} w^{\kappa,d} U_{ij}^{(\alpha-\beta-\kappa)}
  =\int_{E_j}\frac{1}{4\pi\abs{\bm{x}_i-\bm{y}}} N^{\beta,d}\left(t_\alpha-\frac{\abs{\bm{x}_i-\bm{y}}}{c}\right) \diff S_y.
\end{eqnarray*}
Because the support of the B-spline basis $N^{\beta,d}$ is $[t_\beta, t_{\beta+d+1}]$ from Remark~\ref{remark:support}, the coefficient of $q_j$ is non-vanishing if
\begin{eqnarray*}
  t_\alpha-\frac{\abs{\bm{x}_i-\bm{y}}}{c} < t_{\beta+d+1}
  \quad\Longleftrightarrow\quad \beta > \alpha-\left(\frac{\abs{\bm{x}_i-\bm{y}}}{c\Dt}+d+1\right).
\end{eqnarray*}
Here, because the boundary $S$ is finite, the relative distance $\abs{\bm{x}_i-\bm{y}}$ is bounded by $\max_{\bm{x},\bm{y}\in S}\abs{\bm{x}-\bm{y}}$ for any collocation point $\bm{x}_i$. Therefore, we can determine the lower bound $\beta^*$ as follows:
\begin{eqnarray}
  \beta\ge \alpha-\underbrace{\left(\frac{\max_{\bm{x},\bm{y}\in S}\abs{\bm{x}-\bm{y}}}{c\Dt}+d+1\right)}_{\displaystyle \gamma^*}=\alpha-\gamma^*=:\beta^*.
\end{eqnarray}
\begin{remark}{}
  Because $\gamma^*=\alpha-\beta^*$ holds, $\gamma^*$ represents the upper bound of the time difference $\gamma$ ($:=\alpha-\beta$). That is, we may store the coefficients $\mat{U}^{(\gamma)}$ and $\mat{W}^{(\gamma)}$ for $\gamma=1$ to $\gamma^*$. In other words, we may let $\mat{U}^{(\gamma)}$ and $\mat{W}^{(\gamma)}$ be zero if $\gamma>\gamma^*$.
  \label{remark:gamma_bound}
\end{remark}

\begin{remark}{}
  In regard to the time-marching scheme with a constant time-step length, we can in general state that the coefficients $\mat{U}^{(\gamma)}$ and $\mat{W}^{(\gamma)}$ may not be stored for $\gamma>\frac{\Nt}{2}$ by the cast-forward algorithm~\cite{yoshikawa2003phd} without any additional computations. From this and Remark~\ref{remark:gamma_bound}, we may store $\mat{U}^{(\gamma)}$ and $\mat{W}^{(\gamma)}$ for $1\le\gamma\le\mathrm{min}\left(\gamma^*,\frac{\Nt}{2}\right)$. Nevertheless, the conventional BEM requires a very large memory to store those coefficients.
\end{remark}

Further, we can rewrite (\ref{eq:bie_disc}), where the index $\beta$ begins with $\beta^*$ rather than with $0$, in the following simpler form:
\begin{eqnarray}
  \bm{0}=\sum_{\beta=\beta^*+1}^{\alpha}\left(\mat{U}^{(\alpha-\beta+1)} \bm{\uptau}^{\beta-1} - \mat{W}^{(\alpha-\beta+1)} \bm{\upsigma}^{\beta-1}\right)\quad\text{for $\alpha=1,2,\ldots$},
  \label{eq:bie_disc3}
\end{eqnarray}
where, as described in \ref{s:simplify}, we introduced the new boundary variables $\bm{\uptau}$ and $\bm{\upsigma}$ by imposing the summation over the index $\kappa$ to the original boundary variables $\mat{u}$ and $\mat{q}$, respectively; see (\ref{eq:tau}) and (\ref{eq:sigma}).

Analogously to the OBIE, the BMBIE in (\ref{eq:bmbie}) can be reduced to the discretised form in (\ref{eq:bie_disc3}), where the coefficient matrices $\mat{U}^{(\gamma)}$ and $\mat{W}^{(\gamma)}$ are replaced with the following matrices $\overline{\mat{U}}^{(\gamma)}$ and $\overline{\mat{U}}^{(\gamma)}$, respectively:
\begin{subequations}
  \begin{eqnarray}
    \overline{\mat{U}}^{(\gamma)} &:=& \left(\frac{\partial}{\partial n_x}-\frac{1}{c}\frac{\partial}{\partial t}\right)\mat{U}^{(\gamma)},\\
    \overline{\mat{W}}^{(\gamma)} &:=& \left(\frac{\partial}{\partial n_x}-\frac{1}{c}\frac{\partial}{\partial t}\right)\mat{W}^{(\gamma)}.
  \end{eqnarray}%
  \label{eq:U_W_BM}%
\end{subequations}%
It should be noted that the differentiations on the RHSs can be evaluated analytically because we can perform the boundary (spatial) integrals in $\mat{U}^{(\gamma)}$ and $\mat{W}^{(\gamma)}$ analytically in the same manner as for $d=1$, i.e. the piecewise-linear temporal basis, which was investigated in \cite{yoshikawa2003phd}, for example. The details of the analytical integration are described in \ref{s:element_integral}.

\subsection{Linear equations}

At the current time step $t_\alpha$ ($\alpha=1,2,\ldots$), we solve the discretised BIE in (\ref{eq:bie_disc3}) (or the corresponding one for the BMBIE) for the unknown components in the vectors $\mat{u}^{\alpha-1}$ and $\mat{q}^{\alpha-1}$. To this end, we reduce (\ref{eq:bie_disc3}) to the following linear equations with respect to the unknown vector at $t=t_\alpha$, which is denoted by $\mat{x}^{\alpha-1}$ ($\subset\bbbr^{\Ns}$):
\begin{eqnarray}
  \mat{A}\mat{x}^{\alpha-1}
  =-\sum_{\beta=\beta^*+1}^{\alpha}\left(\mat{W}^{(\alpha-\beta+1)}\tilde{\bm{\upsigma}}^{\beta-1} - \mat{U}^{(\alpha-\beta+1)} \tilde{\bm{\uptau}}^{\beta-1}\right)=:\mat{b}^{\alpha-1},
  \label{eq:linear}
\end{eqnarray}
which is derived in \ref{s:linear}. Here, we defined the coefficient matrix $\bm{A}$ of the unknown vector $\bm{x}^{\alpha-1}$ as follows:
\begin{eqnarray*}
  A_{ij}:=w^{0,d}\left(W_{ij}^{(1)} \phi_j - U_{ij}^{(1)}(1-\phi_j)\right),
\end{eqnarray*}
where 
\begin{eqnarray*}
  \phi_j&:=&
  \begin{cases}
    1 & \text{if $S_j\subset S_q$; i.e. $u$ is unknown on $E_j$}\\
    0 & \text{if $S_j\subset S_u$; i.e. $q$ is unknown on $E_j$}
  \end{cases},\\
  x_j^{\alpha-1}&:=&u_j^{\alpha-1}\phi_j+q_j^{\alpha-1}(1-\phi_j).
\end{eqnarray*}
In addition, we define $\tilde{\bm{\upsigma}}^{\beta}$ and $\tilde{\bm{\uptau}}^{\beta}$ as follows:
\begin{subequations}
\begin{eqnarray}
    \tilde{\tau}_j^{\beta}
    &:=&\begin{cases}
    w^{0,d}q_j^\beta\phi_j+\displaystyle\sum_{\kappa=1}^{\min(d+1,\beta-\beta^*)}w^{\kappa,d} q_j^{\beta-\kappa} & \text{if $\beta=\alpha-1$}\\
    \tau_j^{\beta} & \text{if $\beta\in[0,\alpha-2]$}
    \end{cases},\\
    \tilde{\sigma}_j^{\beta}
    &:=&\left\{\begin{array}{ll}
                w^{0,d}u_j^\beta(1-\phi_j)+\displaystyle\sum_{\kappa=1}^{\min(d+1,\beta-\beta^*)}w^{\kappa,d} u_j^{\beta-\kappa} & \text{if $\beta=\alpha-1$}\\
     \sigma_j^\beta & \text{if $\beta\in[0,\alpha-2]$}
              \end{array}\right..
  \end{eqnarray}%
  \label{eq:sigma_tau_tilde}%
\end{subequations}
In other words, $\tilde{\bm{\upsigma}}^{\beta}$ and $\tilde{\bm{\uptau}}^{\beta}$ can be obtained from $\bm{\upsigma}^{\beta}$ in (\ref{eq:sigma}) and $\bm{\uptau}^{\beta}$ in (\ref{eq:tau}), respectively, by letting all their unknown components (at the current time step $t_\alpha$) be zero.

The solution of the linear equations in (\ref{eq:linear}) requires $O(\Ns\Nt)$ computational complexity through all the time steps $\Nt$, because the coefficient matrix $\bm{A}^{(1)}$ is sparse, and therefore (\ref{eq:linear}) is solvable with $O(\Ns)$ cost for every time step. However, regardless of $d$, the conventional time-marching algorithm requires $O(\Ns^2\Nt)$ complexity to evaluate the RHS of (\ref{eq:linear}) because the RHS consists of at most $\gamma^*$ ($\sim\Nt^0$) matrix-vector products, and the coefficient matrices $\mat{U}^{(\gamma)}$ and $\mat{W}^{(\gamma)}$ become dense as the time difference $\gamma$ becomes large.

\def\tc#1{[#1]}
\def\Kt{K_{\rm t}}

\section{Interpolation-based FMM}\label{s:fmm2}

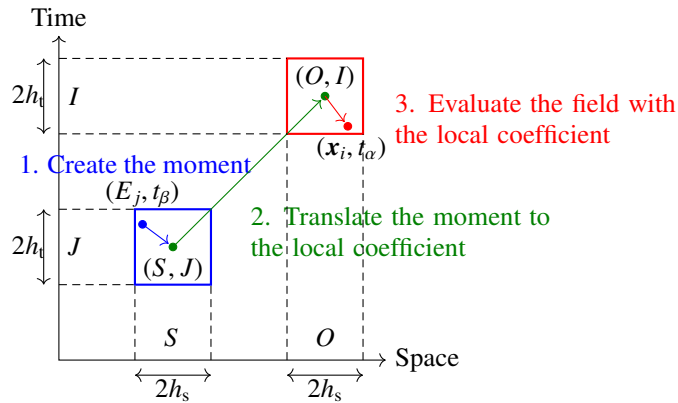
\begin{figure}[H]
  \centering
  \begin{tikzpicture}[scale=1.0]
  \draw [<->] (0.0, 4.3)--(0.0, 0.0)--(4.3, 0.0);
  \draw [anchor=west] (4.3, 0.0) node {Space};
  \draw [anchor=south] (0.0, 4.3) node {Time};

  \def\rectanglepath{-- +(1,0) -- +(1,1) -- +(0,1) -- cycle} 

  \draw [thick, blue] (1,1) \rectanglepath;
  \draw [fill, blue] (1.1,1.8) circle(0.05);
  \draw (1.1,2.2) node {$(E_j,t_\beta)$};

  \def\darkgreen{green!50!black};

  \draw (1.0,2.6) [blue] node {1. Create the moment};
  \draw (2.4,1.7) [anchor=west, text width=4cm, \darkgreen] node {2. Translate the moment to the local coefficient};
  \draw (4.3,3.2) [anchor=west, text width=4cm, red] node {3. Evaluate the field with the local coefficient};

  \draw [fill, \darkgreen] (1.5,1.5) circle(0.05);
  \draw (1.5,1.2) node {$(S,J)$};

  \draw [thick, red] (3,3) \rectanglepath;
  \draw [fill, red] (3.8,3.1) circle(0.05);
  \draw (3.85,2.8) node {$(\bm{x}_i,t_\alpha)$};

  \draw [fill, \darkgreen] (3.5,3.5) circle(0.05);
  \draw (3.5,3.75) node {$(O,I)$};

  \draw [->, shorten >=2pt, blue] (1.1,1.8)--(1.5,1.5);
  \draw [->, shorten >=2pt, \darkgreen] (1.5,1.5)--(3.5,3.5);
  \draw [->, shorten >=2pt, red] (3.5,3.5)--(3.8,3.1);

  \draw [densely dashed] (3,0)--(3,3);
  \draw [densely dashed] (4,0)--(4,3);
  \draw (3.5,0.3) node {$O$};
  \draw [<->] (3.0,-0.2) -- (4.0,-0.2);
  \draw (3.5,-0.4) node {$2\Hs$};

  \draw [densely dashed] (1,0)--(1,1);
  \draw [densely dashed] (2,0)--(2,1);
  \draw (1.5,0.3) node {$S$};
  \draw [<->] (1.0,-0.2) -- (2.0,-0.2);
  \draw (1.5,-0.4) node {$2\Hs$};

  \draw [densely dashed] (0,4)--(3,4);
  \draw [densely dashed] (0,3)--(3,3);
  \draw (0.2,3.5) node {$I$};
  \draw [<->] (-0.2,3) -- (-0.2,4);
  \draw (-0.4,3.5) node {$2\Ht$};

  \draw [densely dashed] (0,2)--(1,2);
  \draw [densely dashed] (0,1)--(1,1);
  \draw (0.2,1.5) node {$J$};
  \draw [<->] (-0.2,1) -- (-0.2,2);
  \draw (-0.4,1.5) node {$2\Ht$};

\end{tikzpicture}
  \caption{Evaluation of the far-field interaction from the source cluster $S\times J$ to the observation cluster $O\times I$ through the moment and local coefficient via the multipole-to-local translation (M2L). Here, the one-dimensional space is considered for ease of explanation.}
  \label{fig:config}
\end{figure}

The present time-domain FMM is a fast approximation method to evaluate the RHS of (\ref{eq:linear}), that is, the contribution from the passed time steps to the future time steps. To see the essential formulation of the FMM, we may consider a \textit{far-field interaction} between two \textit{clusters} $(O,I)$ and $(S,J)$ in space-time, where $O$ ($\subset\bbbr^3$) and $S$ ($\subset\bbbr^3$) are called observation and source cubes (called \textit{cells}), respectively, such that they are non-overlapping (well separated) (Figure~\ref{fig:config}). The side length of $O$ and $S$ is denoted as $2\Hs$. Further, $I$ ($\subset\bbbr$) and $J$ ($\subset\bbbr$) are the observation and source \textit{time intervals} such that $I\ge J$, which means that $t\ge s$ for any $t\in I$ and $s\in J$. The length (duration) of $I$ and $J$ is denoted as $2\Ht$. Then, the far-field interaction can be expressed as follows:
\begin{eqnarray*}
  \Psi(\bm{x}_i,t_\alpha):=\sum_{\{\beta\;\mid\;t_\beta\in J\}}\sum_{\{j\;\mid\;E_j\subset S\}} \left( U_{ij}^{(\alpha-\beta+1)}\tilde{\tau}_j^{\beta-1} - W_{ij}^{(\alpha-\beta+1)}\tilde{\sigma}_j^{\beta-1} \right)
  \quad\text{for $\bm{x}_i\in O$ and $t_\alpha\in I$}.
\end{eqnarray*}
Although we can analytically evaluate the spatial integrals in $U_{ij}^{(\alpha-\beta+1)}$ and $W_{ij}^{(\alpha-\beta+1)}$, as shown in \ref{s:element_integral}, the resulting expressions are too complicated to apply the interpolation-based FMM to those expressions. Hence, we consider the far-field interaction in the following form:
\begin{eqnarray}
  \Psi(\bm{x}_i,t_\alpha):=\frac{1}{4\pi(c\Dt)^d}\sum_{\{\beta\;\mid\;t_\beta\in J}\sum_{\{j\;\mid\;E_j\subset S\}}\int_{E_j} \left( U(\bm{x}_i,\bm{y},t_\alpha,t_{\beta-1})\tau_j^{\beta-1} - \partial_n^y U(\bm{x}_i,\bm{y},t_\alpha,t_{\beta-1})\sigma_j^{\beta-1} \right)\diff S_y
  \quad\text{for $\bm{x}_i\in O$ and $t_\alpha\in I$}.
  \label{eq:target}
\end{eqnarray}
Here, we recall that the single-layer kernel $U$ was defined in (\ref{eq:U}) and note that the double-layer kernel can be expressed as
\begin{eqnarray*}
 \partial_n^y U=\bm{n}_y\cdot\nabla_y U=\bm{n}_y\cdot\bm{W},
\end{eqnarray*}
where $\bm{W}$ was defined in (\ref{eq:W}).

To realise the separation of variables, which is the key to constructing an FMM, we interpolate the function $U(\bm{x},\bm{y},t,s)$ in terms of all the eight variables, i.e. $x_1$, $x_2$, $x_3$, $y_1$, $y_2$, $y_3$, $t$, and $s$. For example, the interpolation in terms of $x_1$ can be written as
\begin{eqnarray*}
  U(\bm{x},\bm{y},t,s)\approx\sum_{a_1=0}^{\Ps-1}U(\bar{O}+\Hs\bm{\omega}_{a_1}^\Ps,\bm{y},t,s)\ell_{a_1}\left(\frac{\bm{x}-\bar{O}}{\Hs}\right),
\end{eqnarray*}
where $\ell_{a_1}$ denotes the prescribed interpolation function; $\Ps$ denotes the number of spatial interpolation points, which are specified with the normalised nodes $\omega_{a_1}^\Ps\in[-1,1]$; and $\bar{O}$ denotes the centre of $O$. Similarly, we can interpolate $U$ for all the variables as follows:
\begin{eqnarray}
  U(\bm{x},\bm{y},t,s)
  \approx\sum_{a<\Ps}\sum_{b<\Ps}\sum_{m<\Pt}\sum_{n<\Pt}U_{a,b,m,n}(O,S,I,J)\ell_a\left(\frac{\bm{x}-\bar{O}}{\Hs}\right)\ell_m\left(\frac{t-\bar{I}}{\Ht}\right)\ell_b\left(\frac{\bm{y}-\bar{S}}{\Hs}\right)\ell_n\left(\frac{s-\bar{J}}{\Ht}\right),
  \label{eq:U_intp}
\end{eqnarray}
where we defined
\begin{eqnarray}
  U_{a,b,m,n}(O, S, I, J):=U\left(\bar{O}+\Hs\bm{\omega}_a^\Ps, \bar{S}+\Hs\bm{\omega}_b^\Ps, \bar{I}+\Ht\omega_m^\Pt, \bar{J}+\Ht\omega_n^\Pt\right)
  \label{eq:Uabmn}
\end{eqnarray}
and introduced the following notations regarding the spatial variables:
\begin{eqnarray*}
  \sum_{\nu<\Ps}:=\sum_{\nu_1=0}^{\Ps-1}\sum_{\nu_2=0}^{\Ps-1}\sum_{\nu_3=0}^{\Ps-1}, \quad
  \ell_\nu(\bm{\xi}):=\ell_{\nu_1}(\xi_1)\ell_{\nu_2}(\xi_2)\ell_{\nu_3}(\xi_3),\quad
  \bm{\omega}_\nu^\Ps:=(\omega_{\nu_1}^\Ps, \omega_{\nu_2}^\Ps, \omega_{\nu_3}^\Ps)\quad\text{for $\nu=a$, $b$}.
\end{eqnarray*}
Here, we define $\bm{\xi}:=(\xi_1,\xi_2,\xi_3)$ such that $-1\le\xi_1,\xi_2,\xi_3\le 1$. Following the previous study~\cite{takahashi2014}, we use the cubic Hermite interpolation such that the first derivative of the interpolated function is approximated with finite differences.

By plugging (\ref{eq:U_intp}) into the layer potential in (\ref{eq:target}), we have
\begin{eqnarray*}
  \Psi(\bm{x}_i,t_\alpha)
  &=&\sum_{a<\Ps}\sum_{m<\Pt}\ell_a\left(\frac{\bm{x}_i-\bar{O}}{\Hs}\right)\ell_m\left(\frac{t_\alpha-\bar{I}}{\Ht}\right)
   \sum_{b<\Ps}\sum_{n<\Pt}U_{a,b,m,n}(O,S,I,J)\nonumber\\
   &&\times\underbrace{\frac{1}{4\pi(c\Dt)^d}\sum_{\beta-1}\sum_j \int_{E_j}\left(\ell_b\left(\frac{\bm{y}-\bar{S}}{\Hs}\right)\tau_j^{\beta-1}-\partial_n^y\ell_b\left(\frac{\bm{y}-\bar{S}}{\Hs}\right)\sigma_j^{\beta-1}\right)\diff S_y\ell_n\left(\frac{t_{\beta-1}-\bar{J}}{\Ht}\right)}_{\displaystyle M_{b,n}(S,J)},
\end{eqnarray*}
where $M_{b,m}(S,J)$ denotes the multipole moment, which contains the information of the source cluster $(S,J)$. By performing the summations over $b$ and $n$, we can obtain the following expression:
\begin{eqnarray}
  \Psi(\bm{x}_i,t_\alpha)
  &=&\sum_{a<\Ps}\sum_{m<\Pt}\ell_a\left(\frac{\bm{x}_i-\bar{O}}{\Hs}\right)\ell_m\left(\frac{t_\alpha-\bar{I}}{\Ht}\right)\underbrace{\sum_{b<\Ps}\sum_{n<\Pt}U_{a,b,m,n}(O,S,I,J) M_{b,m}(S,J)}_{\displaystyle L_{a,m}(O,I)}\nonumber\\
  &=&\sum_{a<\Ps}\sum_{m<\Pt}\ell_a\left(\frac{\bm{x}_i-\bar{O}}{\Hs}\right)\ell_m\left(\frac{t_\alpha-\bar{I}}{\Ht}\right) L_{a,m}(O,I).
     \label{eq:target_m2l}
\end{eqnarray}
Here, $L_{a,m}(O,I)$ denotes the local coefficient and gives the M2L translation from $(S,J)$ to $(O,I)$; that is,
\begin{eqnarray}
  L_{a,m}(O,I) = \sum_{b<\Ps}\sum_{n<\Pt}U_{a,b,m,n}(O,S,I,J) M_{b,m}(S,J)
  \quad\Longleftrightarrow\quad
  \matL(O,I) = \matU(O,S,I,J) \mat{M}(S,J),
  \label{eq:m2l}
\end{eqnarray}
where the latter expression is the matrix form of the former one. That is, $\matL$ and $\mat{M}$ are $\Ps^3\Pt$-dimensional vectors, and $\matU$ is a $\Ps^3\Pt$-dimensional square matrix; to construct them, we may define a row (respectively, column) index as $a_1+\Ps(a_2+\Ps(a_3+\Ps m)$ (respectively, $b_1+\Ps(b_2+\Ps(b_3+\Ps n)$), for example.

In the case of the BMBIE in (\ref{eq:bmbie}), we need to apply the differential operator $\left(\frac{\partial}{\partial n_x}-\frac{1}{c}\frac{\partial}{\partial t}\right)$ to $\Psi$ in (\ref{eq:target}) or (\ref{eq:target_m2l}). Because the operator is related to the point $\bm{x}\in O$ and time $t\in I$, we may apply the operator to the product $\ell_a\ell_m$ in (\ref{eq:target_m2l}). Therefore, the FMM for the BMBIE is different from that for the OBIE only in the last stage of the FMM; i.e. the evaluation with the local coefficients at leaf cells, and all the other FMM's operations, including the M2L, are common to both BIEs.

As seen above, the order $d$ of the temporal basis appears only in the function $U$ in (\ref{eq:U}), and thus, the overall algorithm of the present FMM for $d\ge 2$ is basically the same as that for $d=1$ investigated by the previous study~\cite{takahashi2014}. However, to preserve the computational complexity of $O(\Ns^{1+\delta}\Nt)$ even for $d\ge 2$, we need to modify the M2L appropriately. We will mention the details of the M2L in the next section.

\section{M2L}\label{s:m2l}

We propose an algorithm for the M2L that works for any $d$ ($\ge 1$) with a computational complexity of $O(\Nt)$. In Section~\ref{s:m2l_naive}, we point out that the M2L in (\ref{eq:m2l}) requires $O(\Nt^2)$ computational complexity if it is performed directly. To reduce the computational complexity, we introduce the near- and distant-future M2L in Section~\ref{s:m2l_near_distant}, but the latter is still computationally expensive. To achieve $O(\Nt)$ complexity, we reduce the distant-future M2L to recurrence form with the help of the Taylor expansion (actually, binomial expansion) in Section~\ref{s:m2l_efficient}.

\subsection{Naive approach}\label{s:m2l_naive}

To formulate the M2L for any order $d$, let us consider that the time axis is segmented by $\Kt$ time intervals so that each of interval has length $2\Ht$; then, $\Kt\equiv\Nt\Dt/(2\Ht)$ holds. Under this actual setting, the M2L for an observation cluster $(O,I_k)$, where we let $I_k$ be the current time interval, needs to consider \textit{all} the passed time intervals $I_0,\ldots,I_{k-1}$ (see Remark~\ref{remark:passed} below) as well as all the source cells $S$ in the interaction list of $O$, denoted by $\IL(O)$. Therefore, the M2L can be expressed as follows\footnote{It is unnecessary to compute the local coefficient for the last time interval $I_{\Kt-1}$ because the local coefficient for a time interval is cast to the next time interval or more.}:
\begin{eqnarray}
  \mat{L}(O,I_k) = \sum_{l=0}^{k-1} \sum_{S\in\IL(O)} \matU(O,S,I_k,I_l) \mat{M}(S,I_l)
  \quad\Longleftrightarrow\quad
  \mat{L}_k = \sum_{l=0}^{k-1} \matU_{k,l}\mat{M}_l
  \quad\text{for $k=0,\ldots,\Kt-2$}.
  \label{eq:m2l_actual}
\end{eqnarray}
In the latter expression, we omitted the dependency on $O$ and $S$ for brevity. The subscripts $k$ and $l$ denote the indices of the current time interval $I_k$ and source time interval $I_l$, respectively.

\begin{remark}\label{remark:passed}
  Contrary to Remark~\ref{remark:gamma_bound} for the conventional time-marching algorithm (which is applied to the near-field interaction in the FMM), we cannot ignore any contributions from passed time steps (thus, time intervals) in the M2L (thus, the far-field interaction calculation of the FMM). This is because the FMM handles multiple time steps collectively, as a time interval, and thus, a single time interval can be related to the coefficient matrices $\mat{U}^{(\gamma)}$ and $\mat{W}^{(\gamma)}$ such that both $\gamma\le\gamma^*$ and $\gamma>\gamma^*$. This mixing state makes it difficult to use the property $\mat{U}^{(\gamma)}=\mat{W}^{(\gamma)}=\mat{0}$ for $\gamma>\gamma^*$ mentioned in Remark~\ref{remark:gamma_bound}.\footnote{If we do not apply the M2L to the previous time interval $I_{k-1}$ (then the upper bound of the summation over $l$ becomes $k-2$ in (\ref{eq:m2l_actual})), we can avoid the mixing state. However, this increases the number of source cells $S(\subset\IL(O))$ that directly interact with the observation cell $O$. As a result, the computational cost becomes $O(\Ns^2)$, which is obviously undesirable.} 
\end{remark}

In the actual algorithm, a moment $\mat{M}_k$ of the current time interval $I_k$ is cast to the local coefficient of the \textit{future} time intervals $I_{k+1},I_{k+2},\ldots,I_{\Kt-2}$ once the moment is computed at the end of $I_k$. In this regard, we can reformulate the M2L in (\ref{eq:m2l_actual}) so that the local coefficient of each future time interval is updated sequentially according to the following formula:
\begin{eqnarray}
  \mat{L}_l = \mat{L}^{\tc{k-1}}_l + \mat{U}_{l,k}\mat{M}_k^{\tc{k}}\quad\text{for $l=k+1,\ldots,\Kt-2$},
  \label{eq:m2l_naive}
\end{eqnarray}
where, to clarify the order of computations, a number in square brackets denotes the time interval where the quantity with the number is computed. In addition, all the local coefficients are initially zero; we let $\mat{L}_0^{\tc{-1}}=\ldots=\mat{L}_{\Kt-1}^{\tc{-1}}=\mat{0}$. Note that we do not specify $\tc{\cdot}$ to $\mat{U}_{l,k}$ because the matrix is precomputed.

\begin{figure}[H]
  \centering
  \includegraphics[width=.8\textwidth]{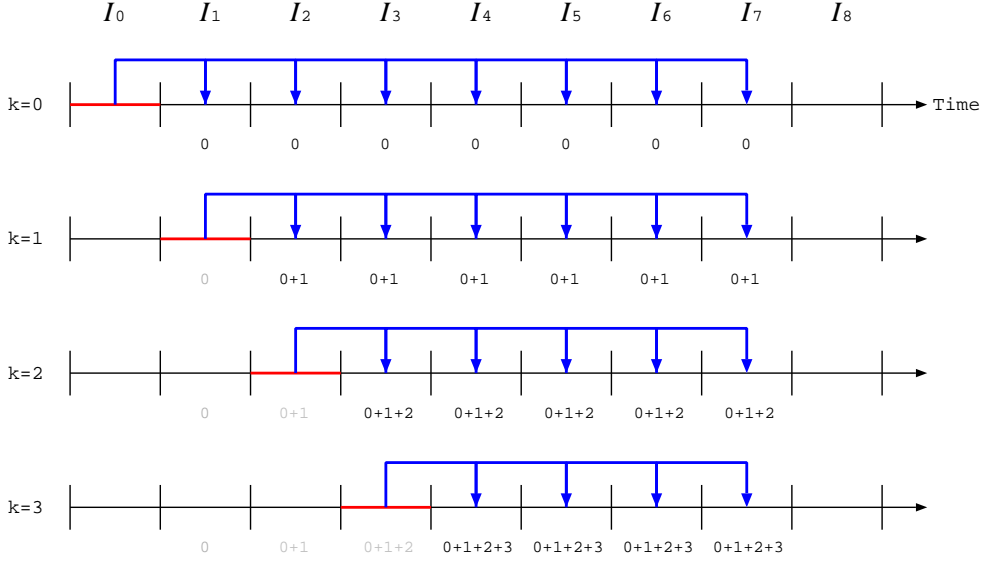}
  \caption{M2L according to (\ref{eq:m2l_naive}), i.e. the naive approach. We consider nine time intervals, i.e. $I_0,\ldots,I_8$. The current time interval indexed by $k$ is coloured in red. The blue arrows from the current time interval $I_k$ to a future time interval $I_l$ (where $l=k+1,\ldots,7$) denote accumulation of $\mat{U}_{l,k}\mat{M}_k$ to the local coefficient of $I_l$, i.e. $\mat{L}_l$. The numbers in a time interval denote the indices of the time intervals that have been accounted in the previous time intervals, i.e. $I_0,\ldots,I_{k-1}$. For example, ``0+1+2'' of $I_4$ at $k=2$ means that the contributions from $I_0$, $I_1$, and $I_2$ are considered in $\mat{L}_4^{\tc{2}}$ at the end of $I_2$. The spatial translation (from a source cell $S$ to the underlying observation cell $O$) is not illustrated here but it is considered in practice.}
  \label{fig:m2l-naive}
\end{figure}

Figure~\ref{fig:m2l-naive} illustrates how the local coefficients are updated according to (\ref{eq:m2l_naive}). Clearly, the computational cost in this naive approach is $O(\Nt^2)$, which is higher than that of the conventional algorithm, and thus unacceptable.

\subsection{Near- and distant-future M2L}\label{s:m2l_near_distant}

To realise the M2L with $O(\Nt)$ complexity, we consider (i) the Taylor expansion of the local coefficients and (ii) rewriting the expanded local coefficients in a recurrence form. These were considered for $d=1$ in the previous study. To generalise the case of $d=1$ to $d\ge 2$, we first split the set of the future time intervals $\{I_k,\ldots,I_{\Kt-2}\}$ into the \textit{near-future time intervals} $\{I_k,\ldots,I_{k+\mu+1}\}$ and \textit{far-future time intervals} $\{I_{k+\mu+2},\ldots,I_{\Kt-2}\}$. Here, the number $\mu$ is determined so that the truncated power function $(c(t-s)-r)_+^d$ of $U$ in (\ref{eq:U}), where $r:=\abs{\bm{x}-\bm{y}}$, can be regarded as an ordinary power function $(c(t-s)-r)^d$, which is infinitely differentiable with respect to the time $t$. Because the exponent $d$ of both power functions does not matter for determining $\mu$, we may follow the case of $d=1$ investigated in the previous study~\cite{takahashi2014} and can state that $\mu\ge 8$ is necessary in conjunction with the construction of the space-time hierarchy. In fact, we will use the lower bound, that is, $\mu=8$ (in all the levels of the space-time hierarchy).

\begin{figure}[H]
  \centering
  \includegraphics[width=.8\textwidth]{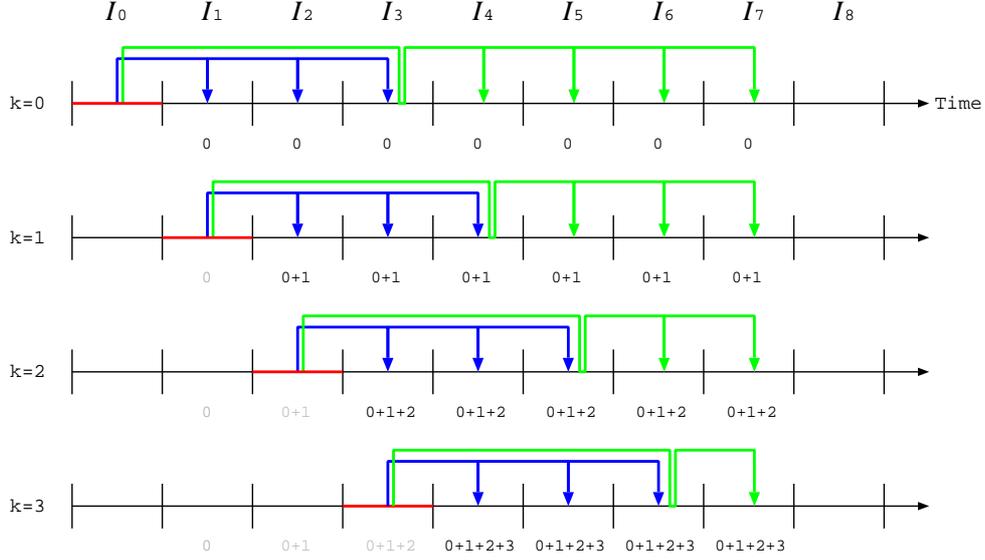}
  \caption{M2L according to both the near-future M2L in (\ref{eq:m2l_near}) and the distant-future M2L in (\ref{eq:m2l_expand}). The current time interval $I_k$ is coloured red. The arrows in blue correspond to the near-future M2L from $I_k$ to a near-future time interval $I_l$, where $l=k+1,\ldots,k+\mu+1$. Here, $\mu$ is set to $2$ (instead of $8$) for ease of explanation. On the other hand, the arrows in green correspond to the distant-future M2L from $I_k$ to all the distant-future time intervals $I_l$ (where $l=k+\mu+2,\ldots,7$) via $I_{k+\mu+1}$.}
  \label{fig:m2l_near}
\end{figure}

We then perform the M2L for the near-future intervals according to (\ref{eq:m2l_naive}), that is,
\begin{eqnarray}
  \mat{L}_{l}^{\tc{k}}=\mat{L}_{l}^{\tc{k-1}}+\mat{U}_{l,k}\mat{M}_k^{\tc{k}}\quad\text{\rm for $l=k+1,\ldots,k+\mu+1$}.
  \label{eq:m2l_near}
\end{eqnarray}
We call this the \textit{near-future M2L}, which is represented by the blue arrows in Figure~\ref{fig:m2l_near}. The computational complexity of the near-future M2L is $O(\mu\Nt)=O(\Nt)$.

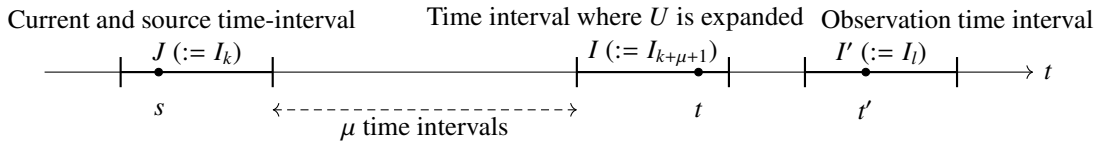
\begin{figure}[H]
  \centering
  \begin{tikzpicture}[scale=1]
  
  \draw [->] (-1, 0)--(12, 0);
  \draw (12.2, 0.0) node {$t$};      
  
  \draw (0.5, 0.0) circle(0.05) [fill];
  \draw (0.5, -0.5) node {$s$};

  \draw (0.0, -0.2)--(0.0, 0.2) [thick];
  \draw (2.0, -0.2)--(2.0, 0.2) [thick];
  \draw (0.0, 0.0)--(2.0, 0.0) [thick];
  \draw (1.0, 0.25) node {$J$ ($:=I_k$)};
  \draw (0.8, 0.7) node {Current and source time-interval};

  \draw (6.5, 0.7) node {Time interval where $U$ is expanded};

  \draw (6.0, -0.2)--(6.0, 0.2) [thick];
  \draw (8.0, -0.2)--(8.0, 0.2) [thick];
  \draw (6.0, 0.0)--(8.0, 0.0) [thick];
  \draw (7.0, 0.25) node {$I$ ($:=I_{k+\mu+1}$)};
  \draw (7.6, 0.0) circle(0.05) [fill];
  \draw (7.6, -0.5) node {$t$};

  \draw (9.0, -0.2)--(9.0, 0.2) [thick];
  \draw (11.0, -0.2)--(11.0, 0.2) [thick];
  \draw (9.0, 0.0)--(11.0, 0.0) [thick];
  \draw (11.0, 0.7) node {Observation time interval};
  \draw (10.0, 0.25) node {$I'$ ($:=I_l$)};
  \draw (9.8, 0.0) circle(0.05) [fill];
  \draw (9.8, -0.5) node {$t'$};


  \draw [<->, densely dashed] (2.0, -0.5)--(6.0, -0.5);
  \draw (4.0, -0.75) node {$\mu$ time intervals};
\end{tikzpicture}
  \caption{Notations related to the expansion of $U$ in (\ref{eq:U_expand}).}
  \label{fig:U_expand}
\end{figure}

The local coefficients of the distant-future time intervals are computed by the Taylor expansion of $\mat{L}_{k+\mu+1}$, which is possible because $\mu$ is chosen so that $\mat{L}_{k+\mu+1}$ is differentiable with respect to time.  To formulate the expansion, we first consider the Taylor expansion of the function $U$ in (\ref{eq:U}). To simplify the notation, we let $J:=I_k$ (i.e. the current and thus source time interval), $I:=I_{k+\mu+1}$ (i.e. the time interval where $U$ is expanded), and $I':=I_l$ (i.e. the observation time interval where $U$ is evaluated; $l\ge k+\mu+2$); see Figure~\ref{fig:U_expand}. Further, we let $s\in J$, $t\in I$ and $t'\in I'$. Then, because $I$ is chosen so that $c(t'-s)-r>0$ holds for any $\bm{x}\in O$ and $\bm{y}\in S\subset\IL(O)$ (where $r:=\abs{\bm{x}-\bm{y}}$), the truncated power function $(c(t'-s)-r)_+^d$ can be treated as the $d$th-order polynomial, i.e. $(c(t'-s)-r)^d$. Therefore, we can expand $U$ with respect to $t$ as follows:
\begin{eqnarray}
  U(\bm{x},\bm{y},t',s)
  &=&\frac{\left(c(t-s)-r+c(t'-t)\right)^d}{r}\nonumber\\
  &=&\sum_{p=0}^\infty \frac{1}{p!}\frac{\partial^p U(\bm{x},\bm{y},t,s)}{\partial t^p}(t'-t)^p\nonumber\\
  &=&\sum_{p=0}^d \underbrace{\frac{1}{c^p p!}\frac{\partial^p U(\bm{x},\bm{y},t,s)}{\partial t^p}}_{\displaystyle U^{(p)}(\bm{x},\bm{y},t,s)}\left(c(t'-t)\right)^p.
    \label{eq:U_expand}
\end{eqnarray}
We simply call $U^{(p)}$ in (\ref{eq:U_expand}) the $p$th \textit{derivative} of $U$, although it differs from the exact temporal derivative $\frac{\partial^p U}{\partial t^p}$ by the factor $c^p p!$. We note that the above Taylor expansion can be also obtained by the binomial expansion of the numerator of $U$.

Using the Taylor expansion of $U$ in (\ref{eq:U_expand}), we can write the M2L from $(S,J)$ to $(O,I')$ in the first expression of (\ref{eq:m2l}) as follows:
\begin{eqnarray*}
  L_{a,m}(O,I')
  &=&\sum_{b<\Ps}\sum_{n<\Pt}U_{a,b,m,n}(O,S,I',J) M_{b,n}(S,J)\nonumber\\
  &=&\sum_{b<\Ps}\sum_{n<\Pt}U\left(\bar{O}+h_s\bm{\omega}_a^\Ps,\bar{S}+h_s\bm{\omega}_b^\Ps,\bar{I}'+h_t\omega_m^\Pt,\bar{J}+h_t\omega_n^\Pt\right) M_{b,n}(S,J)\qquad\text{\hfill $\because$ (\ref{eq:Uabmn})}\nonumber\\
  &=&\sum_{b<\Ps}\sum_{n<\Pt}U\left(\bar{O}+h_s\bm{\omega}_a^\Ps,\bar{S}+h_s\bm{\omega}_b^\Ps,\bar{I}'-\bar{I}+\bar{I}+h_t\omega_m^\Pt,\bar{J}+h_t\omega_n^\Pt\right) M_{b,n}(S,J)\qquad\text{\hfill $\because$ Add and subtract $\bar{I}$}\nonumber\\
  &=&\sum_{b<\Ps}\sum_{n<\Pt}\sum_{p=0}^d U^{(p)}\left(\bar{O}+h_s\bm{\omega}_a^\Ps,\bar{S}+h_s\bm{\omega}_b^\Ps,\bar{I}+h_t\omega_m^\Pt,\bar{J}+h_t\omega_n^\Pt\right) \left(c(\bar{I}'-\bar{I})\right)^p M_{b,n}(S,J)\qquad\text{\hfill $\because$ (\ref{eq:U_expand})}\nonumber\\
  &=&\sum_{b<\Ps}\sum_{n<\Pt}\sum_{p=0}^d U^{(p)}_{a,b,m,n}\left(O,S,I,J\right) \left(c(\bar{I}'-\bar{I})\right)^p M_{b,n}(S,J).\qquad\text{\hfill $\because$ (\ref{eq:Uabmn})}
\end{eqnarray*}  
This can be re-expressed in a form that is similar to (\ref{eq:m2l_near}) as follows:
\begin{eqnarray}  
  \mat{L}_{l}^{\tc{k}} = \underbrace{\mat{L}_{l}^{\tc{k-1}}}_{\displaystyle\mat{0}} + \sum_{p=0}^d \underbrace{\mat{U}_{k+\mu+1,k}^{(p)}\mat{M}_{k}^{\tc{k}}}_{\displaystyle \mat{L}_{k+\mu+1,k}^{(p)\tc{k}}}(l-k-\mu-1)^pT^p = \sum_{p=0}^d \mat{L}_{k+\mu+1,k}^{(p)\tc{k}}(l-k-\mu-1)^pT^p\quad\text{for $l\ge k+\mu+2$},
  \label{eq:m2l_expand}
\end{eqnarray}  
where $\mat{L}_{l}^{\tc{k-1}}\equiv\mat{0}$ was used because $I_k$ is the first time interval where we substitute a (non-zero) value into $\mat{L}_{l}$. In addition, we used the following notations:
\begin{itemize}
\item Matrices $\mat{U}_{l,k}^{(p)}$: Matrix representation of the $p$th-order derivative $U_{a,b,m,n}^{(p)}(O,S,I_l,I_k)$, where $p=1,\ldots,d$ and $\mat{U}_{l,k}^{(0)}\equiv\mat{U}_{l,k}$. These matrices are precomputed and are thus free from the ordering index $\tc{\cdot}$; see Remark~\ref{remark:precomp} below.
\item Vectors $\mat{L}_{l,k}^{(p)}$: Contribution to the derivatives of the local coefficient $\mat{L}_l^{(p)}$ from the moment $\mat{M}_k$ of the current time interval $I_k$, i.e.
  \begin{eqnarray}
    \mat{L}_{l,k}^{(p)}:=\mat{U}_{l,k}^{(p)}\mat{M}_k,
    \label{eq:L_l,k}
  \end{eqnarray}
  where $p=1,\ldots,d$ and $\mat{L}_{l,k}^{(0)}\equiv \mat{L}_{l,k}$.
\item Scalar constant $T:=c(2\Ht)$: Length of the time interval multiplied by the wave speed $c$.
\end{itemize}
We call the formula in (\ref{eq:m2l_expand}) the \textit{distant-future M2L} from $I_k$ to $I_l$ (via $I_{k+\mu+1}$), which is represented by the green arrows in Figure~\ref{fig:m2l_near}. As we can see in this figure, the total computational cost is still $O(\Nt^2)$ because the number of the distant-future intervals is $O(\Kt)=O(\Nt)$ for every time interval $I_k$.

\begin{remark}\label{remark:precomp}
  The matrices $\mat{U}_{l,k} (\equiv\mat{U}_{l,k}^{(0)})$ and $\mat{U}_{l,k}^{(p)}$  can be precomputed because they depend on the difference of the time intervals $I_l$ and $I_k$ as well as on the relative position of cells $O$ and $S$. We may precompute $\mu+1$ matrices for each $p\in[1,d]$ and the possible 316 pairs of $O$ and $S$.
\end{remark}

\subsection{M2L with $O(\Nt)$ complexity}\label{s:m2l_efficient}

\begin{figure}[H]
  \centering
  \includegraphics[width=.8\textwidth]{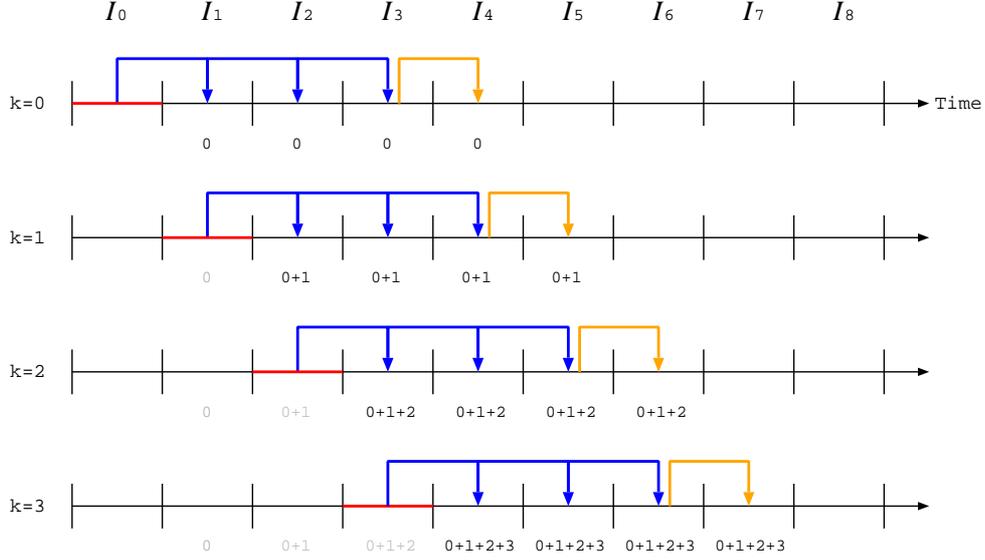}
  \caption{M2L with $O(\Nt)$ complexity given by Formula~\ref{formula:m2l}. The blue arrows correspond to the near-future M2Ls, whereas the orange arrows correspond to the distant-future M2Ls in the recurrence form in (\ref{eq:m2l_rec}).}
  \label{fig:m2l-rec}
\end{figure}

To achieve $O(\Nt)$ complexity, we derive a recurrence formula of the local coefficient $\mat{L}_{k}$ with respect to the index $k$. This is represented by the orange arrows in Figure~\ref{fig:m2l-rec}, while the blue arrows correspond to the near-future M2L mentioned above. Because we consider $\mu+2$ time intervals for each current time interval $I_k$, the computational complexity is indeed $O((\mu+2)\Kt)=O(\Nt)$.

By using both the near-future M2L in (\ref{eq:m2l_near}) and the distant-future M2L in (\ref{eq:m2l_expand}), we can construct the recurrence formula of the M2L with $O(\Nt)$ complexity as in the following Formula~\ref{formula:m2l}.

\begin{formula}[M2L with $O(\Nt)$ complexity]\label{formula:m2l}
  Let $I_k$ be the current time interval, where $k=0,1,\ldots,\Kt-2$. First, we compute $\mat{L}_{k+1},\ldots,\mat{L}_{k+\mu+1}$ according to the near-future M2L in (\ref{eq:m2l_near}), after which we can recursively compute the local coefficient of the {\rm latest} time interval, i.e. $\mat{L}_{k+\mu+2}$, by
  \begin{eqnarray}
    \mat{L}_{k+\mu+2}^{\tc{k}}=\mat{L}_{k+\mu+1}^{\tc{k}}+\sum_{p=1}^d \mat{L}_{k+\mu+1}^{(p)\tc{k}} T^p,
    \label{eq:m2l_rec}
  \end{eqnarray}
  where $\bm{L}_{k+\mu+1}^{(p)}$ corresponds to the $p$th-order {\rm derivative of local-coefficient} $\mat{L}_{k+\mu+1}$ (where $p=1,\ldots,d$) and is defined as
  \begin{eqnarray}
    \mat{L}_{k+\mu+1}^{(p)} := \sum_{l=0}^k \mat{L}_{l+\mu+1,l}^{(p)}\left((k+1-l)^p-(k-l)^p\right).
    \label{eq:L_derivative}
  \end{eqnarray}
  This derivative can be computed by the following recurrence formula:
  \begin{eqnarray}
    \mat{L}_{k+\mu+1}^{(p)\tc{k}} = \mat{L}^{(p)\tc{k-1}}_{(k-1)+\mu+1} + \mat{L}_{k+\mu+1,k}^{(p)\tc{k}}+\sum_{m=0}^{p-2}c_m^p(k) \mat{L}_{(k-1)+\mu+1}^{(p,m)\tc{k-1}},
    \label{eq:m2l_expand_rec1}
  \end{eqnarray}
  where $L_{k+\mu+1}^{(p,m)}$ is called the {\rm auxiliary local-coefficient} defined by
  \begin{eqnarray}
    \mat{L}_{k+\mu+1}^{(p,m)}:=\sum_{l=0}^{k} \mat{L}_{l+\mu+1,l}^{(p)}l^m.
    \label{eq:L_aux}
  \end{eqnarray}
  The auxiliary local coefficient can be computed by the following recurrence formula:
  \begin{eqnarray}
    \mat{L}_{k+\mu+1}^{(p,m)\tc{k}}=\mat{L}_{(k-1)+\mu+1}^{(p,m)\tc{k-1}} + \mat{L}_{k+\mu+1,k}^{(p)\tc{k}} k^m.
    \label{eq:m2l_expand_rec2}
  \end{eqnarray}
\end{formula}

\begin{proof}
First, we can prove the recurrence formula in (\ref{eq:m2l_rec}) --- not rigorously but deductively --- as shown in \ref{s:m2l_rec}.

Second, we derive the second recurrence formula (with respect to the index $k$) in (\ref{eq:m2l_expand_rec1}). To this end, we derive the recurrence relation of the term $(k+1-l)^p-(k-l)^p$ in (\ref{eq:L_derivative}) as follows:
\begin{eqnarray}
  (k+1-l)^p-(k-l)^p
  &=&(k-l)^p-(k-1-l)^p+\underline{(k-1-l)^p-2(k-l)^p+(k+1-l)^p}\nonumber\\
  &=&(k-l)^p-(k-1-l)^p+\sum_{m=0}^{p-2}c_m^p(k)l^m,
     \label{eq:tmp1}
\end{eqnarray}
where $c_m^p(k)$ is defined as the coefficient of the underlined polynomial in terms of the index $l$ and computed as follows:
\begin{eqnarray}
  c_m^p(k)=\binom{p}{m}\left((k-1)^{p-m}-2k^{p-m}+(k+1)^{p-m}\right)(-1)^m.
  \label{eq:cmp}
\end{eqnarray}
It should be noted that $c_p^p(k)=c_{p-1}^p(k)=0$ for any $k$, and thus, the upper limit of the summation over $m$ is $p-2$ in (\ref{eq:tmp1}). With this expression, we can rewrite (\ref{eq:L_derivative}) as follows:
\begin{eqnarray*}
  \mat{L}^{(p)\tc{k}}_{k+\mu+1}
  &=&\sum_{l=0}^{k-1} \mat{L}_{l+\mu+1,l}^{(p)}\left((k+1-l)^p-(k-l)^p\right)+\mat{L}_{k+\mu+1,k}^{(p)}\quad\text{(The last term $l=k$ was separated.)}\\
  &=&\sum_{l=0}^{k-1} \mat{L}_{l+\mu+1,l}^{(p)}\left((k-l)^p-(k-1-l)^p\right)
     +\sum_{l=0}^{k-1} \mat{L}_{l+\mu+1,l}^{(p)}\sum_{m=0}^{p-2}c_m^p(k)l^m+\mat{L}_{k+\mu+1,k}^{(p)}\quad\text{((\ref{eq:tmp1}) was substituted.)}\\
  &=&\mat{L}^{(p)\tc{k-1}}_{(k-1)+\mu+1}+\sum_{m=0}^{p-2}c_m^p(k)\left(\sum_{l=0}^{k-1} \mat{L}_{l+\mu+1,l}^{(p)}l^m\right)+\mat{L}_{k+\mu+1,k}^{(p)}.\quad\text{((\ref{eq:L_derivative}), where $k$ is replaced with $k-1$, was used.)}
\end{eqnarray*}
Then, by replacing the summation over $l$ with (\ref{eq:L_aux}), where $k$ is replaced with $k-1$, we can obtain the underlying formula in (\ref{eq:m2l_expand_rec1}).

The third recurrence formula for the auxiliary local coefficient in (\ref{eq:m2l_expand_rec2}) follows from its definition in (\ref{eq:L_aux}): we may separate the last term of $l=k$, i.e. $\mat{L}_{k+\mu+1,k}^{(p)}k^m$, from the summation over $l$.

\end{proof}

\begin{remark}[Coefficients $c_m^p(k)$ in (\ref{eq:cmp})]
  The coefficients $c_m^p(k)$ (where $p=1,2,\ldots$ and $0\le m \le p-2$) can be computed as in the table below.
  \begin{center}
    \begin{tabular}{|c|c|c|c|c|}
      \hline
      $p~\backslash~m$ & $0$ & $1$ & $2$ & $3$\\
      \hline
      $1$ & - & - & - & -\\
      \hline
      $2$ & $2$ & - & - & -\\
      \hline
      $3$ & $6k$ & $-6$ & - & -\\
      \hline
      $4$ & $2+12k^2$ & $-24k$ & $12$ & -\\
      \hline
      $5$ & $10k+20k^3$ & $-10-60k^2$ & $60k$ & $-20$\\
      \hline
    \end{tabular}
  \end{center}

\end{remark}

\begin{remark}\label{remark:aux}
The number of auxiliary local coefficients $\mat{L}_k^{(p,m)}$ is $n_d:=\sum_{p=2}^d(p-1)=\frac{d(d-1)}{2}$. For example, $n_1=0$, $n_2=1$, $n_3=3$, $n_4=6$, and $n_5=10$.
\end{remark}

\begin{remark}\label{remark:d=1}
  In the previous study~\cite{takahashi2014} based on $d=1$, the first-order derivative of local coefficients $\mat{L}^{(1)}_k$ (denoted by $\dot{L}_k$ in \cite{takahashi2014}) is actually used, but the auxiliary local coefficients $\mat{L}^{(p,m)}_k$ never appear. This is consistent with the present result.
\end{remark}

\subsection{Algorithm of M2L}\label{s:m2l_algo}

The algorithm of the M2L in Formula~\ref{formula:m2l} is shown in Algorithm~\ref{algo:m2l}, which is performed at the end of the underlying or current time interval $I_k$, where $k=0,1,\ldots,\Kt-2$, at every cell $O$ in every level of the octree after we have computed the moment of $I_k$, i.e. $\mat{M}_k$, for all the cells of the level. In the algorithm, all the lengths of the loops are independent of $k$ owing to the recurrence formulae in (\ref{eq:m2l_near}), (\ref{eq:m2l_expand_rec1}), and (\ref{eq:m2l_expand_rec2}). To be precise, the length of the loop over $l$ is $\mu+1$ ($=9$), and the lengths of the loops over $p$ and $m$ are up to $d$. In addition, although the loop over the source cells, which is included in the first two loops, is omitted in the algorithm, the loop length is at most 189, as in the ordinary FMMs~\cite{greengard1987}. Therefore, because $k\lesssim \Kt \sim \Nt$, the total computational complexity of the M2L is $\Order(\Nt)$.

\begin{algorithm}[H]
  \caption{M2L for the time interval $I_k$. Some quantities are coloured to clarify the data dependency.}
  \label{algo:m2l}
  \begin{algorithmic}
    \STATE // Initialisation at $k=0$
    \IF{$k==0$}
    
    \STATE $\mat{L}_i=\mat{L}_i^{(p)}=\mat{0}$ for $i\in[0,\mu+1]$ and $p\in[1,d]$
    
    \STATE $\mat{L}_\mu^{(p,m)}=\mat{0}$ for $p\in[2,d]$ and $m\in[0,p-2]$
    
    \ENDIF
    
    \STATE // Perform the near-future M2L from $I_k$ to the near-future time intervals, i.e. $I_{k+1},\ldots,I_{k+\mu+1}$, according to (\ref{eq:m2l_near}).

    \FOR{$l=1$ \TO $\mu+1$}
    \STATE $\red{\mat{L}_{k+l}} = \red{\mat{L}_{k+l}} + \mat{U}_{k+l,k} \magenta{\mat{M}_k}$ // The loop over source cells is omitted.
    \ENDFOR

    \STATE // Compute the derivatives $\darkgreen{\mat{L}_{k+\mu+1,k}^{(p)}}$ in (\ref{eq:L_l,k}), i.e. the contribution to the derivatives of the local coefficients from $I_k$. These are temporarily used in the following computations.
    \FOR{$p=1$ \TO $d$}
    \STATE $\darkgreen{\mat{L}_{k+\mu+1,k}^{(p)}}=\mat{U}_{k+\mu+1,k}^{(p)}\magenta{\mat{M}_k}$ // The loop over source cells is omitted.
    \ENDFOR

    \STATE // Update the derivatives of the local coefficient of $I_{k+\mu+1}$ (and the current cell) according to (\ref{eq:m2l_expand_rec1}).
    \FOR{$p=1$ \TO $d$}
    \STATE $\blue{\mat{L}_{k+\mu+1}^{(p)}}=\blue{\mat{L}_{k+\mu}^{(p)}}+\darkgreen{\mat{L}_{k+\mu+1,k}^{(p)}}+\displaystyle\sum_{m=0}^{p-2}c_m^p(k) \cyan{\mat{L}_{k+\mu}^{(p,m)}}$
    \ENDFOR

    \STATE // Compute the local coefficient of $I_{k+\mu+2}$ (and the current cell) according to (\ref{eq:m2l_rec}).
    \STATE $\red{\mat{L}_{k+\mu+2}}=\red{\mat{L}_{k+\mu+1}}+\displaystyle\sum_{p=1}^d \blue{\mat{L}_{k+\mu+1}^{(p)}} T^p$

    \STATE // Update the auxiliary local coefficients of $I_{k+\mu+2}$ (and the current cell) according to (\ref{eq:m2l_expand_rec2}).
    \FOR{$p=2$ \TO $d$}
    \FOR{$m=0$ \TO $p-2$}
    \STATE $\cyan{\mat{L}_{k+\mu+1}^{(p,m)}} = \cyan{\mat{L}_{k+\mu}^{(p,m)}} + \darkgreen{\mat{L}_{k+\mu+1,k}^{(p)}} k^m$
    \ENDFOR
    \ENDFOR

  \end{algorithmic}
\end{algorithm}

\begin{remark}\label{remark:m2l_complexity}
In the previous study~\cite{takahashi2014}, the complexity of the M2L and ``Another'' M2L, which corresponds to (\ref{eq:L_l,k}) where $p=1$, are theoretically estimated as $Z^{-4/3}\Ps^3\Pt(\log\Ps+\log\Pt)\Ns^{1+\delta}\Nt$ and $Z^{-4/3}\Ps^3(\Pt+\log\Ps)\Ns^{1+\delta}\Nt$, respectively, where $\delta=1/3$ or $1/2$ (recall Section~\ref{s:intro}), and the specified number $Z$ denotes the maximum number of boundary elements per leaf. Here, it is assumed that the 4D-FFT is used to perform the calculation. From Remarks~\ref{remark:aux} and \ref{remark:d=1}, these estimates can be unified as $d^2 Z^{-4/3}\Ps^3\Pt(\log\Ps+\log\Pt)\Ns^{1+\delta}\Nt$ in the case of the present M2L.
\end{remark}

\def\prerr#1{{\nprounddigits{2}\npproductsign{\times}\textrm{\numprint{#1}}}}
\def\prtime#1{{\nprounddigits{0}\textrm{\numprint{#1}}}}
\def\prmem#1{{\nprounddigits{0}\npproductsign{\times}\textrm{\numprint{#1}}}}

\def\prerrx#1{\textcolor{red}{\prerr{#1}}}
\def\prtimex#1{\textcolor{black}{--}} 
\def\prmemx#1{\textcolor{black}{--}} 

\def\prerrG#1{\darkgreen{\prerr{#1}}}
\def\prtimeG#1{\darkgreen{\prtime{#1}}}
\def\prmemG#1{\darkgreen{\prmem{#1}}}

\def\prerrxnan#1{\textcolor{red}{nan}}
\def\prtimexnan#1{\textcolor{black}{\prtimex{#1}}}
\def\prmemxnan#1{\textcolor{black}{\prmemx{#1}}}

\section{Numerical verification}\label{s:num}

We wrote the TDBEM program based on the BMBIE and the $d$th-order B-spline temporal basis, where $d$ is $1$, $2$, and $3$. We numerically verified the proposed TDBEM from several aspects, i.e. the kind of BIE, algorithm, FMM's precision parameter ($\Ps$ and $\Pt$), order $d$, and problem size and type.

\subsection{Example 1: Sphere}\label{s:sphere}

\subsubsection{Problem setting}\label{s:sphere_setting}

\begin{figure}[H]
  \centering
  \includegraphics[width=.4\textwidth]{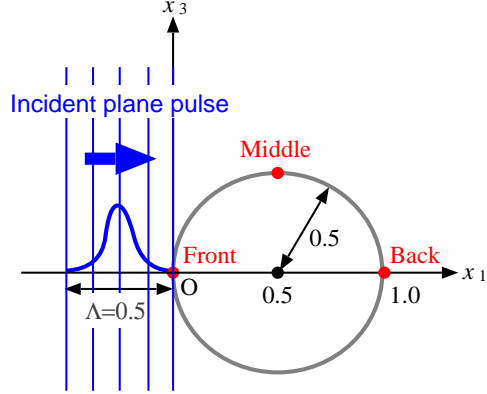}
  \caption{Sphere model. The boundary condition is $u=0$ or $q=0$ on the surface. The profiles of $u$ and $q$ shown in Figures~\ref{fig:profile_dirichlet} and \ref{fig:profile_neumann} were computed at the three red points.}
  \label{fig:sphere}
\end{figure}

As in the previous study~\cite{takahashi2014}, we considered scattering problems involving a spherical scatterer with  radius $0.5$ and centre $(0.5, 0.0, 0.0)$ (Figure~\ref{fig:sphere}). We gave the incidental field $u^{\rm in}$ as a plane pulse that propagates in the $+x_1$-direction, that is,
\begin{eqnarray} 
  u^{\rm in}(\bm{x},t)=0.5\left(1-\mathrm{Cos}\frac{2\pi}{\Lambda}(ct-x_1)\right),
  \label{eq:uin}
\end{eqnarray} 
where $\Lambda$ denotes a pulse length of $0.5$ and
\begin{eqnarray} 
  \mathrm{Cos}(x):=\begin{cases}
    \cos(x) & 0\le x\le 2\pi\\
    1 & \textrm{otherwise}
    \end{cases}.
\end{eqnarray} 
With regard to the boundary condition, we gave either the homogeneous Neumann or Dirichlet boundary condition, i.e. $q=0$ and $u=0$, respectively, on the surface of the sphere.

Both problems were solved by the OBIE and BMBIE. The order $d$ of the B-spline temporal basis was set to $1$, $2$, or $3$. The algorithm used to solve the BIEs was the conventional algorithm (CONV for short) or the present interpolation-based FMM. In the latter case, we let the numbers of the interpolation nodes, i.e. $\Ps$ and $\Pt$, be $8$ or $12$, where the former is expected to be faster but less accurate than the latter. We call the FMM using $\Ps=\Pt=8$ FAST8 and the FMM using $\Ps=\Pt=12$ FAST12. In addition, we generated the octree in such a way that the number of boundary elements in each leaf is $100$ or less.

The surface of the sphere was discretised with $2880$, $5120$, or $11520$ triangular boundary elements; then, the edge length (denoted by $\Ds$) ranged from $0.0461$ to $0.0562$, $0.0350$ to $0.0423$, and $0.0231$ to $0.0283$, respectively. Correspondingly, we let the time-step length $\Dt$ be $0.04$, $0.03$, and $0.02$ so that $\Dt$ can satisfy $\Ds/(c\Dt)\lesssim 1$, which is not a sufficient condition for stabilisation of the solution but is often used in the literature. Correspondingly, we let the number of time steps (i.e. $\Nt$) be $240$, $320$, and $480$, respectively, so that the analysis time, i.e. $\Nt\Dt$, is $9.6$ in all the cases.

To measure the numerical accuracy, we computed the relative $\ell^2$-error of the TDBEM's solution $v_i(t_\alpha)$ (where $v$ is $u$ and $q$ for the Neumann and Dirichlet problems, respectively) to the reference solution, denoted by $\bar{v}(\bm{x},t)$, i.e.
\begin{eqnarray*}
  \textrm{Error}:=\frac{\sqrt{\displaystyle\sum_i\sum_{0\le\alpha<\Nt}\left(v_i(t_\alpha)-\bar{v}(\bm{x}_i,t_\alpha)\right)^2}}{\sqrt{\displaystyle\sum_i\sum_{0\le\alpha<\Nt}\left(\bar{v}(\bm{x}_i,t_\alpha)\right)^2}},
\end{eqnarray*}
where the collocation (evaluation) points $\bm{x}_i$ were chosen so that they were close to $33$ equidistant points on the arc of the semi-circle, such as $x_2=0$ and $x_3\ge 0$ on the sphere. We computed the reference solutions semi-analytically; they were obtained by applying the numerical Laplace transform to the analytical solutions in the frequency domain~\cite{Bowman}.

\subsubsection{Results}\label{s:sphere_results}

\begin{figure}[H]
  \centering
  \includegraphics[height=.295\textheight]{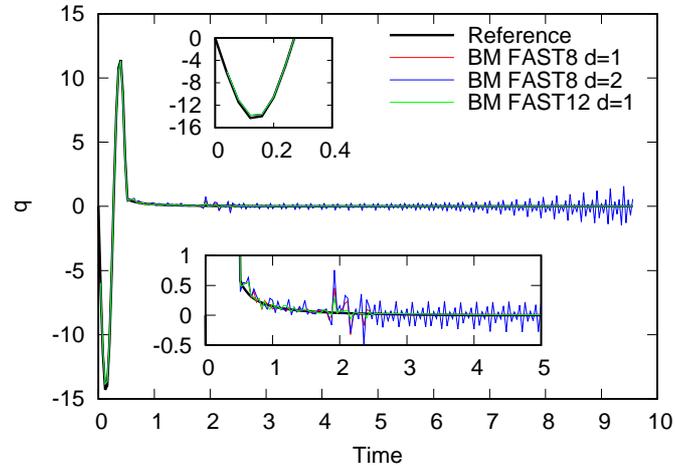}\\
  (a)~Front ($x_1=0.0$)\\
  \includegraphics[height=.295\textheight]{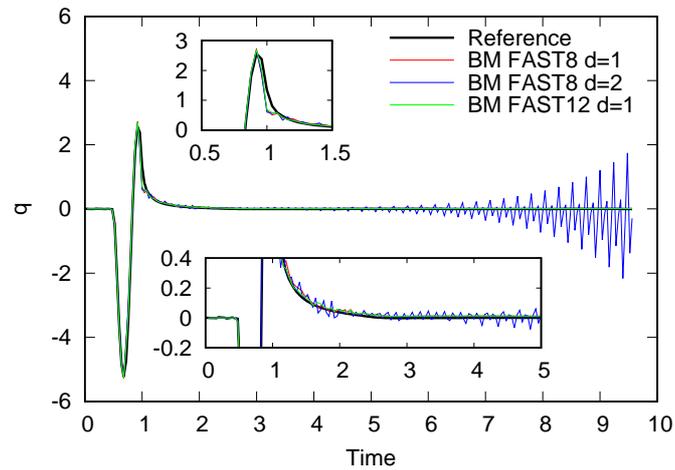}\\
  (b)~Middle ($x_1=0.5$)\\
  \includegraphics[height=.295\textheight]{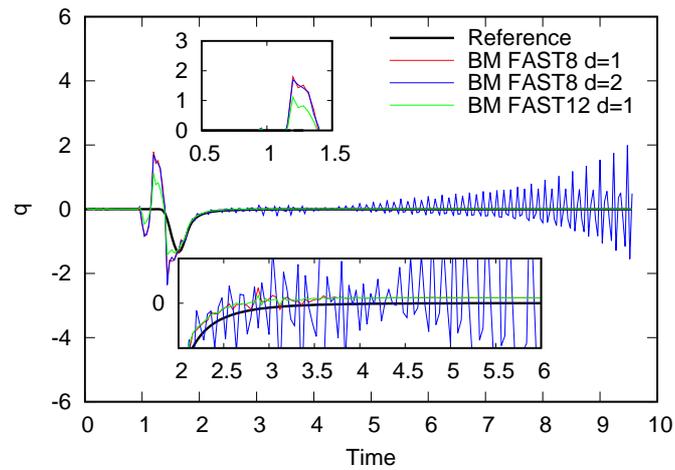}\\
  (c)~Back ($x_1=1.0$)
  \caption{Profiles of $q$ in the homogeneous \textbf{Dirichlet} problem in the case of BM (which stands for BMBIE) and $\Ns=2880$. FAST12 was unstable with $d=2$ as well as $d=1$.}
  \label{fig:profile_dirichlet}
\end{figure}

Tables~\ref{tab:dirichlet} and ~\ref{tab:neumann} show the results of the Dirichlet and Neumann problems, respectively. In these tables, `Error', `Time', and `Mem' denote the relative error mentioned above, the total computation time in seconds, and the memory consumption in GB, respectively. In addition, to visualise the errors of the numerical solutions, we drew the profile of $q$ and $u$ at specific points in Figures~\ref{fig:profile_dirichlet} and \ref{fig:profile_neumann}, respectively. Here, we selected three points on the sphere, such as $x_1=0.0$ (i.e. the front of the sphere), $x_1=0.5$ (i.e. the middle), and $x_1=1.0$ (i.e. the back); see Figure~\ref{fig:sphere}.

With regard to the Dirichlet problem in Table~\ref{tab:dirichlet}, the OBIE was unstable and diverged for any $d$, problem size, and algorithm. Meanwhile, the BMBIE successfully ran for $d=1$. These results are consistent with the semi-analytical studies on the stability of the TDBEM~\cite{fukuhara2019, chiyoda2019}. However, the high orders, i.e. $d=2$ and $3$, did not work even for the BMBIE, although FAST8 and FAST12 did not diverge but caused large errors when $d=2$ and $\Ns=2880$. These results of the BMBIE are indeed reflected in the profiles in Figure~\ref{fig:profile_dirichlet}.

\begin{table}[H]
  \begin{center}
  \caption{Numerical result for the homogeneous \textbf{Dirichlet} problem. An error in \red{red} designates an unstable or obviously inaccurate result, whereas an error in \darkgreen{green} means that it is less than the relative error of 1 (i.e. 100\%) but the corresponding solution appears to diverge when it is drawn, as shown in Figure~\ref{fig:profile_dirichlet}. In the column of BIE, OR and BM stand for the OBIE and BMBIE, respectively.}
  \label{tab:dirichlet}
  %
  %
  \renewcommand{\arraystretch}{1.15}
  \begin{tabular}{|c|c|c|l|r|r|l|r|r|l|r|r|}
    \hline
    BIE  & Algo. & $d$ & \multicolumn{3}{c|}{$\Ns=2880$, $\Nt=240$} & \multicolumn{3}{c|}{$\Ns=5120$, $\Nt=320$} & \multicolumn{3}{c|}{$\Ns=11520$, $\Nt=480$}\\
    \cline{4-12}    
         &        &   & Error & Time & Mem & Error & Time & Mem & Error & Time & Mem\\
    \hline
    OR & CONV   & 1 & \prerrx{6.682938e+06} & \prtimex{ 43.995481} & \prmemx{9.51572400000000000000}&
\prerrx{1.333051e+08} & \prtimex{ 204.329454} & \prmemx{38.41043600000000000000}&
\prerrx{1.290505e+11} & \prtimex{ 1479.420190} & \prmemx{253.96342800000000000000}\\
    \cline{3-12}
         &        & 2 & \prerrx{5.865050e+132} & \prtimex{ 50.699072} & \prmemx{9.82974400000000000000}&
\prerrxnan{nan} & \prtimexnan{ 203.421122} & \prmemxnan{39.40548800000000000000}&
\prerrxnan{nan} & \prtimexnan{ 1778.679658} & \prmemxnan{258.15156800000000000000}\\
    \cline{3-12}
         &        & 3 &\prerrxnan{-nan} & \prtimexnan{ 52.617927} & \prmemxnan{10.14474000000000000000}&
\prerrxnan{nan} & \prtimexnan{ 197.311690} & \prmemxnan{40.39474400000000000000}&
\prerrxnan{-nan} & \prtimexnan{ 1812.824705} & \prmemxnan{261.05211600000000000000}\\
    \cline{2-12}
         & FAST8  & 1 &\prerrx{7.024663e+04} & \prtimex{ 29.709239} & \prmemx{2.78324400000000000000}&
\prerrx{9.131507e+05} & \prtimex{ 142.832581} & \prmemx{5.46992000000000000000}&
\prerrx{6.438924e+06} & \prtimex{ 284.834478} & \prmemx{8.19636000000000000000}\\
    \cline{3-12}
         &        & 2 &\prerrx{5.864472e+132} & \prtimex{ 33.693454} & \prmemx{3.00716800000000000000}&
\prerrxnan{-nan} & \prtimexnan{ 144.136637} & \prmemxnan{5.81200800000000000000}&
\prerrxnan{-nan} & \prtimexnan{ 334.912983} & \prmemxnan{8.71730000000000000000}\\
    \cline{3-12}
         &        & 3 &\prerrxnan{-nan} & \prtimexnan{ 39.775138} & \prmemxnan{3.26251200000000000000}&
\prerrxnan{-nan} & \prtimexnan{ 175.645375} & \prmemxnan{6.19957600000000000000}&
\prerrxnan{-nan} & \prtimexnan{ 376.769743} & \prmemxnan{9.31081600000000000000}\\
    \cline{2-12}
         & FAST12 & 1 &\prerrx{6.114442e+05} & \prtimex{ 182.631231} & \prmemx{10.40944400000000000000}&
\prerrx{9.113309e+06} & \prtimex{ 1473.677108} & \prmemx{13.71750400000000000000}&
\prerrx{1.107870e+08} & \prtimex{ 3107.745621} & \prmemx{17.55230800000000000000}\\
    \cline{3-12}
         &        & 2 &\prerrx{5.864758e+132} & \prtimex{ 199.091765} & \prmemx{11.36292800000000000000}&
\prerrxnan{-nan} & \prtimexnan{ 1544.818384} & \prmemxnan{14.87169600000000000000}&
\prerrxnan{-nan} & \prtimexnan{ 3302.999403} & \prmemxnan{19.08696800000000000000}\\
    \cline{3-12}
         &        & 3 &\prerrxnan{-nan} & \prtimexnan{ 218.165680} & \prmemxnan{12.36218400000000000000}&
\prerrxnan{nan} & \prtimexnan{ 1616.710424} & \prmemxnan{16.40543600000000000000}&
\prerrxnan{nan} & \prtimexnan{ 3496.198900} & \prmemxnan{21.04776000000000000000}\\
    \hline
    BM & CONV   & 1 & \prerr{2.515757e-02} & \prtime{ 50.222437} & \prmem{9.51584400000000000000}&
\prerr{1.821510e-02} & \prtime{ 189.818241} & \prmem{38.41072000000000000000}&
\prerr{1.171056e-02} & \prtime{ 1749.584227} & \prmem{253.96300400000000000000}\\
    \cline{3-12}
         &        & 2 & \prerrx{1.776316e+00} & \prtimex{ 48.096272} & \prmemx{9.83059200000000000000}&
\prerrx{2.695220e+01} & \prtimex{ 194.879538} & \prmemx{39.40226400000000000000}&
\prerrx{7.832280e+04} & \prtimex{ 1838.714827} & \prmemx{258.14662000000000000000}\\
    \cline{3-12}
         &        & 3 & \prerrx{5.802913e+131} & \prtimex{ 54.203641} & \prmemx{10.14411200000000000000}&
\prerrxnan{nan} & \prtimexnan{ 221.330897} & \prmemxnan{40.39819600000000000000}&
\prerrxnan{nan} & \prtimexnan{ 1991.722917} & \prmemxnan{261.07000000000000000000}\\
    \cline{2-12}
         & FAST8  & 1 & \prerr{7.342976e-02} & \prtime{ 30.163219} & \prmem{2.78277200000000000000}&
\prerr{6.546132e-02} & \prtime{ 131.569108} & \prmem{5.46936000000000000000}&
\prerr{5.743643e-02} & \prtime{ 283.594463} & \prmem{8.19516000000000000000}\\
    \cline{3-12}
         &        & 2 & \prerrG{1.848601e-01} & \prtimeG{ 32.458046} & \prmemG{3.00746400000000000000}&
\prerrx{1.299512e+01} & \prtimex{ 143.289126} & \prmemx{5.81137600000000000000}&
\prerrx{5.648963e+02} & \prtimex{ 304.860143} & \prmemx{8.71708400000000000000}\\
    \cline{3-12}
         &        & 3 & \prerrx{5.802914e+131} & \prtimex{ 35.902669} & \prmemx{3.26348400000000000000}&
\prerrxnan{-nan} & \prtimexnan{ 165.971432} & \prmemxnan{6.20110400000000000000}&
\prerrxnan{-nan} & \prtimexnan{ 355.725479} & \prmemxnan{9.31316000000000000000}\\
    \cline{2-12}
         & FAST12 & 1 & \prerr{4.670451e-02} & \prtime{ 181.437136} & \prmem{10.41415600000000000000}&
\prerr{3.700801e-02} & \prtime{ 1457.177099} & \prmem{13.77927200000000000000}&
\prerr{3.395901e-02} & \prtime{ 3142.662314} & \prmem{17.54205600000000000000}\\
    \cline{3-12}
         &        & 2 & \prerrG{7.341495e-01} & \prtimeG{ 196.825202} & \prmemG{11.35486400000000000000}&
\prerrx{6.136218e+00} & \prtimex{ 1513.124758} & \prmemx{14.88749200000000000000}&
\prerrx{1.478331e+02} & \prtimex{ 3279.502863} & \prmemx{19.08648000000000000000}\\
    \cline{3-12}
         &        & 3 & \prerrx{5.802784e+131} & \prtimex{ 218.781619} & \prmemx{12.35438800000000000000}&
\prerrxnan{-nan} & \prtimexnan{ 1586.535962} & \prmemxnan{16.48002800000000000000}&
\prerrxnan{nan} & \prtimexnan{ 3466.922978} & \prmemxnan{21.06154400000000000000}\\
    \hline
  \end{tabular}
  \end{center}
\end{table}

For $d=1$ in the Neumann problem, Table~\ref{tab:neumann} as well as Figure~\ref{fig:profile_neumann}(a) show that the OBIE worked well regardless of the algorithm. This result is consistent with that of the previous work~\cite{takahashi2014}, where the same Neumann problem was solved. On the other hand, although the BMBIE with $d=1$ did not diverge, its accuracy ($\sim 10^{-1}$, i.e. $10\%$) was worse than that of the OBIE ($\sim 1\%$). 

For $d=2$ and $3$, the OBIE diverged in all the cases. However, the BMBIE with $d=2$ achieved the best accuracy ($\lesssim 10^{-2}$) in all the cases (see Figure~\ref{fig:profile_neumann}(b)), whereas the BMBIE with $d=3$ diverged. It should be noted that, as Table~\ref{tab:neumann} shows, the BMBIE with $d=2$ improved the accuracy without increasing the computation time and memory usage significantly in comparison with the BMBIE with $d=1$.

In all the stably computed cases, it is observed that the error decreases as the problem size increases. This is due to the reduction in the discretisation error.

In the case of the fast TDBEM (only when it did not diverge), we can {observe a trade-off between numerical accuracy and computation time as well as memory usage depending on $\Ps$ and $\Pt$. Because the convergence rate with respect to $\Ps$ and $\Pt$ is extremely small, FAST12 (i.e. fast TDBEM using $\Ps=\Pt=12$) was slower than CONV in all the cases. However, because the conventional algorithm has a computational complexity of $O(\Ns^2\Nt)$ and requires a very large memory, it is impractical to analyse larger-scale problems with the convectional TDBEM.

It should be noted that, because the moments and local coefficients can be treated as $\Ps^3\Pt$-dimensional vectors, the computational complexity of each FMM operator includes the pre-factor $\Ps^3\Pt$ (which is for the creation of moments and the evaluation with local coefficients), $\Ps^6\Pt^2$ (for the M2M and L2L), or $\Ps^3\Pt(\log\Ps+\log\Pt)$ (for the M2L; recall Remark~\ref{remark:m2l_complexity}). This is why FAST8 (where $\Ps=8$ and $\Pt=8$) is much faster than FAST12 (where $\Ps=\Pt=12$).

\begin{table}[H]
  \centering
  \caption{Numerical result for the homogeneous \textbf{Neumann} problem.}
  \label{tab:neumann}
  \renewcommand{\arraystretch}{1.15}
  %
%
\begin{tabular}{|c|c|c|l|r|r|l|r|r|l|r|r|}
  \hline
  BIE  & Algo. & $d$ & \multicolumn{3}{c|}{$\Ns=2880$, $\Nt=240$} & \multicolumn{3}{c|}{$\Ns=5120$, $\Nt=320$} & \multicolumn{3}{c|}{$\Ns=11520$, $\Nt=480$}\\
  \cline{4-12}    
       &        &   & Error & Time & Mem & Error & Time & Mem & Error & Time & Mem\\
  \hline
  OR   & CONV   & 1 & \prerr{3.423881e-02} & \prtime{ 41.667815} & \prmem{9.51553200000000000000} & \prerr{1.775474e-02} & \prtime{ 220.433312} & \prmem{38.41052400000000000000} & \prerr{1.114806e-02} & \prtime{ 1612.395621} & \prmem{253.96332800000000000000}\\
  \cline{3-12}
       &        & 2 & \prerrx{1.119688e+13} & \prtimex{ 49.916084} & \prmemx{9.82917200000000000000} & \prerrx{2.863237e+15} & \prtimex{ 205.185105} & \prmemx{39.40246800000000000000} & \prerrx{9.841636e+20} & \prtimex{ 1737.645206} & \prmemx{258.14404400000000000000}\\
  \cline{3-12}
       &        & 3 & \prerrx{4.463055e+134} & \prtimex{ 49.158494} & \prmemx{10.14272400000000000000} & \prerrxnan{nan} & \prtimex{ 188.572615} & \prmemx{40.39408800000000000000} & \prerrxnan{nan} & \prtimex{ 1884.642814} & \prmemx{261.05238000000000000000}\\
  \cline{2-12}
       & FAST8  & 1 & \prerr{3.782309e-02} & \prtime{ 31.122201} & \prmem{2.78283200000000000000} & \prerr{2.938634e-02} & \prtime{ 135.641399} & \prmem{5.46929200000000000000} & \prerr{2.550533e-02} & \prtime{ 280.669735} & \prmem{8.19614000000000000000}\\
  \cline{3-12}
       &        & 2 & \prerrx{2.725866e+11} & \prtimex{ 31.758443} & \prmemx{3.00714000000000000000} & \prerrx{5.035767e+13} & \prtimex{ 146.182652} & \prmemx{5.81114400000000000000} & \prerrx{3.675531e+15} & \prtimex{ 304.272762} & \prmemx{8.71756400000000000000}\\
  \cline{3-12}
       &        & 3 & \prerrx{4.463098e+134} & \prtimex{ 35.165295} & \prmemx{3.26232400000000000000} & \prerrxnan{-nan} & \prtimex{ 158.015888} & \prmemx{6.19982000000000000000} & \prerrxnan{-nan} & \prtimex{ 355.888574} & \prmemx{9.31112400000000000000}\\
  \cline{2-12}
       & FAST12 & 1 & \prerr{3.178029e-02} & \prtime{ 181.491687} & \prmem{10.40968400000000000000} & \prerr{2.001298e-02} & \prtime{ 1484.024042} & \prmem{13.71455600000000000000} & \prerr{1.226066e-02} & \prtime{ 3145.381449} & \prmem{17.55216400000000000000}\\
  \cline{3-12}
       &        & 2 & \prerrx{3.120689e+12} & \prtimex{ 198.602343} & \prmemx{11.36352400000000000000} & \prerrx{4.218195e+14} & \prtimex{ 1517.216924} & \prmemx{14.87665200000000000000} & \prerrx{2.586086e+18} & \prtimex{ 3236.984402} & \prmemx{19.08688400000000000000}\\
  \cline{3-12}
       &        & 3 & \prerrx{4.463061e+134} & \prtimex{ 217.511305} & \prmemx{12.36208400000000000000} & \prerrxnan{-nan} & \prtimex{ 1620.299565} & \prmemx{16.56409600000000000000} & \prerrxnan{nan} & \prtimex{ 3489.785906} & \prmemx{21.03626400000000000000}\\
  \hline
  BM & CONV     & 1 & \prerr{1.909606e-01} & \prtime{ 50.838090} & \prmem{9.51581200000000000000} & \prerr{1.469617e-01} & \prtime{ 203.388391} & \prmem{38.41052400000000000000} & \prerr{1.002265e-01} & \prtime{ 1843.253859} & \prmem{253.96286400000000000000}\\
  \cline{3-12}
       &        & 2 & \prerr{1.701852e-02} & \prtime{ 51.747418} & \prmem{9.83038400000000000000} & \prerr{9.840038e-03} & \prtime{ 207.408439} & \prmem{39.40199600000000000000} & \prerr{5.101089e-03} & \prtime{ 1974.747267} & \prmem{258.14678800000000000000}\\
  \cline{3-12}
       &        & 3 & \prerrx{3.325767e+45} & \prtimex{ 50.538438} & \prmemx{10.14413600000000000000} & \prerrx{8.132068e+61} & \prtimex{ 197.680388} & \prmemx{40.39472400000000000000} & \prerrx{7.637551e+95} & \prtimex{ 1824.263323} & \prmemx{261.06174000000000000000}\\
  \cline{2-12}
       & FAST8  & 1 & \prerr{1.878164e-01} & \prtime{ 30.542090} & \prmem{2.78262400000000000000} & \prerr{1.450580e-01} & \prtime{ 133.216258} & \prmem{5.46947600000000000000} & \prerr{9.991420e-02} & \prtime{ 284.456034} & \prmem{8.19431600000000000000}\\
  \cline{3-12}
       &        & 2 & \prerr{2.987187e-02} & \prtime{ 31.298734} & \prmem{3.00753200000000000000} & \prerr{2.125243e-02} & \prtime{ 139.494509} & \prmem{5.81124000000000000000} & \prerr{1.789968e-02} & \prtime{ 300.552323} & \prmem{8.71708400000000000000}\\
  \cline{3-12}
       &        & 3 & \prerrx{3.318988e+45} & \prtimex{ 35.871904} & \prmemx{3.26358000000000000000} & \prerrx{8.124275e+61} & \prtimex{ 159.122186} & \prmemx{6.20010800000000000000} & \prerrx{7.408303e+95} & \prtimex{ 308.628512} & \prmemx{9.31270400000000000000}\\
  \cline{2-12}
       & FAST12 & 1 & \prerr{1.899199e-01} & \prtime{ 173.147166} & \prmem{10.41401600000000000000} & \prerr{1.452102e-01} & \prtime{ 1481.631242} & \prmem{13.76982800000000000000} & \prerr{9.989237e-02} & \prtime{ 3149.940574} & \prmem{17.54216400000000000000}\\
  \cline{3-12}
       &        & 2 & \prerr{2.093128e-02} & \prtime{ 195.127808} & \prmem{11.35491200000000000000} & \prerr{1.274179e-02} & \prtime{ 1536.499944} & \prmem{14.90843200000000000000} & \prerr{9.297933e-03} & \prtime{ 3237.570170} & \prmem{19.08647600000000000000}\\
  \cline{3-12}
       &        & 3 & \prerrx{3.340136e+45} & \prtimex{ 221.884852} & \prmemx{12.35417200000000000000} & \prerrx{8.156051e+61} & \prtimex{ 1581.091829} & \prmemx{16.42850400000000000000} & \prerrx{7.643676e+95} & \prtimex{ 3407.525811} & \prmemx{21.00894000000000000000}\\
  \hline
\end{tabular}

\end{table}

\begin{figure}[H]
  \centering
  \includegraphics[height=.4\textheight]{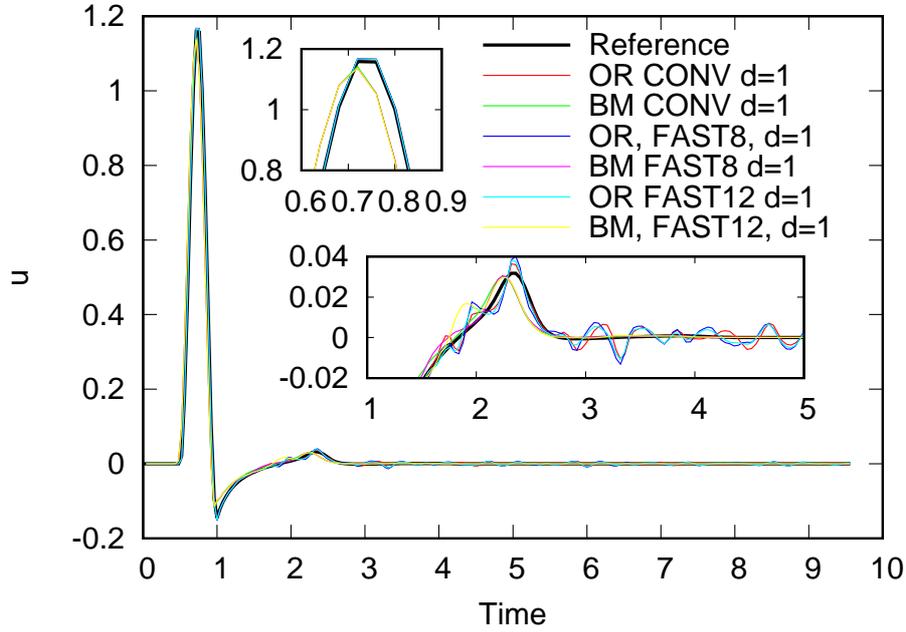}\\
  (a) Comparison for $d=1$.\\
  \includegraphics[height=.4\textheight]{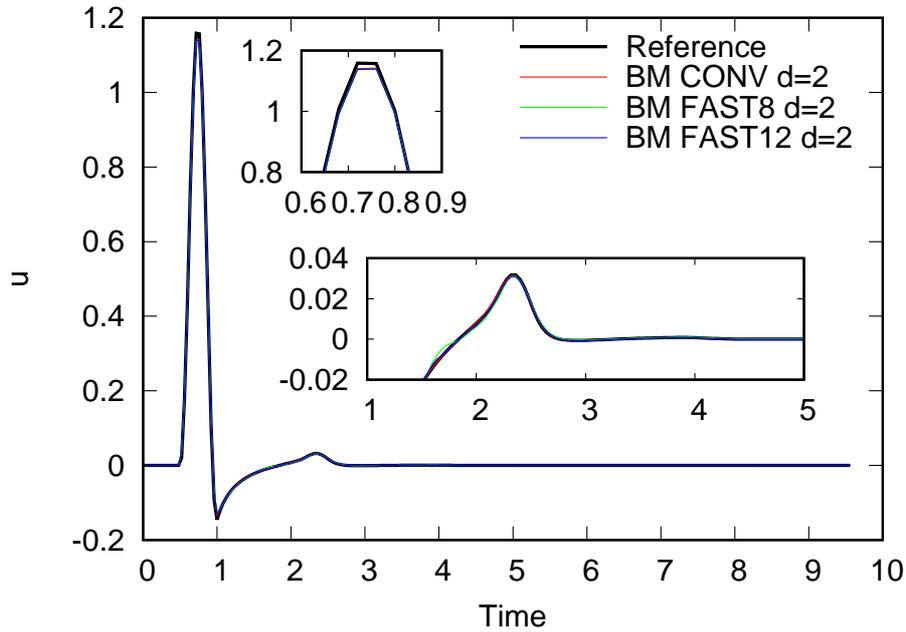}\\
  (b) Comparison for the BMBIE with $d=2$.
  \caption{Profiles of $u$ in the homogeneous \textbf{Neumann} problem in the case of $\Ns=2880$ and the middle point, i.e. $x_1=0.5$ (recall Figure~\ref{fig:sphere}). OR and BM in the legend stand for OBIE and BMBIE, respectively.}
  \label{fig:profile_neumann}
\end{figure}

\subsubsection{Discussion}\label{s:sphere_discuss}

The present TDBEM, which employs the B-spline temporal basis of order $d$, tends to become unstable as $d$ becomes large. We discuss the stability from the arithmetic viewpoint. We suspect that cancellation of significant digits can occur when summing up $d+2$ coefficients $U_{ij}^{(\alpha-\beta-\kappa)}$ and $W_{ij}^{(\alpha-\beta-\kappa)}$ over $\kappa$ in (\ref{eq:discretised_single_layer_potential}) and (\ref{eq:discretised_double_layer_potential}). We recall that the summations over $\kappa$ result from the decomposition of the B-spline basis in (\ref{eq:decompose}). These coefficients are obtained by integrating the kernel functions $U$ in (\ref{eq:U}) and $\bm{W}$ in (\ref{eq:W}), that is,
\begin{eqnarray*}
  U(\bm{x},\bm{y},t,s)
  &=&\frac{(T-r)_+^d}{r},\\
  \bm{W}(\bm{x},\bm{y},t,s)
  &=&\left(d(T-r)^{d-1}_+T-(d-1)(T -r)^d_+\right)\frac{\bm{x}-\bm{y}}{r^3}
  =\left((T-r)^d_++dr(T -r)^{d-1}_+\right)\frac{\bm{x}-\bm{y}}{r^3},
\end{eqnarray*}
where we let $T:=c(t-s)$ and $r:=\abs{\bm{x}-\bm{y}}$ for simplicity. As we can see, these functions behave like the $d$th-order polynomial $\sim T^d$ as the time $t$ becomes large relative to the source time $s$. Therefore, the underlying coefficients $U_{ij}^{(\alpha-\beta-\kappa)}$ and $W_{ij}^{(\alpha-\beta-\kappa)}$ can also increase very much with the time difference between $t_\alpha$ and $t_\beta$, i.e. $(\alpha-\beta)\Dt$. Then, a subtraction of such large numbers, which can be slightly different, in the summations over $\kappa$ can cause the cancellation of significant digits. The larger $d$ is, the easier it is for a cancellation to occur.

Another interesting question is why the BMBIE of $d=2$ was more accurate than that of $d=1$ in the Neumann problem (recall Table~\ref{tab:neumann}). The probable reason is that $d=2$ is more suitable for representing the solution than $d=1$; in fact, the solution is smooth rather than piecewise-linear with respect to time, as observed in Figure~\ref{fig:profile_neumann}. However, this appears to contradict the fact that the BMBIE of $d=2$ was less accurate than that of $d=1$ in the Dirichlet problem (recall Table~\ref{tab:dirichlet}). This can be explained by considering the decay rate of the kernel functions. The kernel function of the BMBIE for the Neumann problem, i.e. $\left(\frac{\partial}{\partial n_x}-\frac{1}{c}\frac{\partial}{\partial t}\right)\bm{W}$, behaves as $T^dr^{-3}$, whereas that for the Dirichlet problem, i.e. $\left(\frac{\partial}{\partial n_x}-\frac{1}{c}\frac{\partial}{\partial t}\right)U$, behaves as $T^dr^{-2}$. 
Hence, the kernel function of the Dirichlet problem decays more slowly than that of the Neumann problem. Therefore, in the case of the BMBIE for the Dirichlet problem, the negative effect of $d=2$, i.e. the cancellation of significant digits, might overwhelm the positive effect, i.e. the high-accuracy interpolation.

It is difficult to provide numerical data to support the above argument. To do so, we need to modify our computer program drastically. This is because the underlying summations of $U_{ij}^{(\alpha-\beta-\kappa)}$ and $W_{ij}^{(\alpha-\beta-\kappa)}$ over $\kappa$ are equivalently imposed on the boundary variables $u_j^{\beta-\kappa}$ and $q_j^{\beta-\kappa}$ (to yield the alternative boundary variables $\tau_j^\beta$ and $\sigma_j^\beta$ defined by (\ref{eq:sigma}) and (\ref{eq:tau})). Rewriting our program so that it can handle the summations of $u_j^{\beta}$ and $q_j^{\beta}$ directly is very time consuming and beyond the scope of the present study. The stability and instability for $d\ge 2$ are, thus, an open question.

\subsection{Example 2: Hollow box with an aperture}\label{s:hollow}

Instead of the sphere in the previous example, we considered a more complicated scatterer, i.e. a hollow box with a small aperture. As shown in Figure~\ref{fig:hollow}, the width ($W$), depth ($D$), and height ($H$) of the box are $1.0$, $0.5$, and $1.0$, respectively. The width and depth of the aperture on the top of the box are $0.1$ and $0.2$, respectively. The centre of the aperture is $0.25$ from both the left-hand-side ($-x_1$ side) and the front-side ($-x_2$ side) walls. In addition, an internal partition of length ($F$) $0.5$ is attached to the top side, and its distance ($G$) from the left-hand-side wall is $0.5$. The thickness ($T$) of all the walls as well as the partition is $0.02$. We let all the walls be rigid; i.e. the boundary condition was $q=0$. We discretised the boundary model of the hollow box using the mesh generator Gmsh~\cite{geuzaine2009Gmsh}, specifying the mesh size as $0.04$. Then, the number of boundary elements ($\Ns$) was $12682$.

\begin{figure}[H]
  \centering
  \includegraphics[width=.5\textwidth]{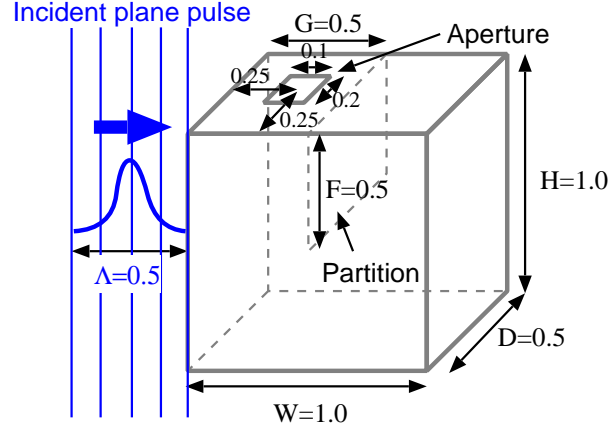}
  \caption{Hollow box with an aperture.}
  \label{fig:hollow}
\end{figure}

We tested the fast TDBEM based on the OBIE or BMBIE, where the order $d$ of the B-spline temporal basis was chosen as $1$ or $2$; thus, we tested four cases. The time-step size ($\Dt$) was $0.02$, and the number of time steps ($\Nt$) was $400$. With respect to the FMM, the numbers of interpolation nodes (i.e. $\Ps$ and $\Pt$) were $8$ in each case.

As a result, only the BMBIE with $d=2$ was stable, whereas the others were unstable. It should be noted that the OBIE with $d=1$ was no longer stable although it was stable in the previous example.

This example as well as the previous does not guarantee that the BMBIE with $d=2$ is always stable for the rigid boundary condition, i.e. $q=0$, which is a practical boundary condition in acoustics. However, we have learned that $d=2$ can be helpful in resolving the instability of the (fast) TDBEM based on $d=1$.

\def\prval#1{{\nprounddigits{3}\npproductsign{\times}\textrm{\numprint{#1}}}}

\section{Application to parameter optimisation}\label{s:app}

By means of the fast TDBEM based on the BMBIE using the B-spline function of order $d=2$, we performed a parameter optimisation for the hollow box in Section~\ref{s:hollow}; recall Figure~\ref{fig:hollow}. This optimisation primarily concerns the stability of the underlying BMBIE in terms of the subtle changes in the boundary shape of the hollow model. In addition, we describe the design of the sound absorbing box called Blast-Wave Eater (BWE)\footnote{The BWEs arraigned in a tunnel can be seen in the following Web page: \url{https://www.jsce.or.jp/prize/tech/files/2014_16.shtml}.}. The BWE is intended to absorb unsteady and low-frequency noise produced by dynamite in tunnel construction. Because the BWE is constructed from plywood boards, the noise is expected to lose its acoustic energy by vibrating the plywood boards. Hence, to maximise the noise reduction, we can design some geometrical parameters of the box so that the sound pressure can excite a structural eigenmode of the box. Strictly speaking, the analysis should be treated as a structural-acoustic coupling problem. However, we approximately consider maximising the sound pressure $u$ on a specified surface, denoted by $Q$, inside the box, so that it can vibrate more or less if it is flexible. To be specific, when an incident field $\uin$ is given, we considered maxmising the following objective function:
\begin{eqnarray*}
  J(\Theta):=\int_0^T\int_{Q}\frac{\abs{u(\bm{x},t)}}{\abs{Q}}\diff S \diff t,
\end{eqnarray*}
where $\Theta$ is a set of parameters to be optimised, $\abs{Q}$ is the area of $Q$, and $T$ denotes the analysis time. In the following analysis, $Q$ was chosen as the ceiling of the right-hand-side room of the cavity. Moreover, we optimised both the length $F$ and the location $G$ of the partition inside the box. We allowed $G$ ($F$, respectively) to vary from $0.4$ to $0.9$ ($0.1$ to $0.9$, respectively); see Figure~\ref{fig:hollow-cross}. Both initial values were set to $0.5$.

To perform the above maximisation, we used the constrained optimisation by linear approximation (COBYLA) method~\cite{powell2007}, which is gradient free and capable of handling inequality constraints. Then, we need to perform the TDBEM for every $\Theta$ that the COBYLA method has determined, and evaluate $J$ from the profile of $u$ on $Q$ computed by the TDBEM.

\begin{figure}[H]
  \centering
  \includegraphics[width=.4\textwidth]{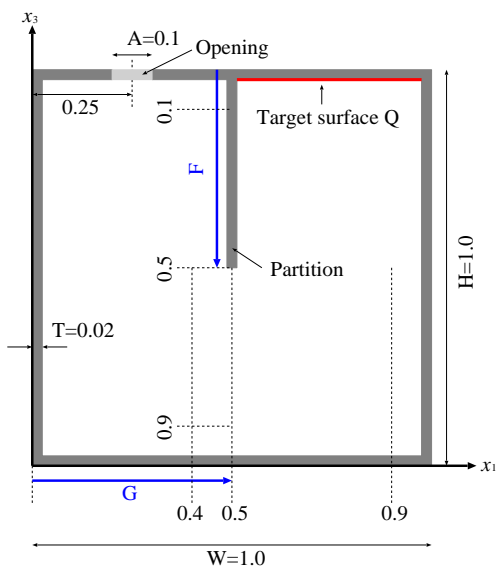}
  \caption{A cross section of the hollow model.}
  \label{fig:hollow-cross}
\end{figure}

The parameters of the TDBEM were the same as those used in the example in Section~\ref{s:hollow}. That is, we let $\Dt=0.02$, $\Nt=400$ and $\Ps=\Pt=8$. Further, the mesh size was about $0.04$ over all the optimisation steps, and the number of boundary elements ($\Ns$) was approximately $13000$, which varied according to the values of $G$ and $F$.

Figure~\ref{fig:bwe-result} plots the history of the objective function $J$ and the design parameters $F$ and $G$ against the number of iteration steps. We terminated the iterations when the relative change of $J$ was less than the prescribed tolerance of $10^{-4}$. Then, after 36 steps, we obtained $J=\prval{2.978834e+00}$ at $F=\prval{5.550683e-01}$ and $G=\prval{6.139577e-01}$. We confirmed that the fast TDBEM using the BMBIE and $d=2$ was stable at any iteration step. 

\begin{figure}[H]     
  \centering          
\iffalse
  \begin{tabular}{cc} 
    \includegraphics[width=.45\textwidth]{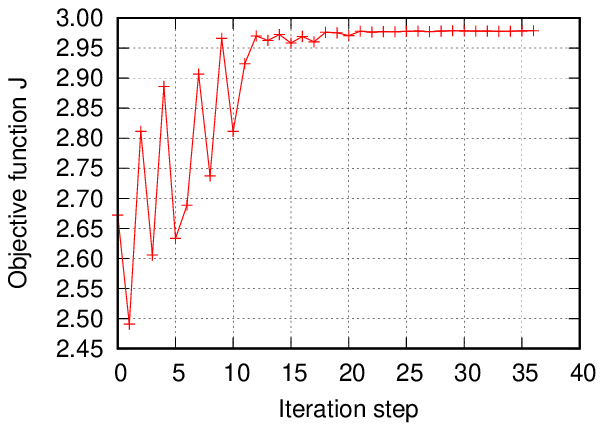}
    &\includegraphics[width=.45\textwidth]{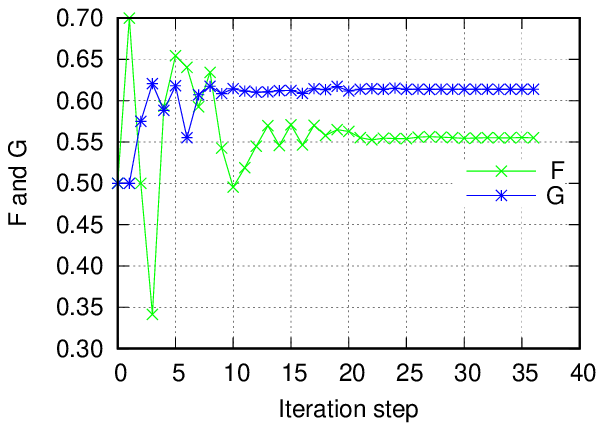}
  \end{tabular}
\else
  \begin{tabular}{cc} 
   \begin{tikzpicture}[scale=0.8] 
  \begin{axis}[
      xlabel={Iteration step},
      ylabel style={align=center},
      ylabel={Objective function $J$},
      xmin=0, xmax=40, xtick={0,5,...,40},
      ymin=2.45, ymax=3.00, ytick={2.45,2.50,2.55,2.60,2.65,2.70,2.75,2.80,2.85,2.90,2.95,3.00},
      y tick label style={
        /pgf/number format/.cd,
        fixed,
        fixed zerofill,
        precision=2,
        /tikz/.cd
      },
      xmajorgrids=true,
      ymajorgrids=true,
      legend pos=south east,
    ]
    \addplot[color=red,mark=*] table [sharp plot, x=Iteration, y=J] {figure/bwefg-RESULT.txt};
  \end{axis}
\end{tikzpicture} & \begin{tikzpicture}[scale=0.8] 
  \begin{axis}[
      xlabel={Iteration step},
      ylabel style={align=center},
      ylabel={$F$, $G$},
      xmin=0, xmax=40, xtick={0,5,...,40},
      ymin=0.30, ymax=0.70, ytick={0.30,0.35,0.40,0.45,0.50,0.55,0.60,0.65,0.70},
      y tick label style={
        /pgf/number format/.cd,
        fixed,
        fixed zerofill,
        precision=2,
        /tikz/.cd
      },
      xmajorgrids=true,
      ymajorgrids=true,
      legend pos=south east,
    ]
    \addplot [color=green, mark=square] table [sharp plot, x=Iteration, y=F] {figure/bwefg-RESULT.txt};
    \addlegendentry{F}
    \addplot [color=blue, mark=triangle] table [sharp plot, x=Iteration, y=G] {figure/bwefg-RESULT.txt};
    \addlegendentry{G}
  \end{axis}
\end{tikzpicture}
  \end{tabular}
\fi
  \caption{History of $J$, $F$ and $G$.}
  \label{fig:bwe-result}
\end{figure}

Figures~\ref{fig:bwe-snapshot0} and \ref{fig:bwe-snapshot} show the sound pressure on the surface at several time steps for the initial and optimum configurations, respectively.

\begin{figure}[H]     
  \centering          
  \begin{tabular}{ccc} 
    \includegraphics[width=.3\textwidth]{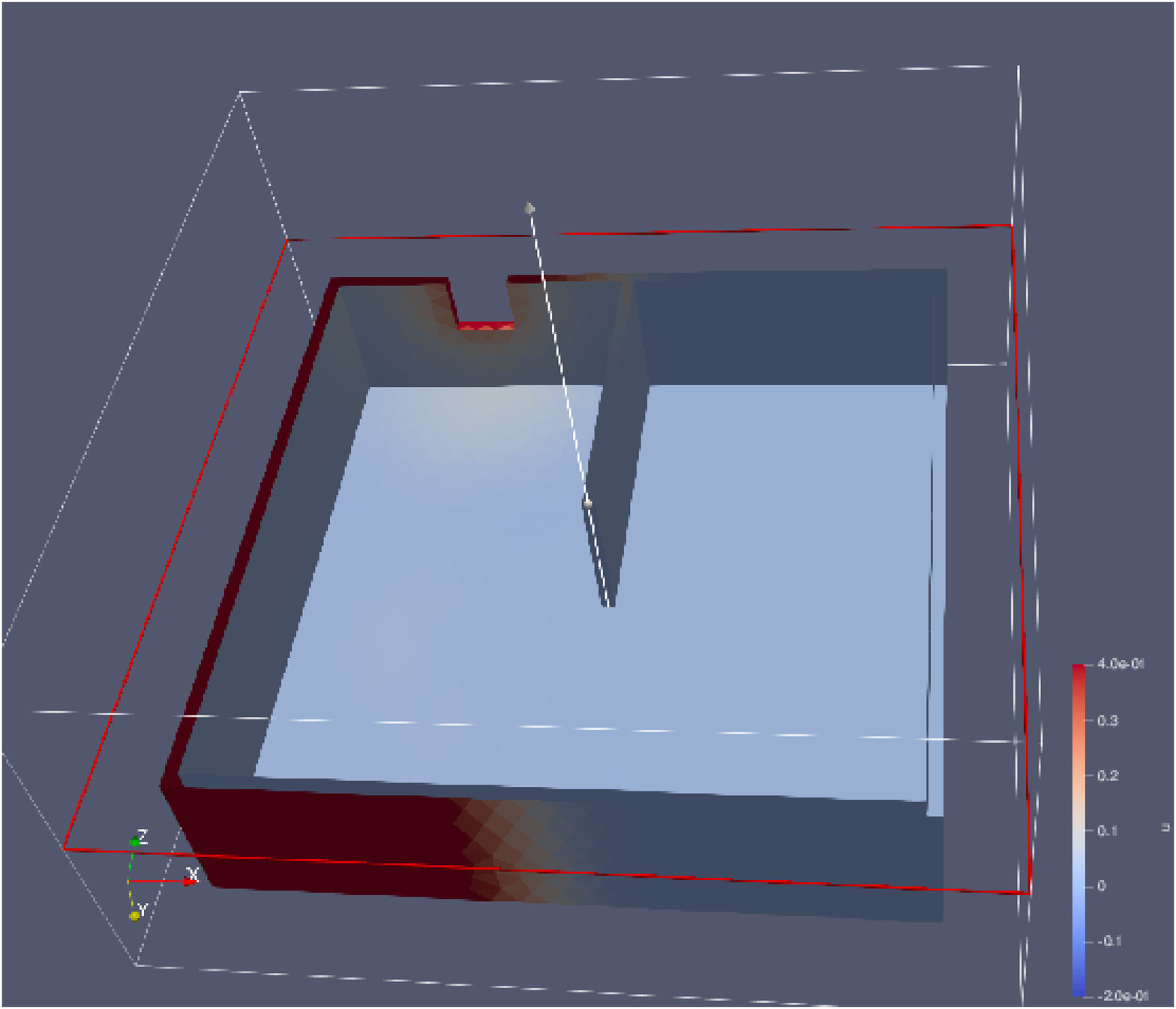}
    &\includegraphics[width=.3\textwidth]{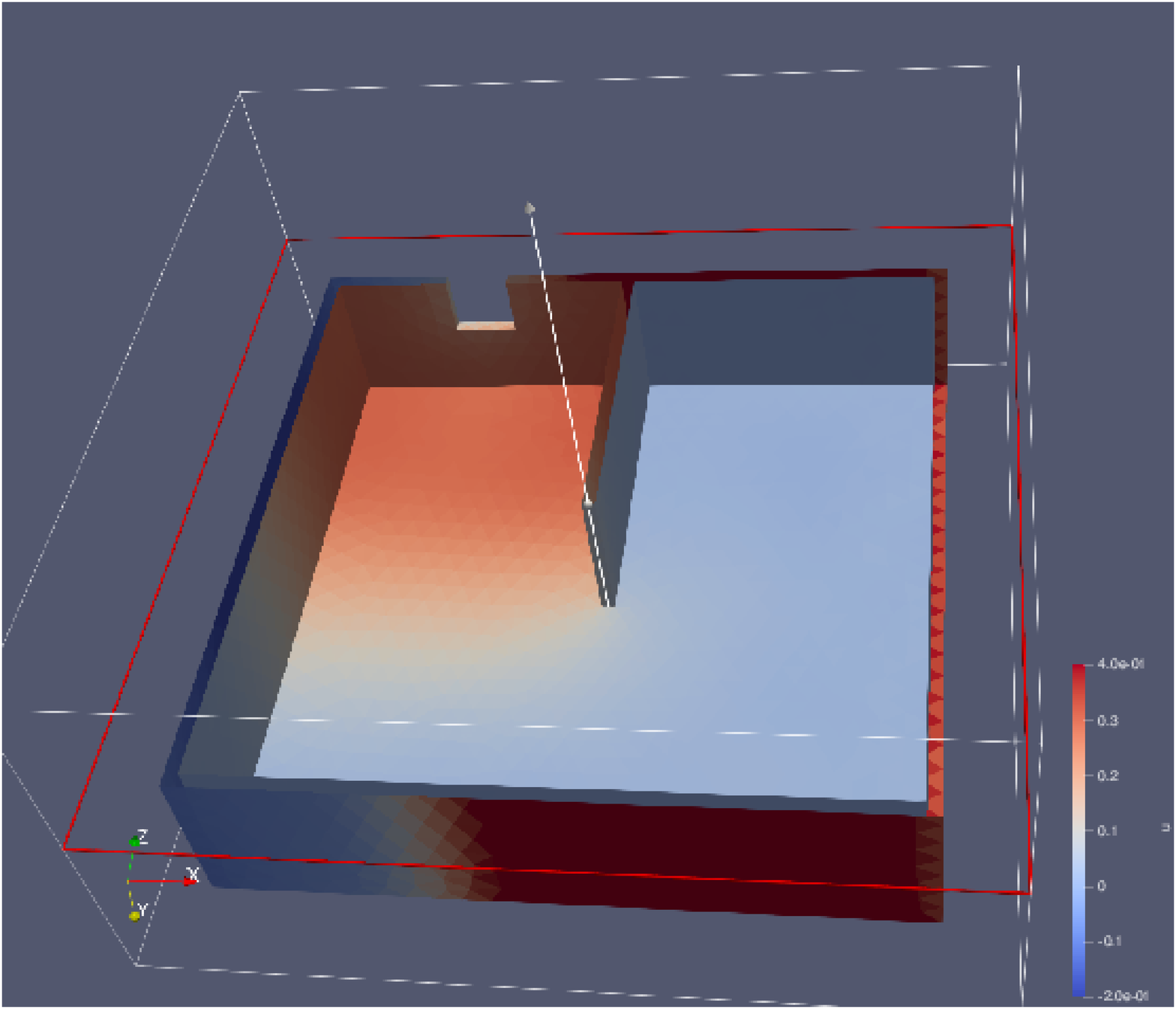}
    &\includegraphics[width=.3\textwidth]{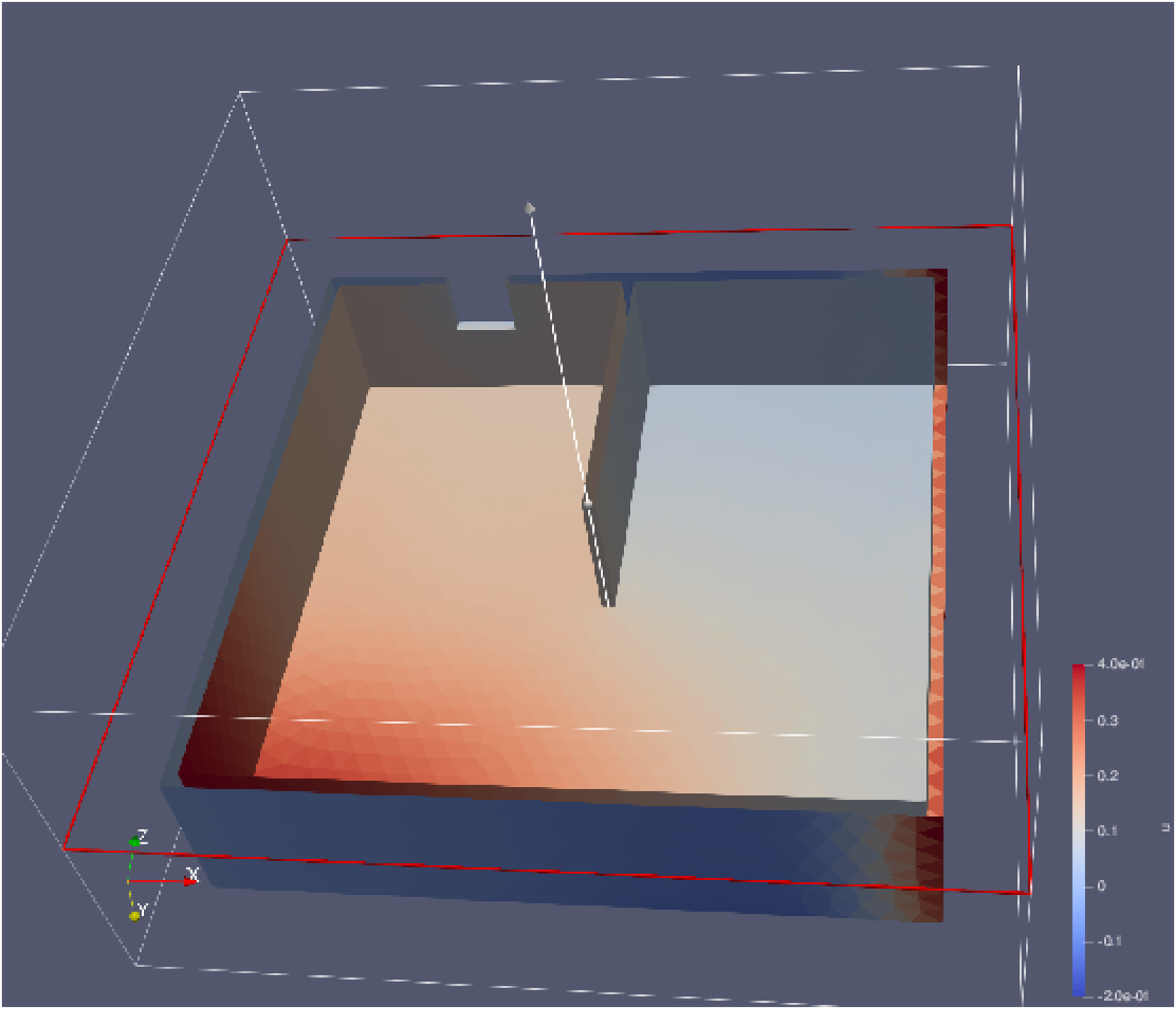}\\
    $t=0.6$ & $t=1.2$ & $t=1.8$\\
    \includegraphics[width=.3\textwidth]{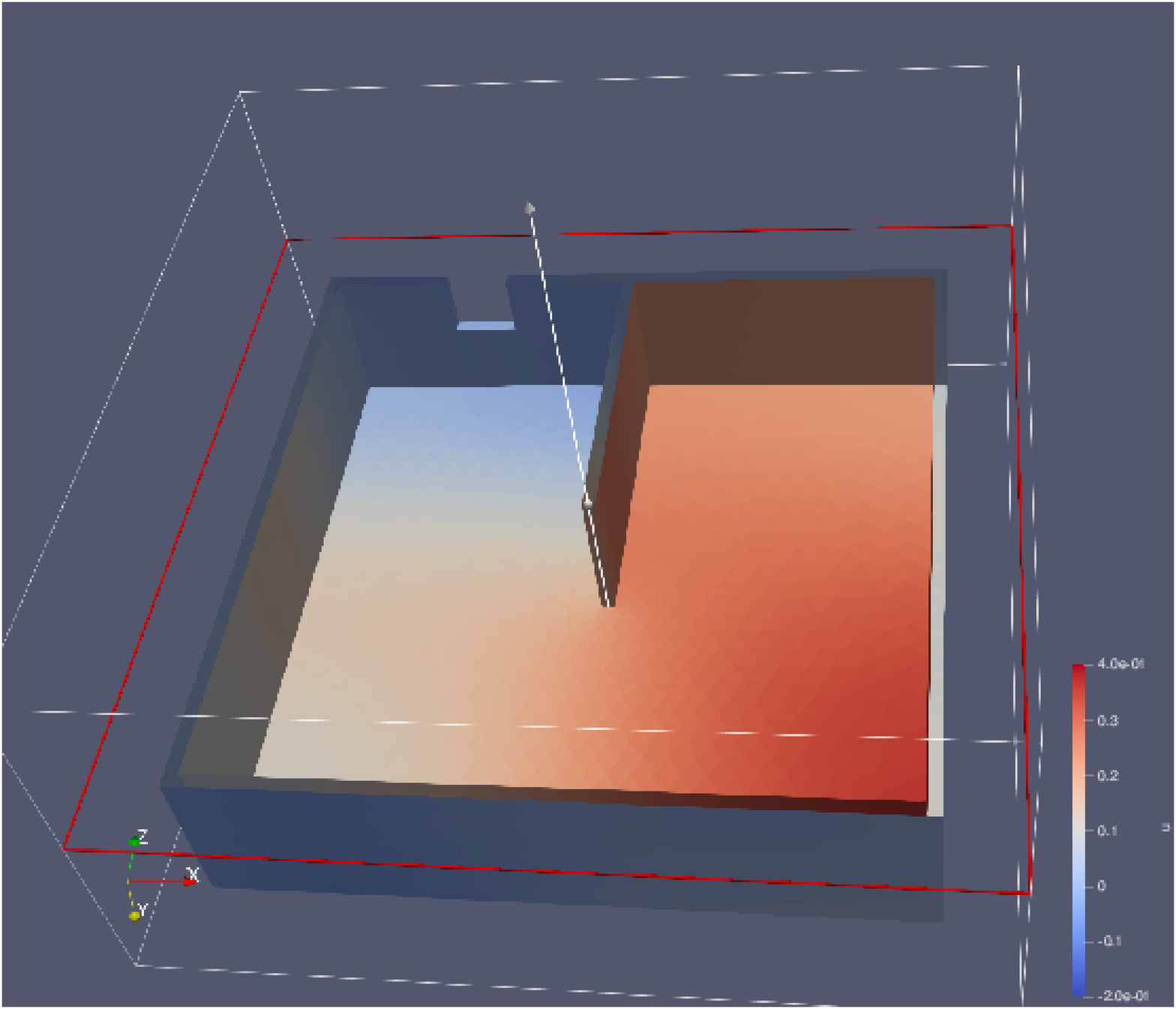}
    &\includegraphics[width=.3\textwidth]{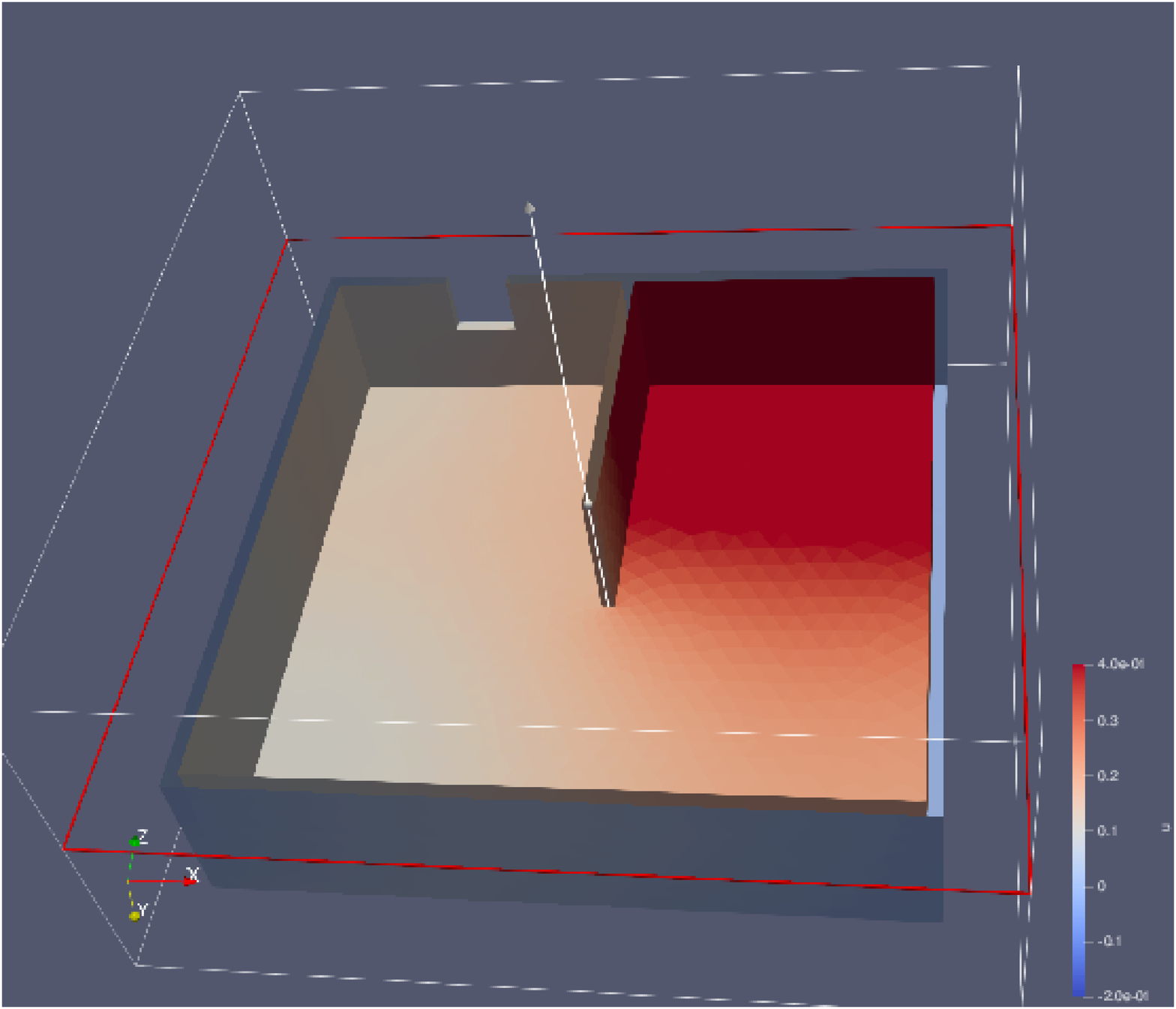}
    &\includegraphics[width=.3\textwidth]{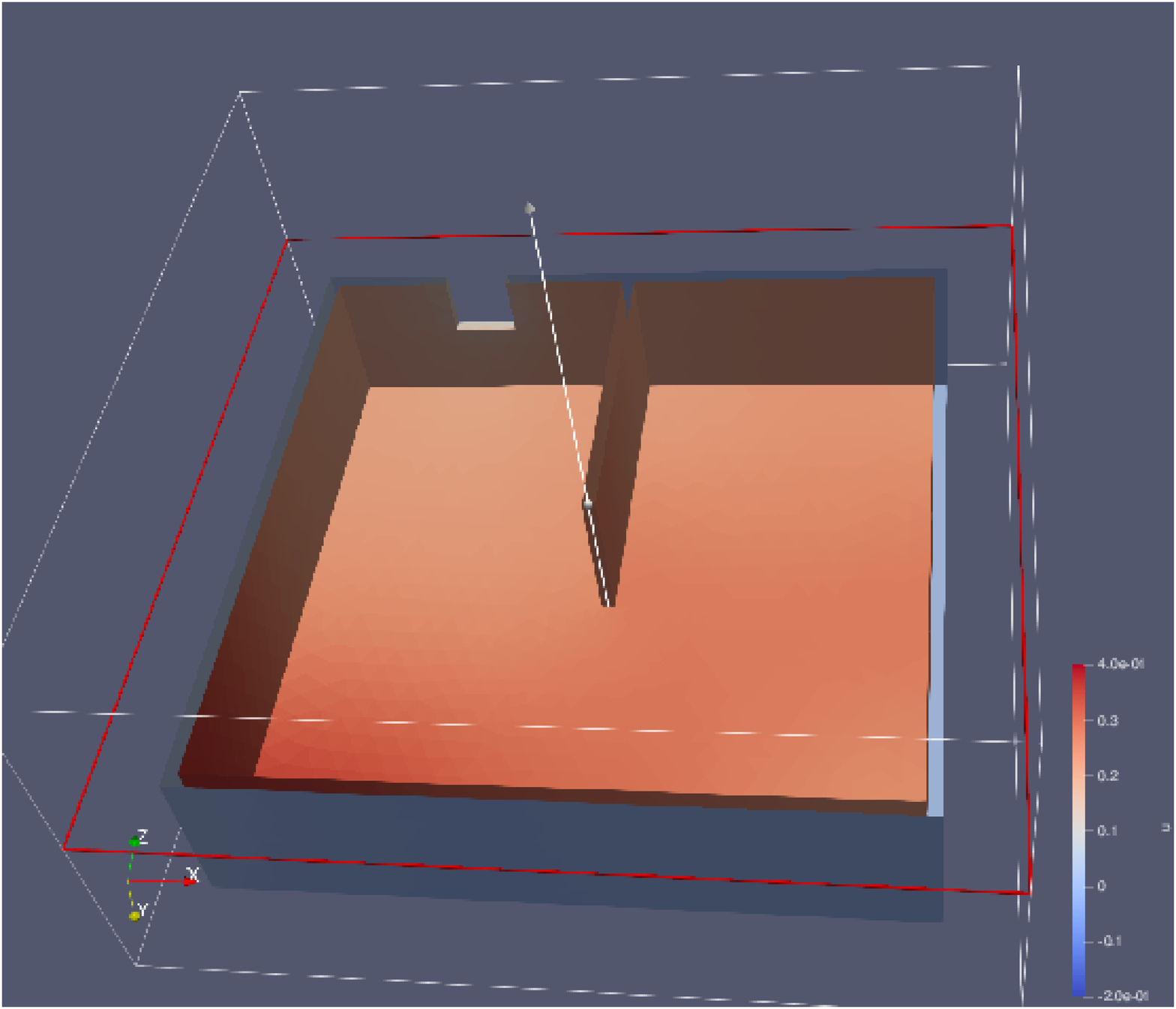}\\
    $t=2.4$ & $t=3.0$ & $t=3.6$\\
    \includegraphics[width=.3\textwidth]{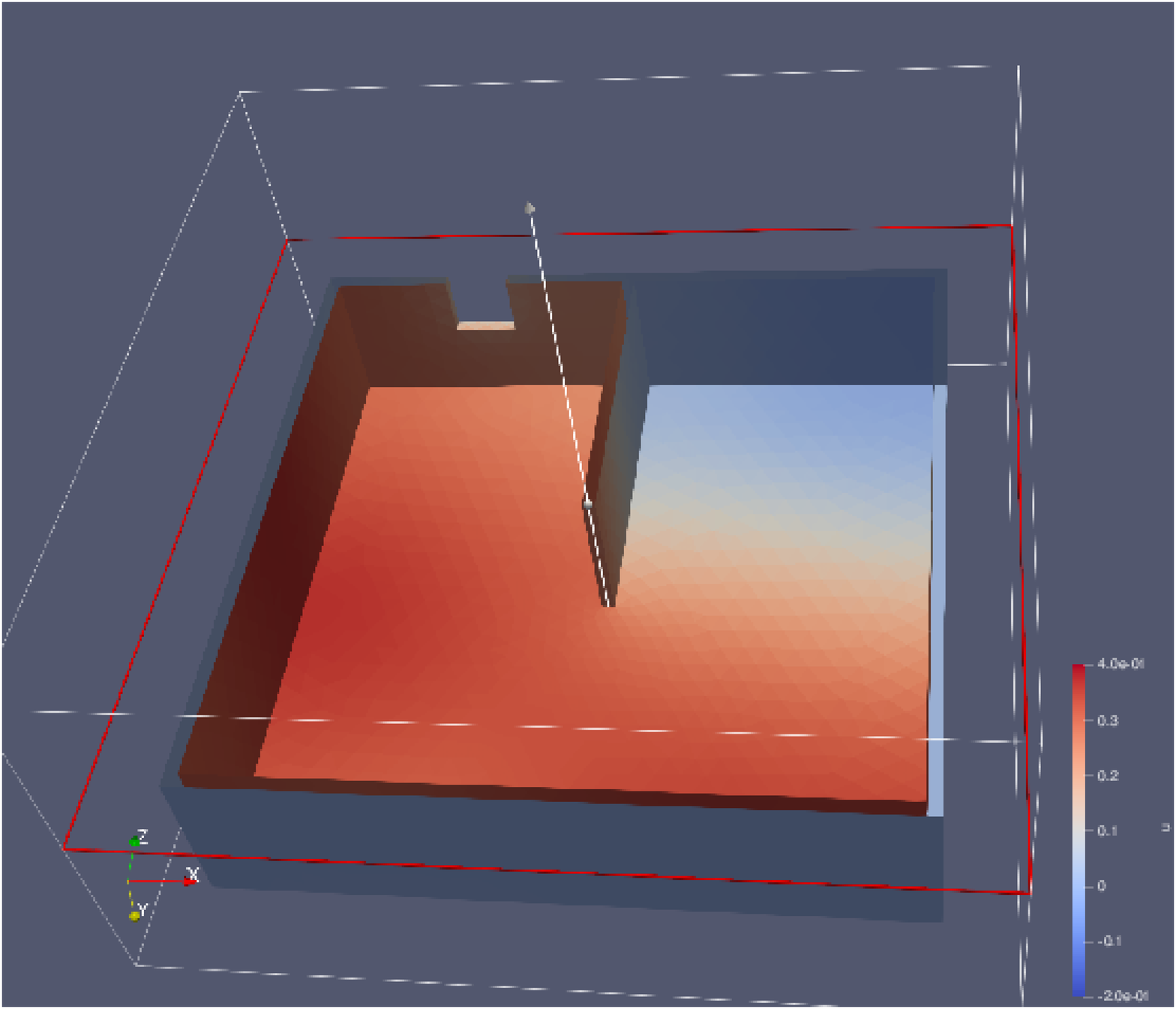}
    &\includegraphics[width=.3\textwidth]{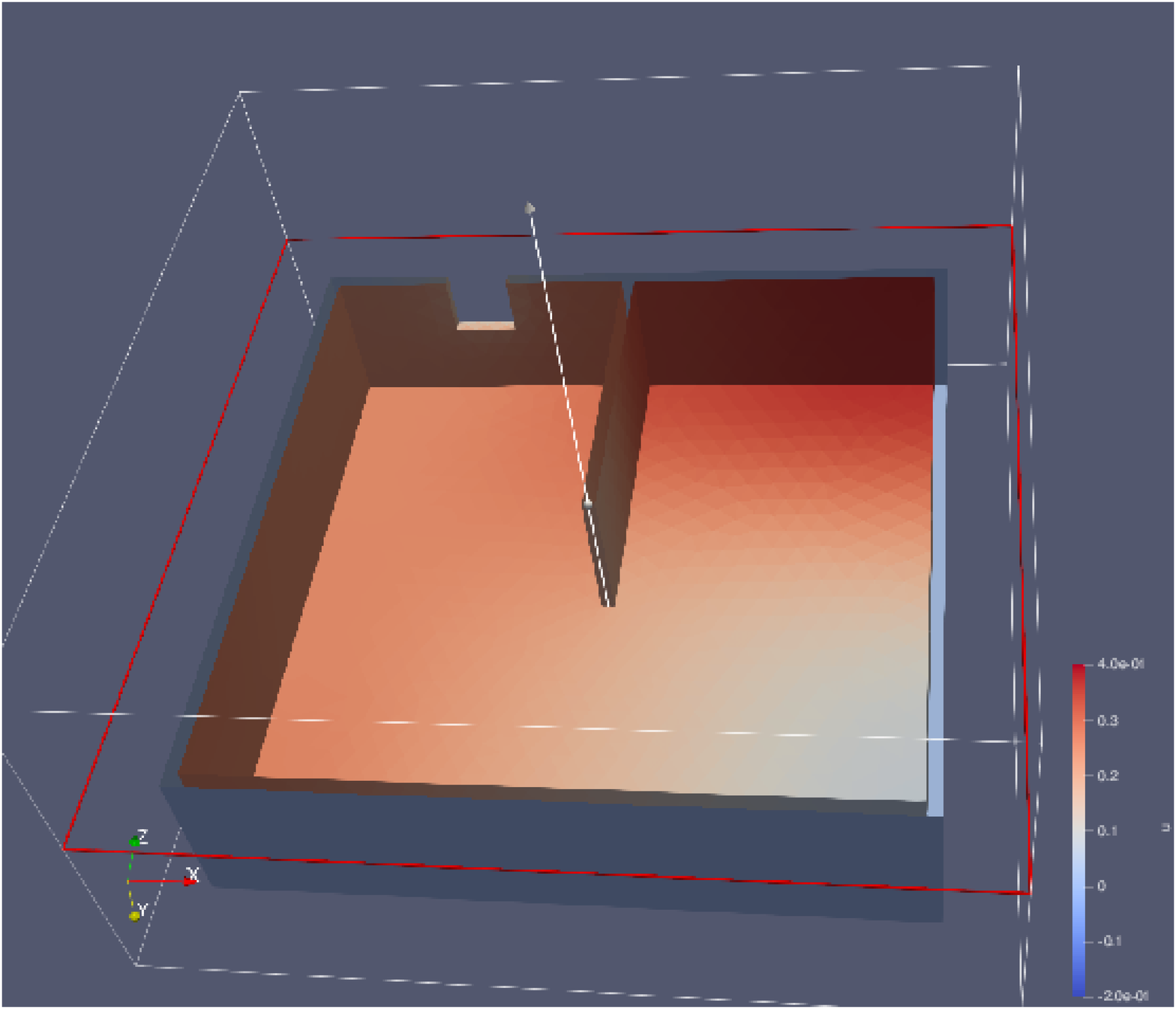}
    &\includegraphics[width=.3\textwidth]{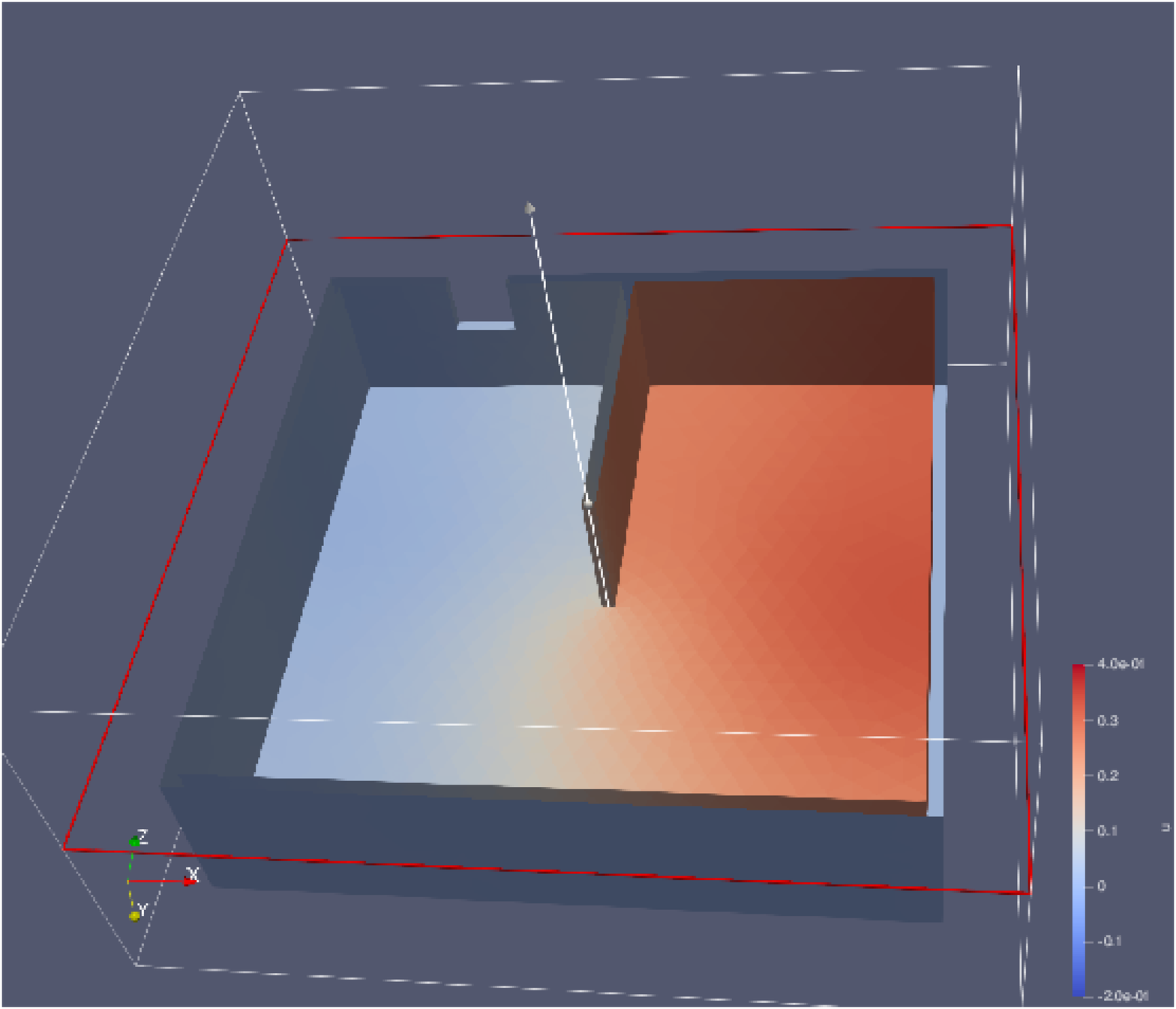}\\
    $t=4.2$ & $t=4.8$ & $t=5.4$\\
    \includegraphics[width=.3\textwidth]{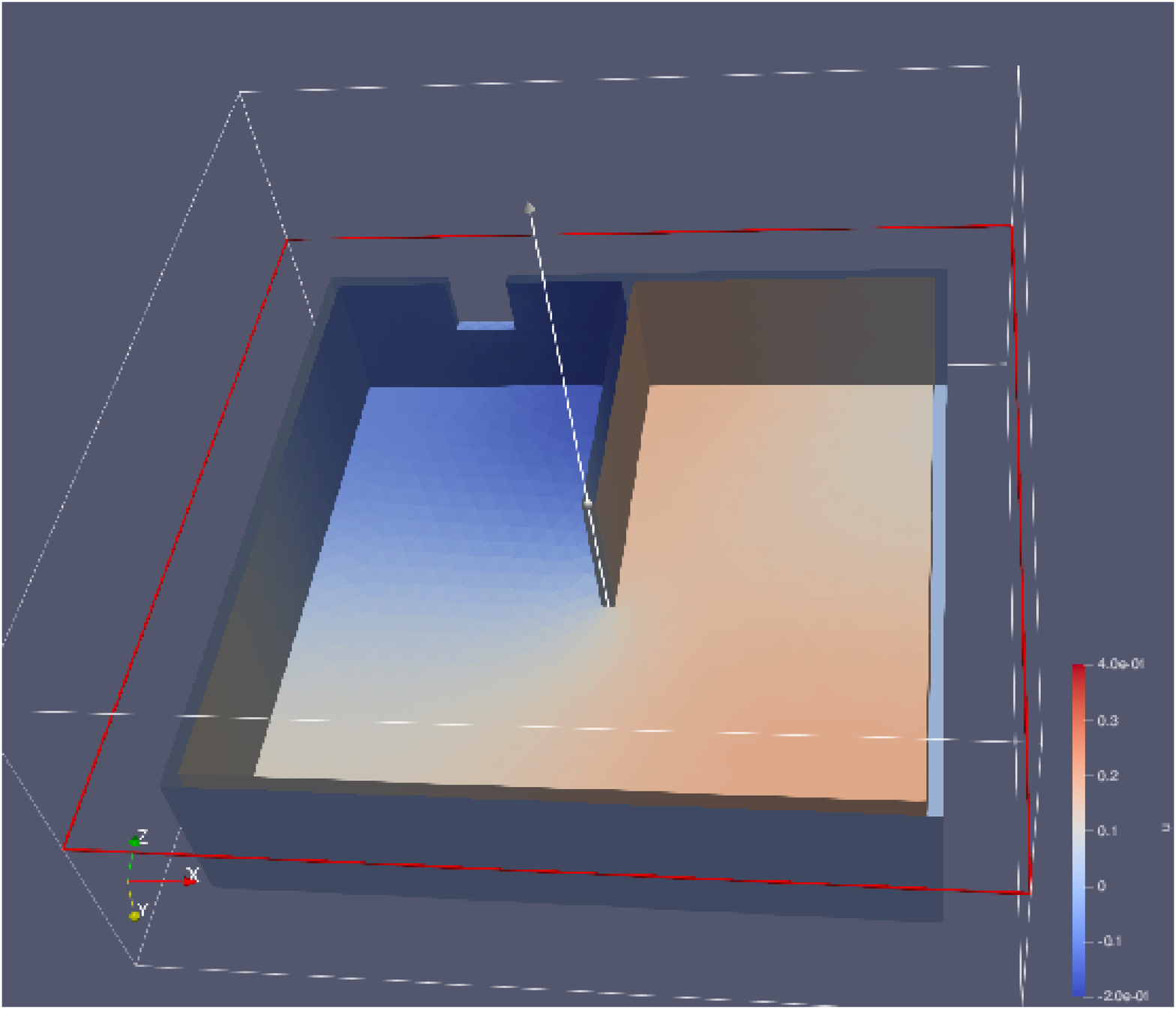}
    &\includegraphics[width=.3\textwidth]{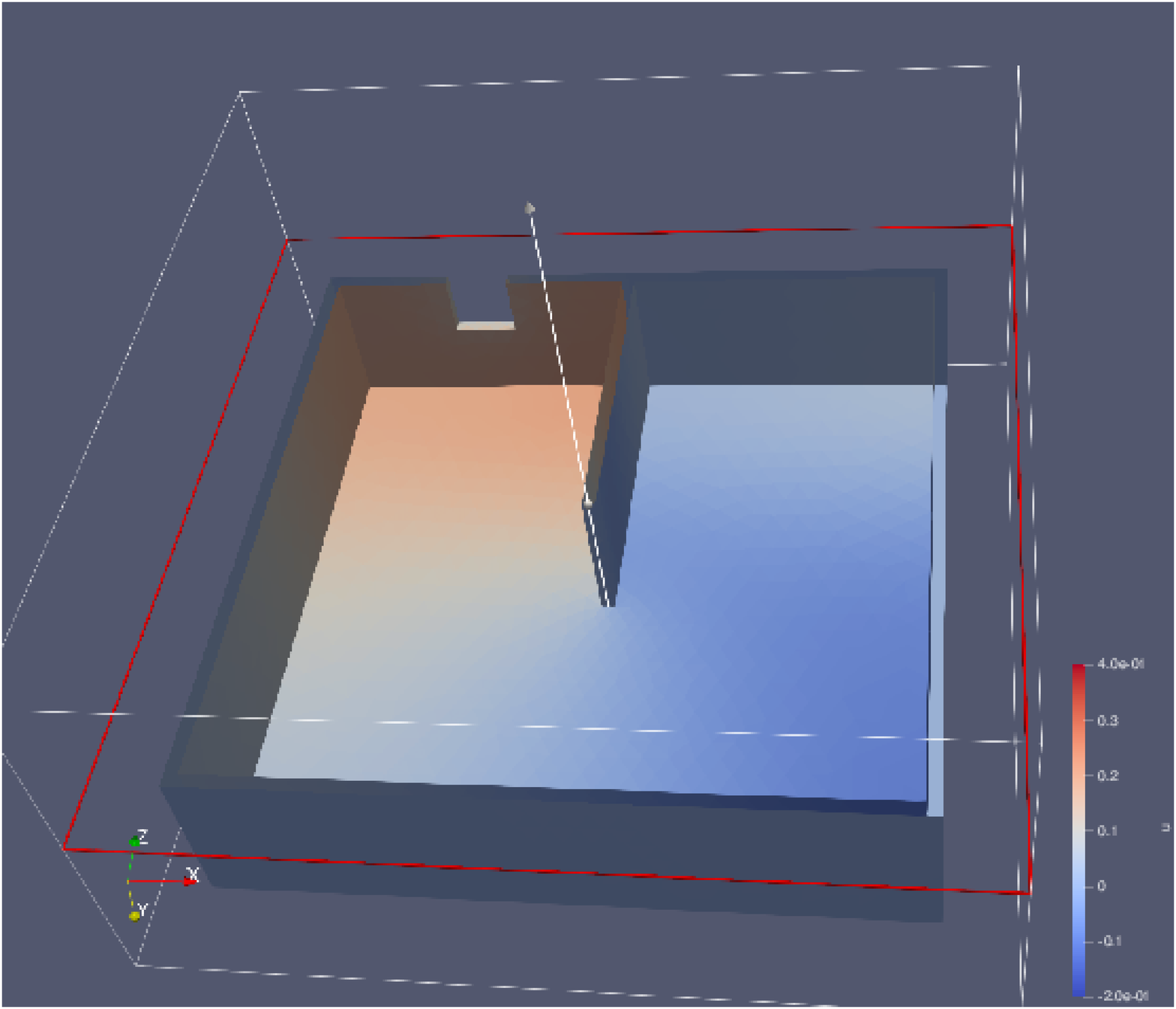}
    &\includegraphics[width=.3\textwidth]{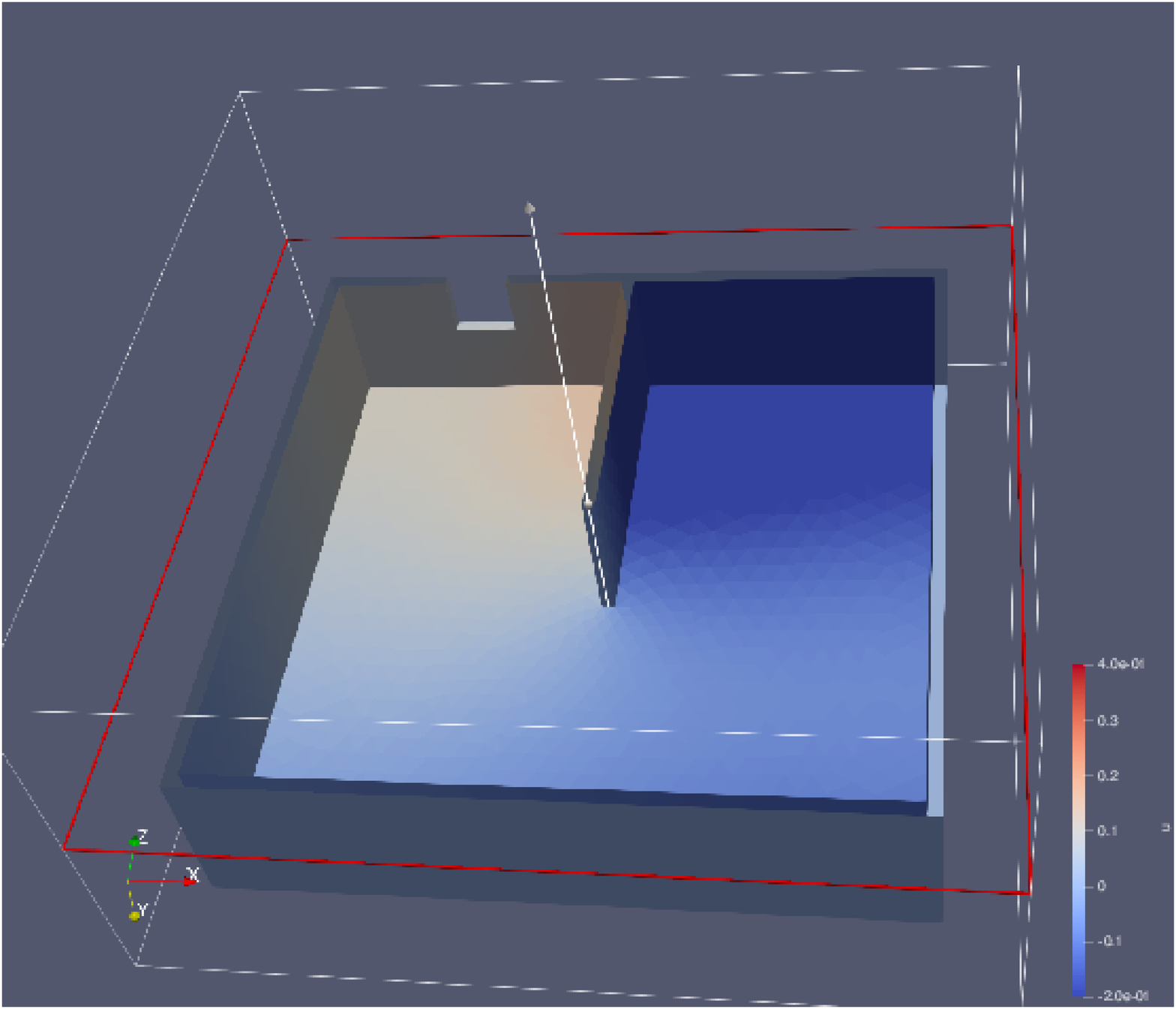}\\
    $t=6.0$ & $t=6.6$ & $t=7.2$
  \end{tabular}
  \caption{Snapshots of the distribution of the sound pressure $u$ on the surface of the hollow model in the \textbf{initial} shape. The value of $u$ over all the time sptes ranges from $\prval{-0.906227}$ to $\prval{1.794499}$ but is truncated from $-0.2$ to $0.4$, which correspond to blue and red, respectively, so that the pressure inside the cavity can be clearly seen.}
  \label{fig:bwe-snapshot0}
\end{figure}

\begin{figure}[H]     
  \centering          
  \begin{tabular}{ccc} 
    \includegraphics[width=.3\textwidth]{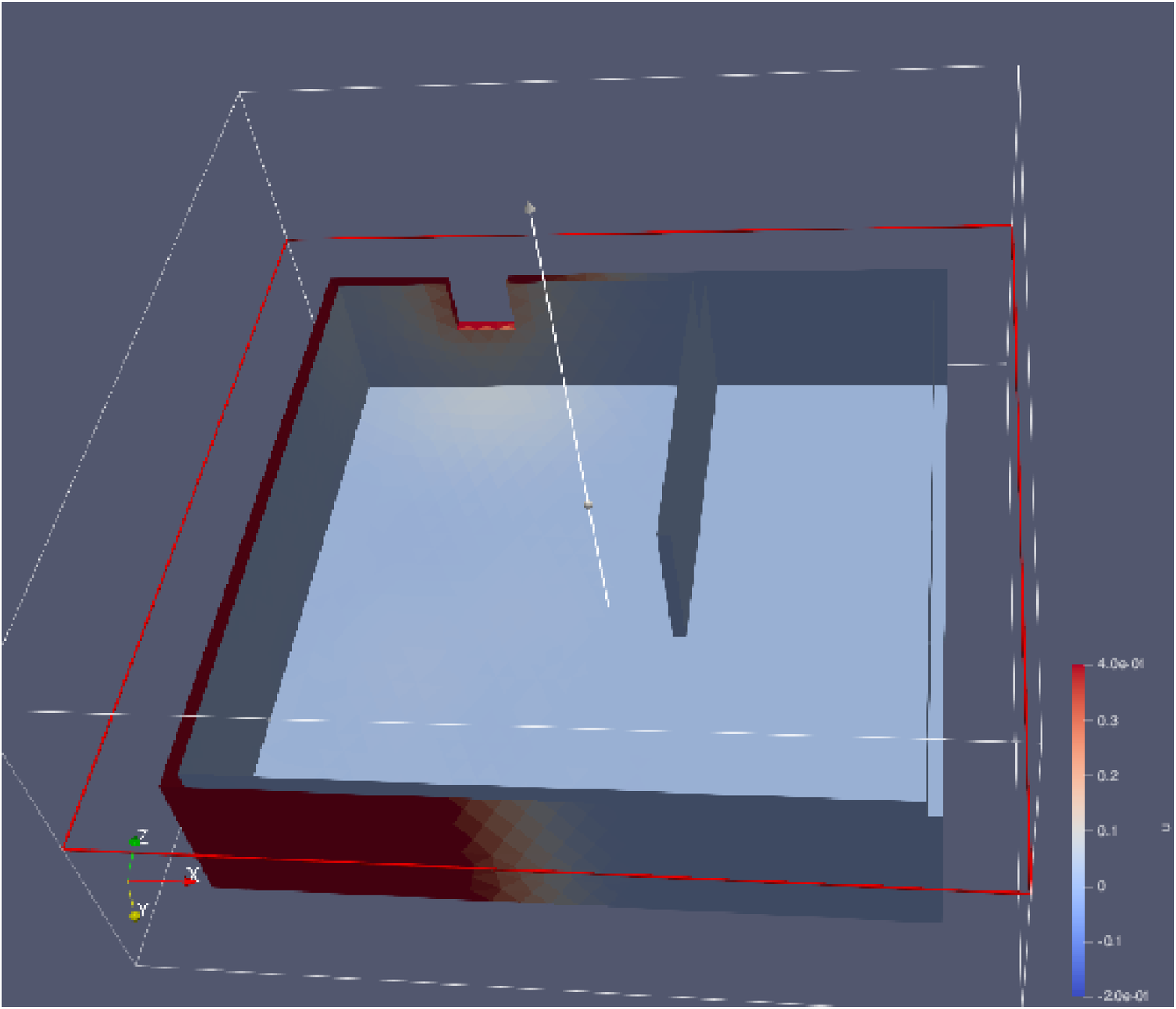}
    &\includegraphics[width=.3\textwidth]{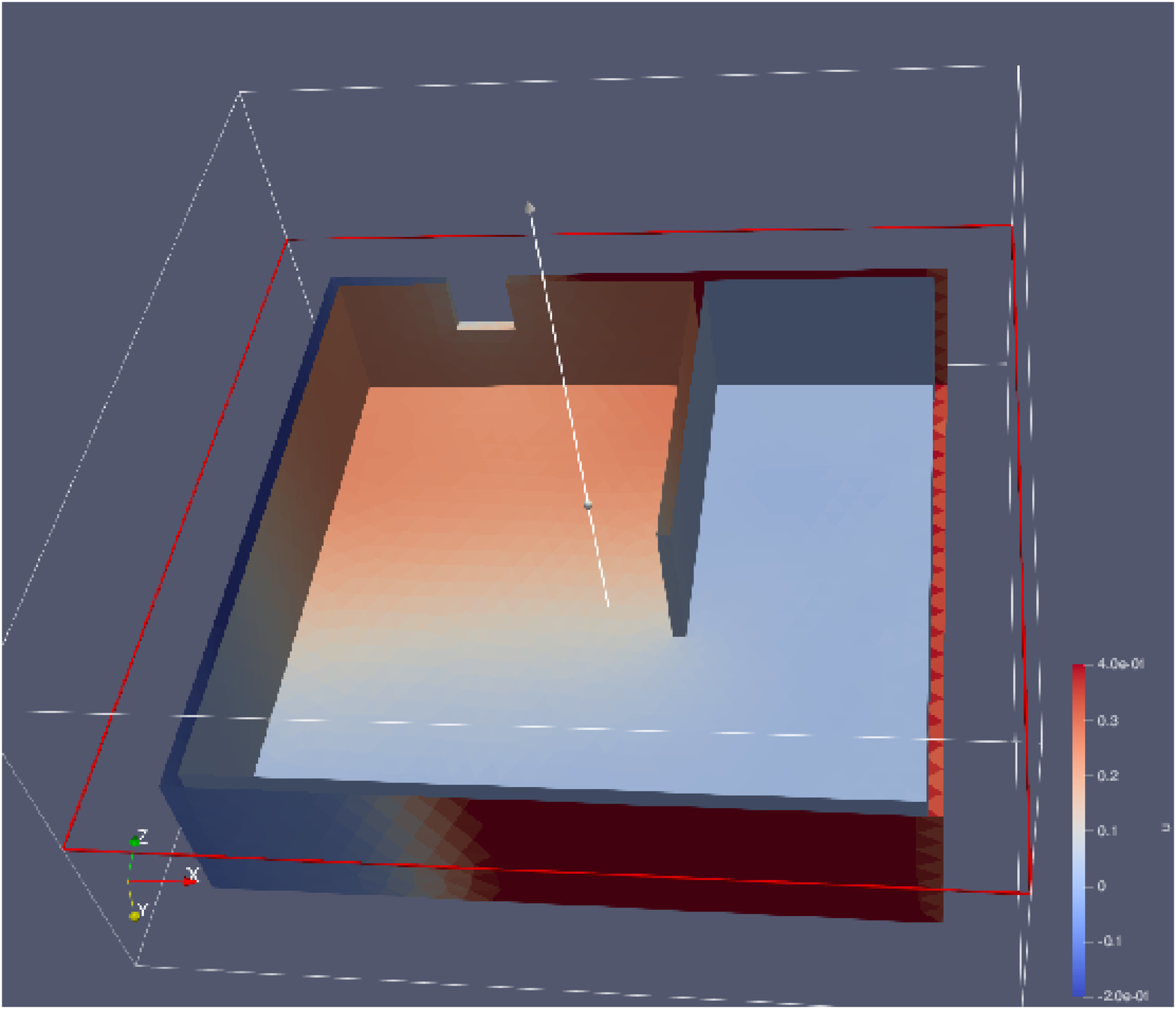}
    &\includegraphics[width=.3\textwidth]{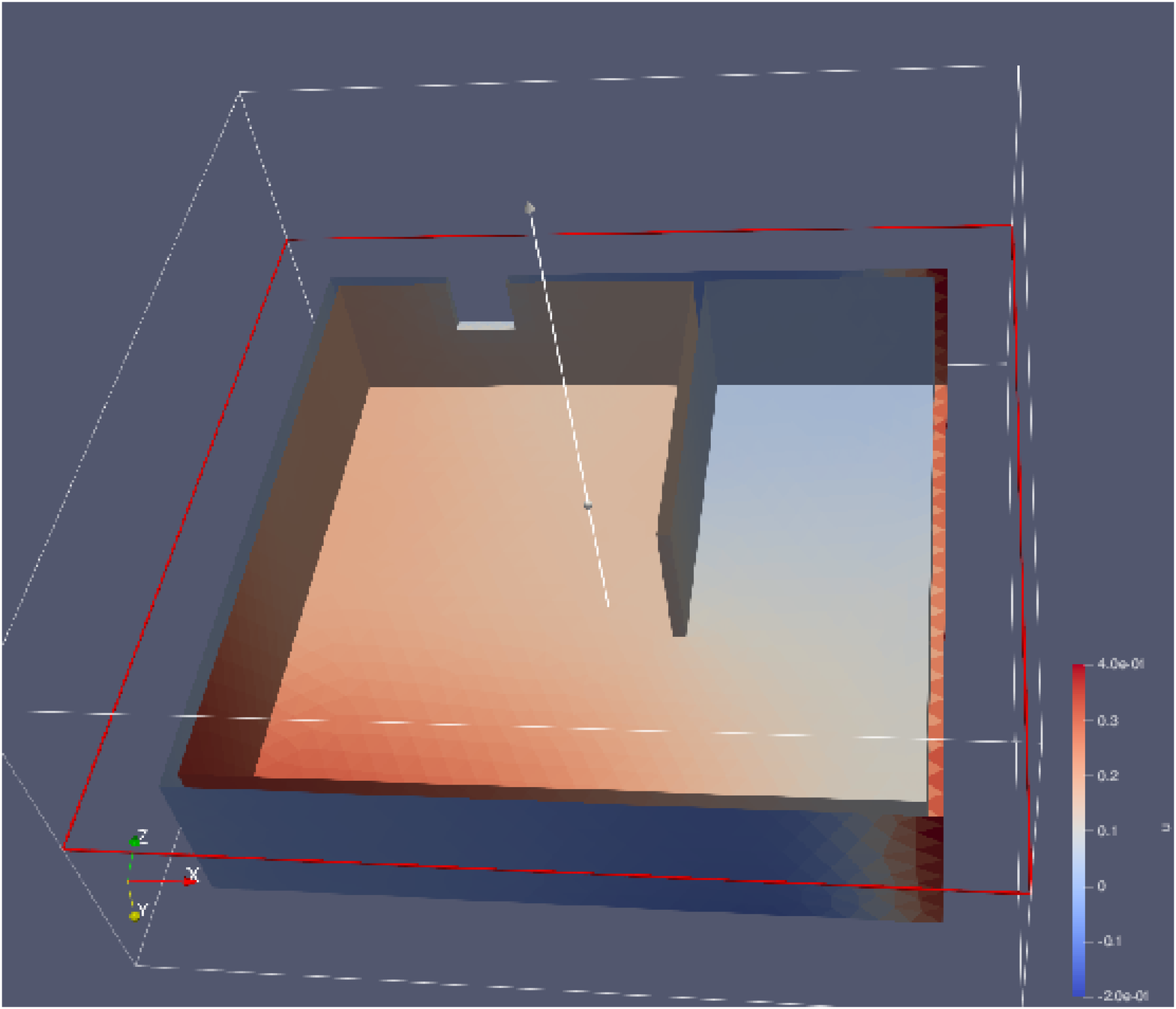}\\
    $t=0.6$ & $t=1.2$ & $t=1.8$\\
    \includegraphics[width=.3\textwidth]{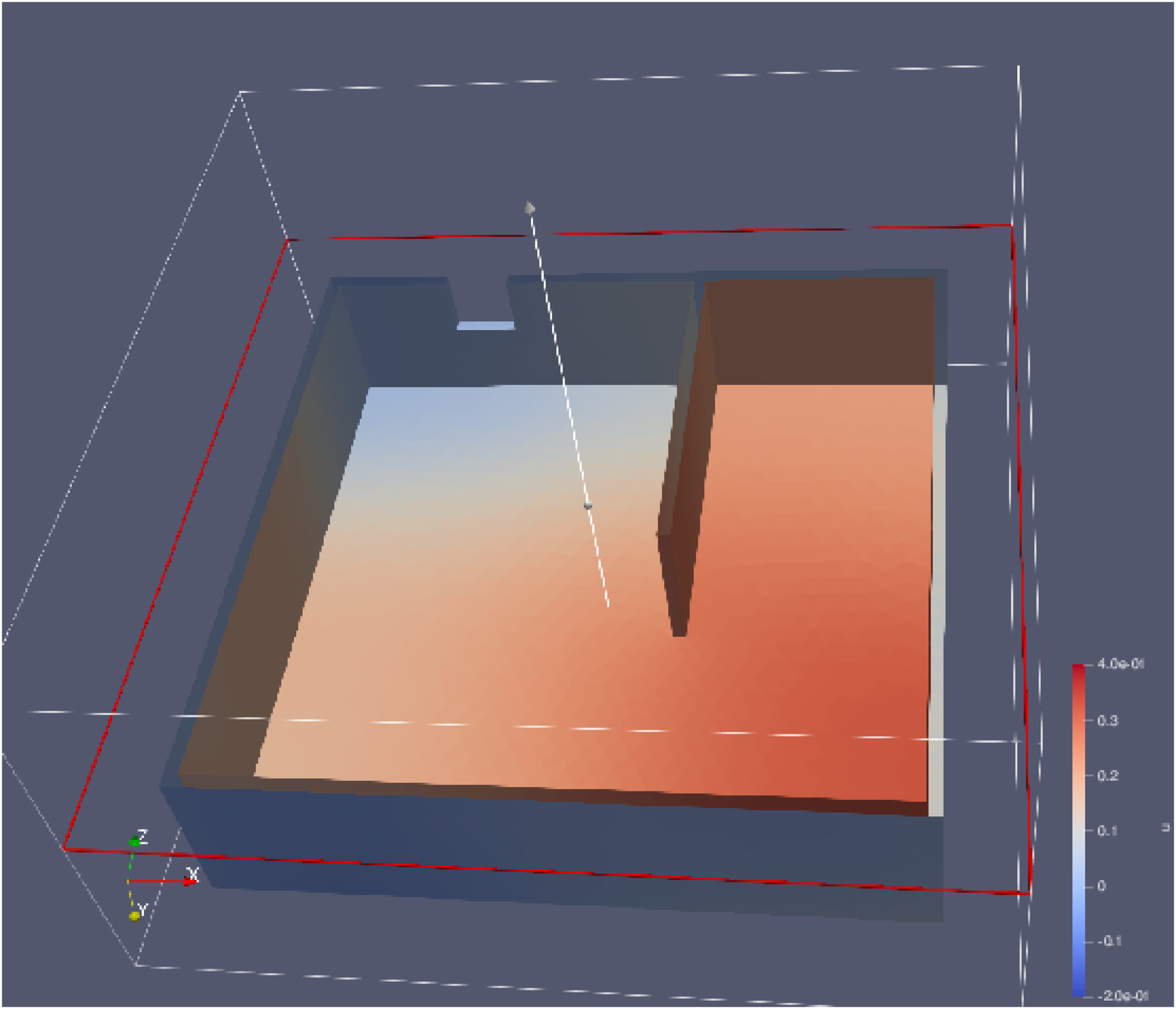}
    &\includegraphics[width=.3\textwidth]{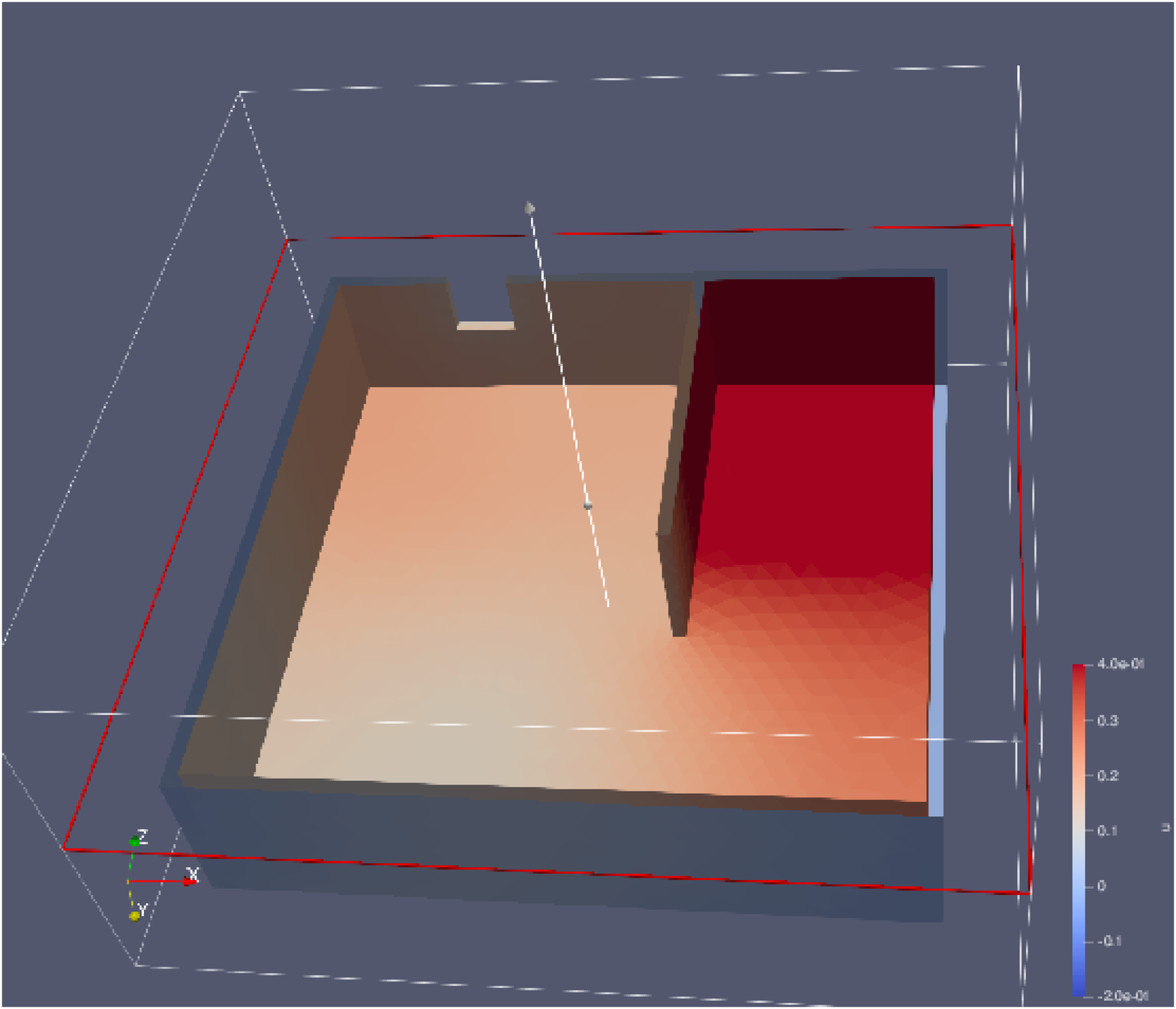}
    &\includegraphics[width=.3\textwidth]{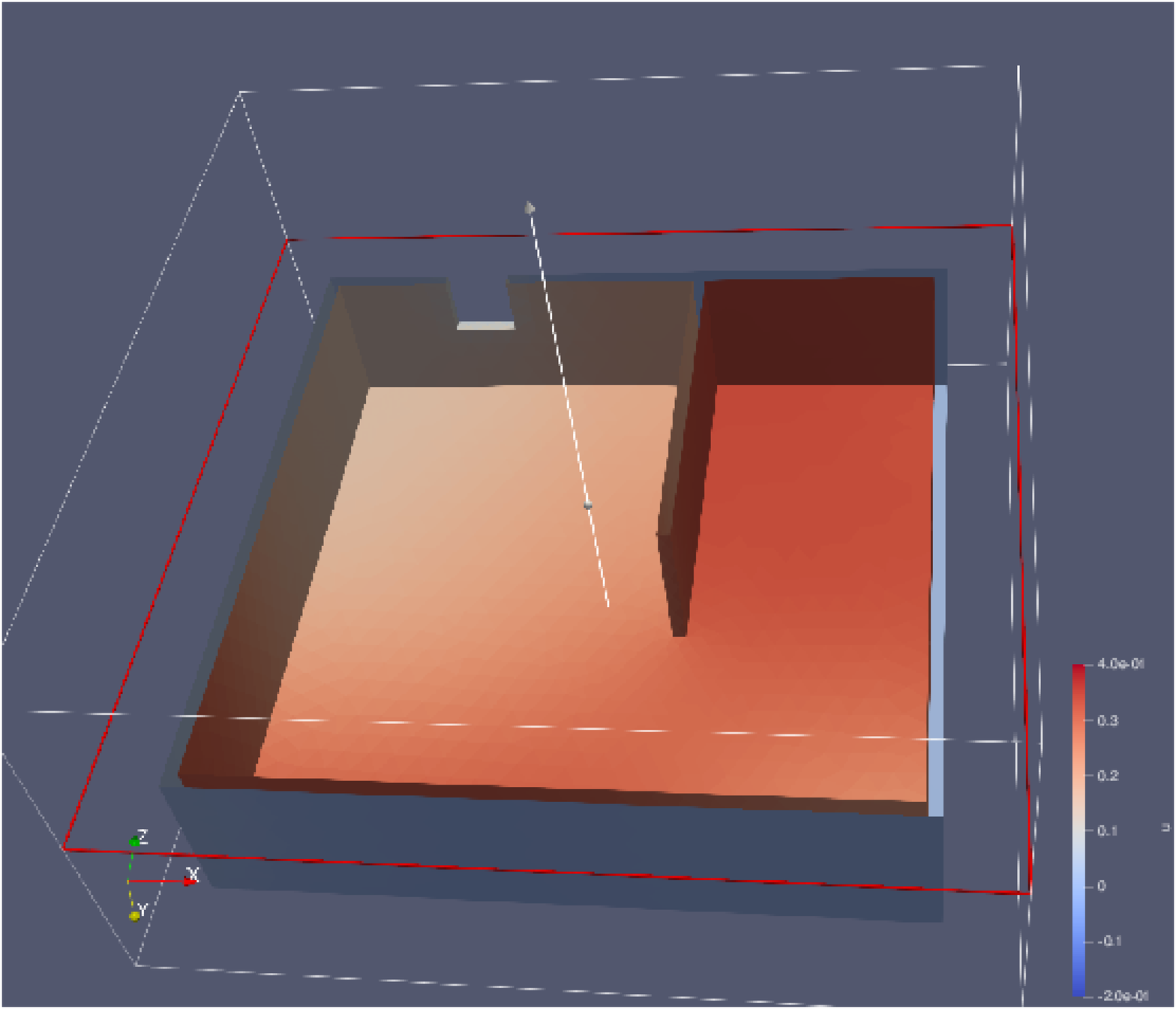}\\
    $t=2.4$ & $t=3.0$ & $t=3.6$\\
    \includegraphics[width=.3\textwidth]{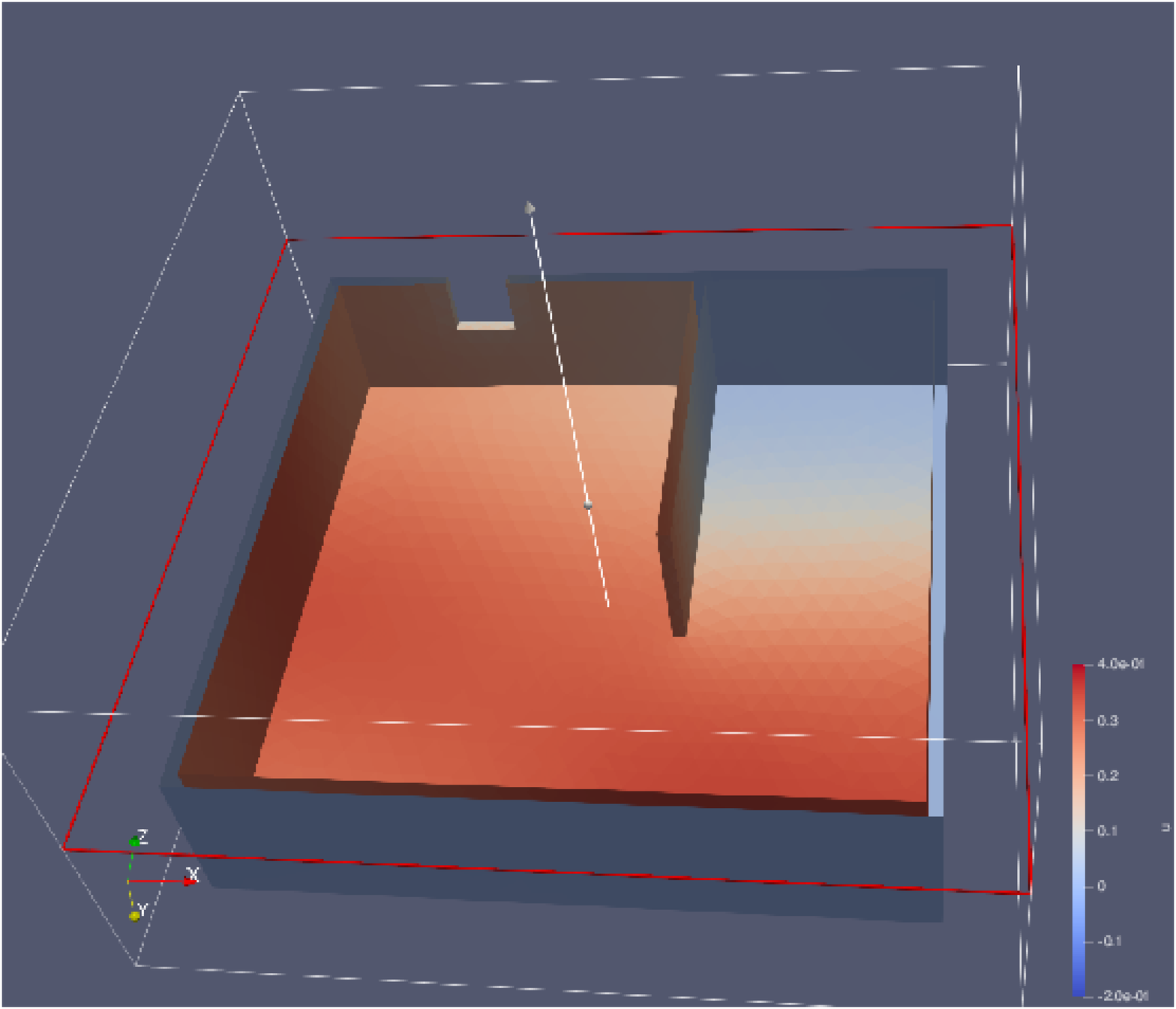}
    &\includegraphics[width=.3\textwidth]{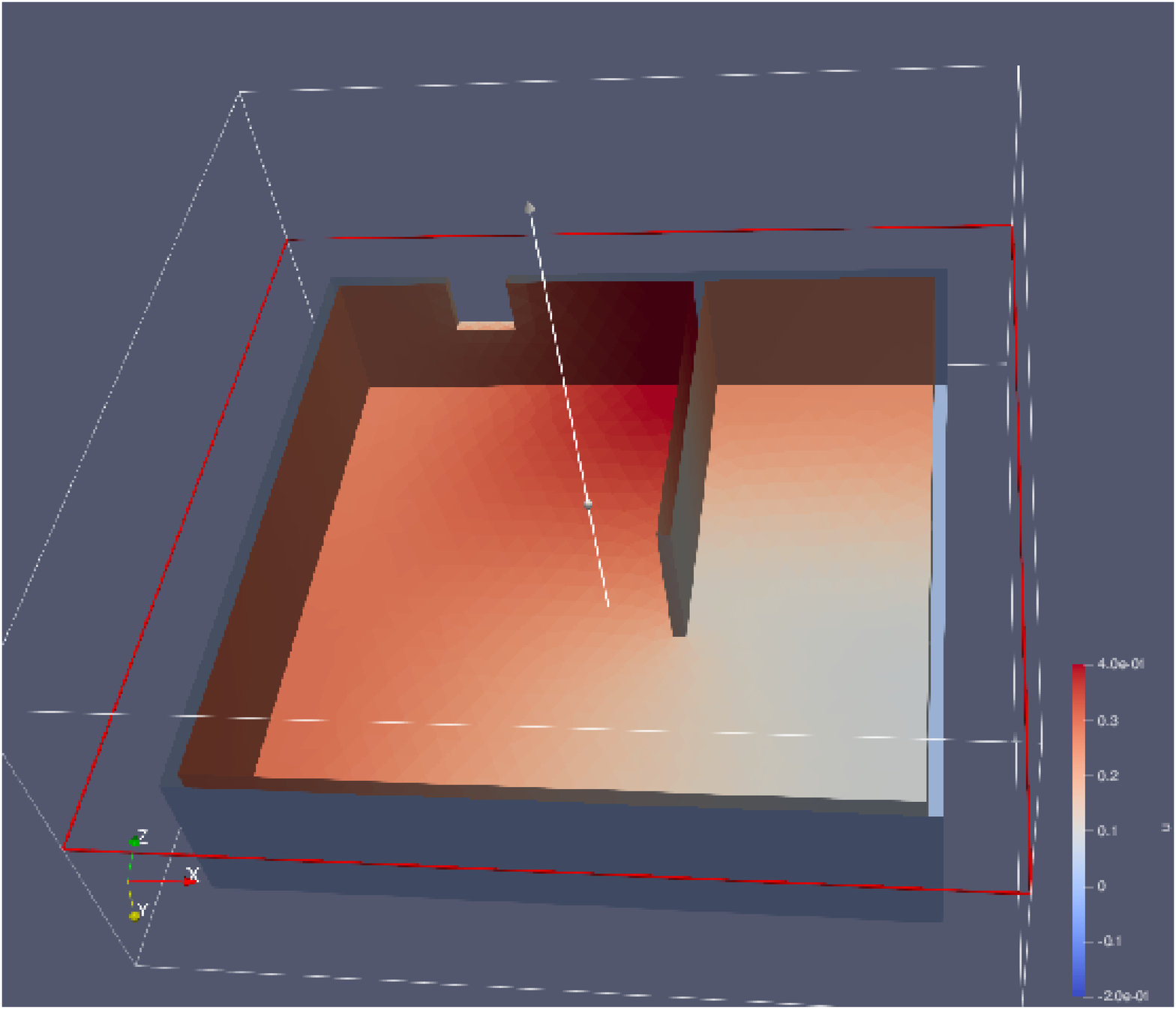}
    &\includegraphics[width=.3\textwidth]{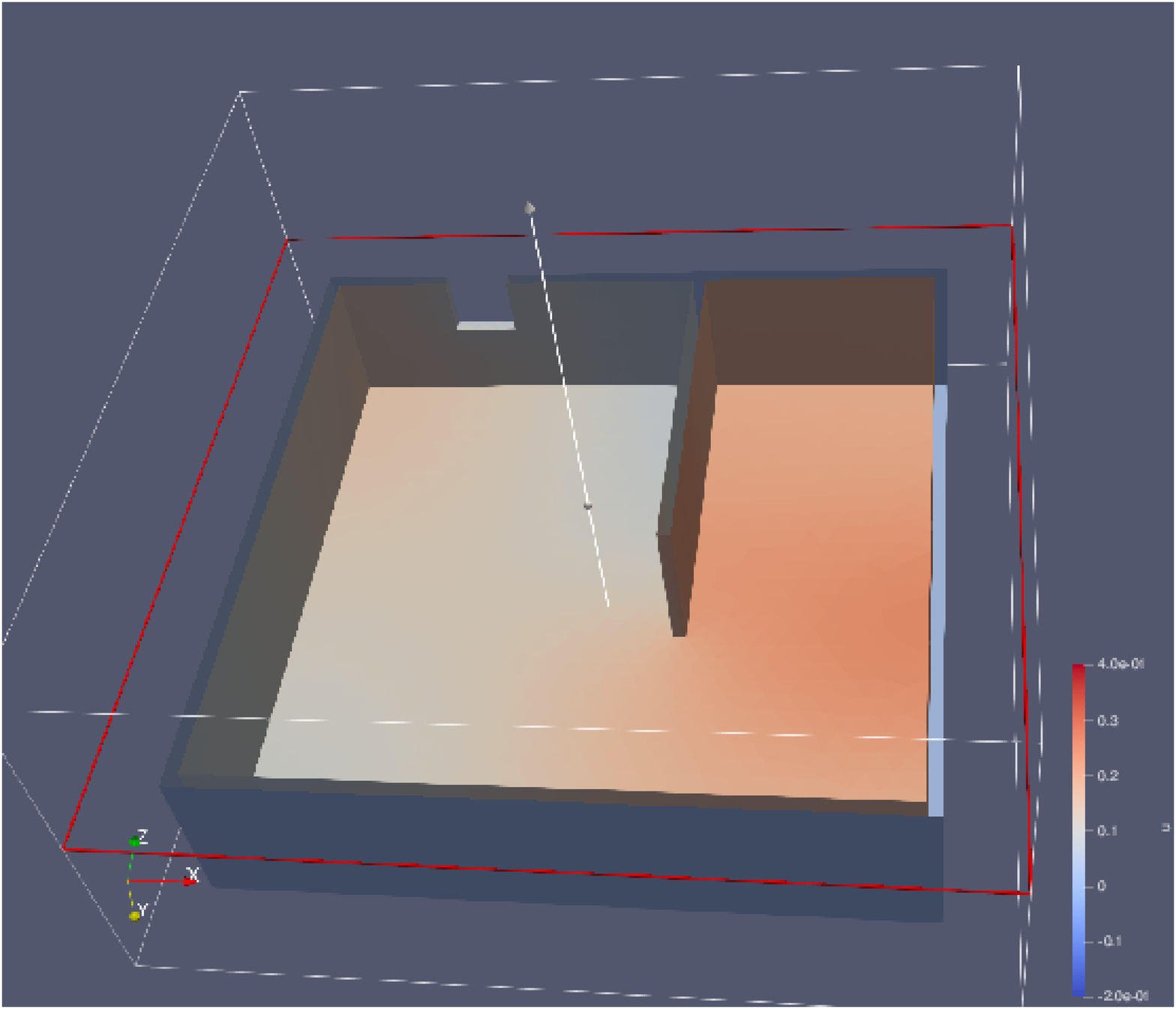}\\
    $t=4.2$ & $t=4.8$ & $t=5.4$\\
    \includegraphics[width=.3\textwidth]{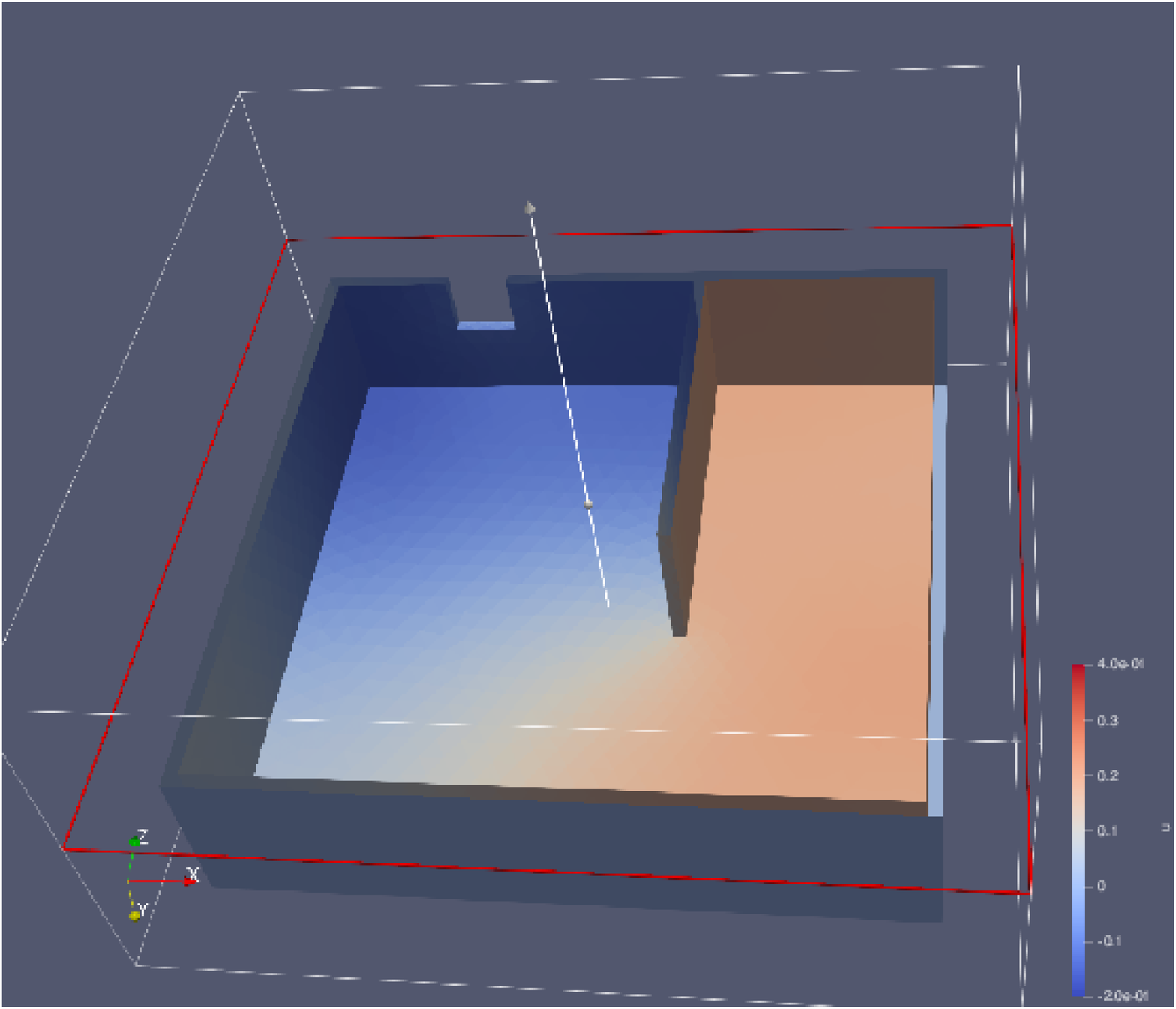}
    &\includegraphics[width=.3\textwidth]{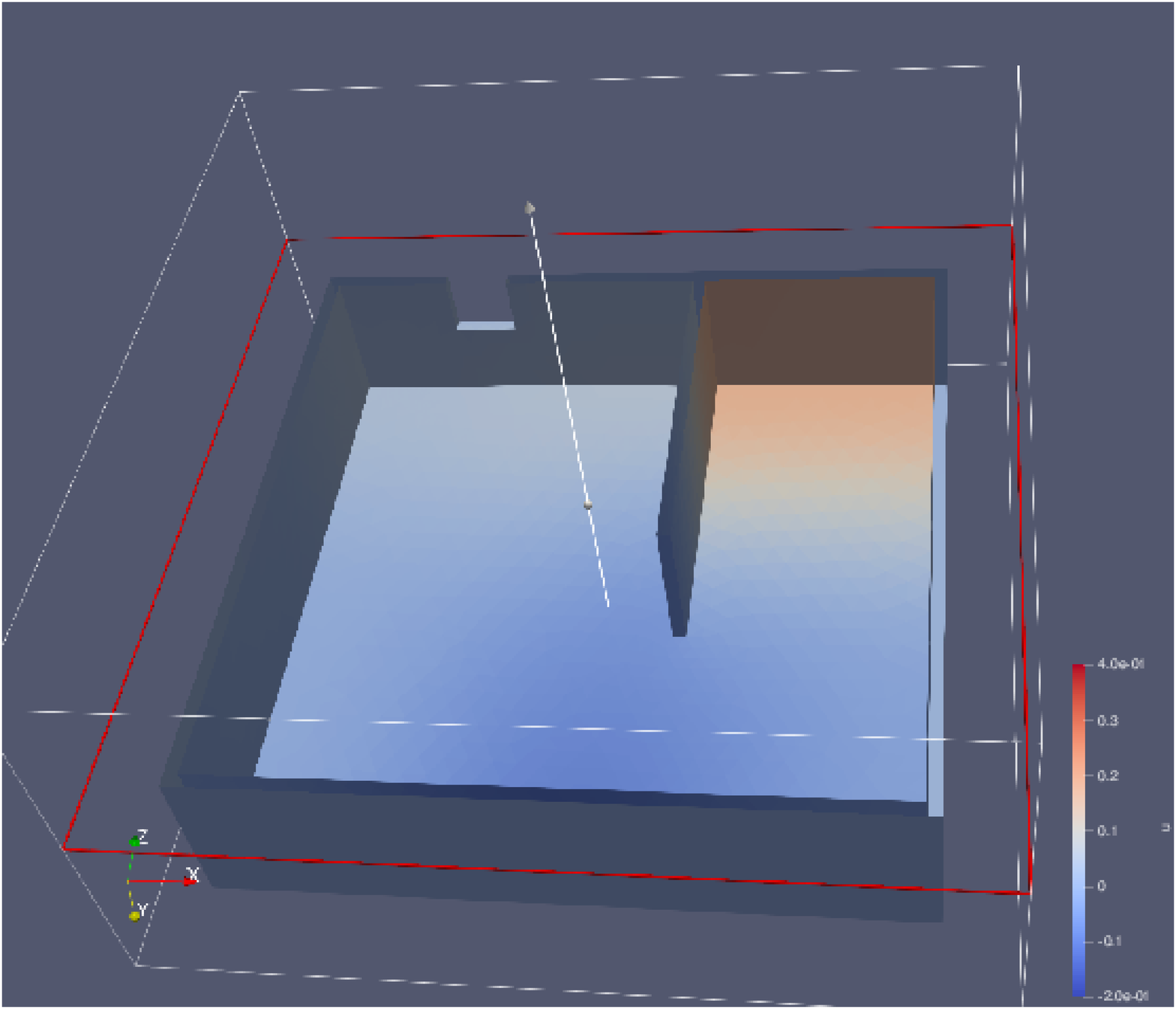}
    &\includegraphics[width=.3\textwidth]{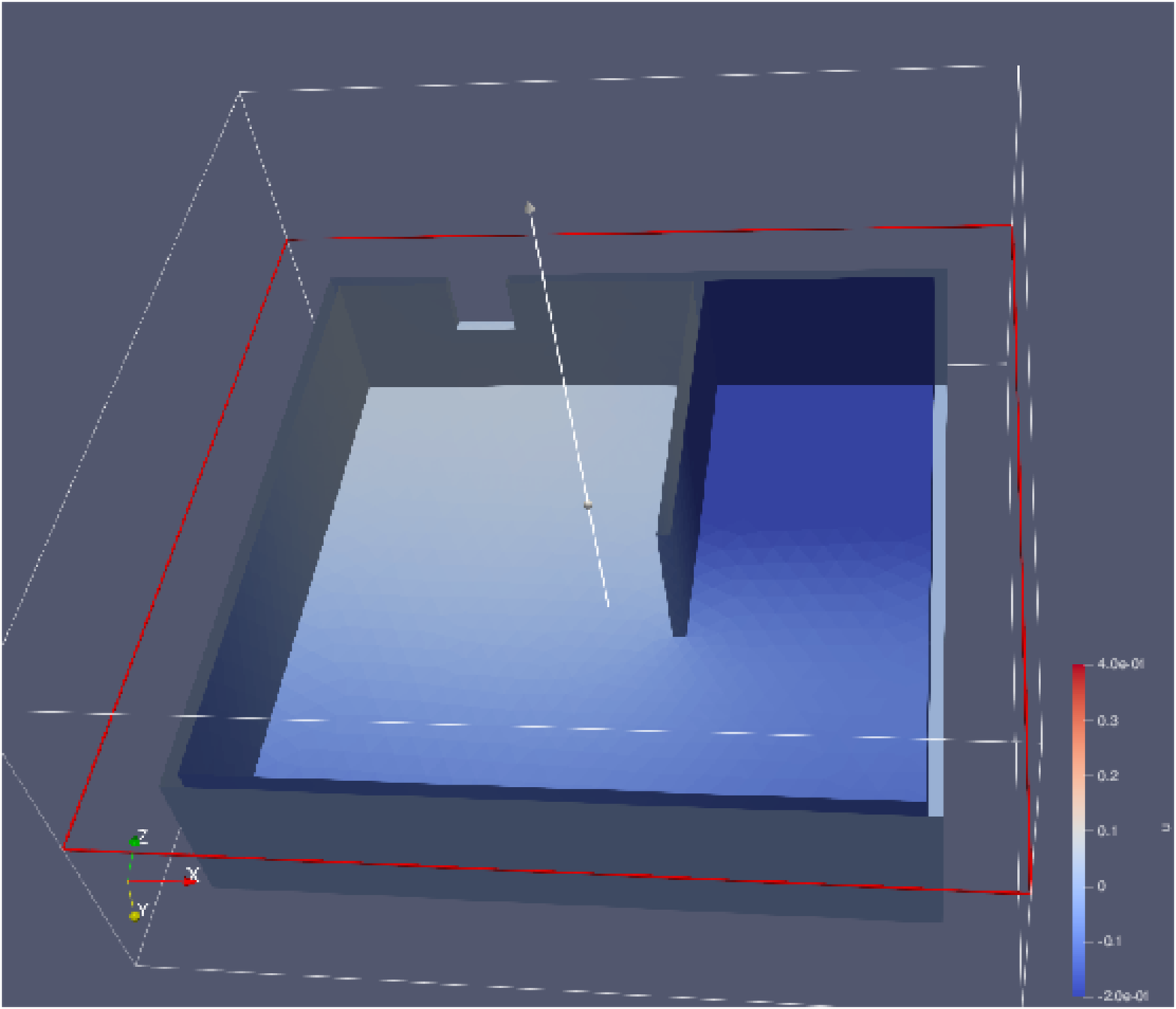}\\
    $t=6.0$ & $t=6.6$ & $t=7.2$
  \end{tabular}
  \caption{Snapshots of the distribution of the sound pressure $u$ on the surface of the hollow model in the \textbf{optimal} shape (at the iteration step $36$). The value of $u$ over all the time steps ranges from $\prval{-0.906107}$ to $\prval{1.794499}$ but is truncated from $-0.2$ to $0.4$, which correspond to blue and red, respectively, so that the pressure inside the cavity can be clearly seen.}
  \label{fig:bwe-snapshot}
\end{figure}

\section{Conclusion}\label{s:conclusion}

The present study enhanced the fast time-domain boundary element method (TDBEM) for the 3D scalar wave equation proposed in the previous study~\cite{takahashi2014}. First, we used the Burton--Miller-type boundary integral equation (BMBIE) instead of the ordinary boundary integral equation (OBIE) to stabilise the TDBEM, following previous studies~\cite{fukuhara2019,chiyoda2019}. Second, we formulated the TDBEM that employs the B-spline function of order $d$ ($\ge 1$) as the temporal basis, whereas only the special case of $d=1$, which corresponds to the piecewise-linear temporal basis, was used in \cite{takahashi2014}. Third, we generalised the interpolation-based FMM from $d=1$ to $d\ge 2$ for both the OBIE and the BMBIE corresponding to the generalisation of the B-spline temporal basis of order $d$ from $d=1$ to $d\ge 2$. In particular, we constructed an $O(\Nt)$ M2L by considering the auxiliary local coefficient as well as the derivatives of local coefficient, and deriving their recurrence formulae in Formula~\ref{formula:m2l}.

We assessed the enhanced (fast) TDBEM through numerical experiments in Section~\ref{s:num} as well as Section~\ref{s:app}, where we considered the typical boundary conditions in time-domain acoustics, that is, $u=0$ (i.e. acoustically soft) and $q=0$ (i.e. acoustically hard). The results indicate the following:
\begin{enumerate}

\item In the case of $u=0$, the OBIE is unstable but the BMBIE is stable, which is consistent with the semi-analytic study on the stability of TDBEM by Fukuhara et al.~\cite{fukuhara2019}. This is true for $d=1$ but not for $d\ge 2$. (We actually considered $d=1$, $2$, and $3$.) The instability is irrelevant to the acceleration by the interpolation-based FMM. 

\item In the case of $q=0$, the BMBIE with $d=2$ is stable and more accurate than the OBIE with $d=1$ when it is available, as seen in Section~\ref{s:sphere}. The OBIE with $d=1$, which was considered in the previous study~\cite{takahashi2014}, can be unstable when the boundary shape is complicated, as shown in Section~\ref{s:hollow}.

\item For $d=3$, the TDBEM was always unstable in any case.

\end{enumerate}

We discussed the instability due to the choice of $d$ from the viewpoint of cancellation of significant digits in calculating the layer potentials of the OBIE and BMBIE, but a more rigorous analysis is necessary in the future. Nevertheless, the present work opens up new avenues for analysing and designing 3D large-scale time-domain exterior acoustic problems stably and efficiently.

Future plans include enhancing the present TDBEM for acoustics to electromagnetics. Because the combined field integral equation, which is known to be stable, contains the second-order time derivative,\footnote{This is the case when we consider a vector $\bm{e}$ such as $\dot{\bm{e}}=\bm{J}$, where $\bm{J}$ denotes the surface-induced current, to remove the time integral, which can prevent the construction of an efficient algorithm with respect to time, in the scalar potential of the electric field integral equation~\cite{jung2003}.} we need a smooth temporal basis such as $d\ge 2$ in the case of the B-spline basis of order $d$. Therefore, the present study is important as it lays a foundation for the electromagnetic TDBEM under consideration.

\appendix

\section{Derivation of the discretised OBIE in (\ref{eq:bie_disc3})}\label{s:simplify}

We derive (\ref{eq:bie_disc3}) by rewriting the RHS in (\ref{eq:bie_disc}), where $\beta=0$ is replaced with $\beta=\beta^*$; that is,
\begin{eqnarray*}
  R:=\sum_{\kappa=0}^{d+1} \sum_{\beta=\beta^*}^{\alpha-1} w^{\kappa,d} \left(\mat{U}^{(\alpha-\beta-\kappa)} \mat{q}^\beta - \mat{W}^{(\alpha-\beta-\kappa)} \mat{u}^\beta \right).
\end{eqnarray*}
To this end, we first split the summation over $\beta$ in $R$ into three parts after introducing a new index $\beta':=\beta+\kappa$ as follows:
\begin{eqnarray*}
 R &=& \sum_{\kappa=0}^{d+1} \sum_{\beta'=\beta^*+\kappa}^{\alpha-1+\kappa} w^{\kappa,d} \left(\mat{U}^{(\alpha-\beta')} \mat{q}^{\beta'-\kappa} - \mat{W}^{(\alpha-\beta')} \mat{u}^{\beta'-\kappa} \right)\nonumber\\
   &=& \sum_{\kappa=0}^{d+1} \left( \sum_{\beta'=\beta^*}^{\alpha-1}-\sum_{\beta'=\beta^*}^{\beta^*+\kappa-1}+\sum_{\beta'=\alpha}^{\alpha-1+\kappa}\right) w^{\kappa,d} \left(\mat{U}^{(\alpha-\beta')} \mat{q}^{\beta'-\kappa} - \mat{W}^{(\alpha-\beta')} \mat{u}^{\beta'-\kappa} \right),
\end{eqnarray*}
where we ignore a summation $\sum_{i=s}^{e}$ if $s>e$. Then, the third summation $\sum_{\beta'=\alpha}^{\alpha-1+\kappa}$ always vanishes for the following reason:
\begin{itemize}
\item For $\kappa=0$, the third summation reduces to $\sum_{\beta'=\beta^*}^{\beta^*-1}$ and, thus, vanishes in accordance with the above convention. (In this case, the second summation also vanishes, whereas the first summation is identical to the original summation, i.e. $\sum_{\beta=\beta^*}^{\alpha-1}$.)

\item For $\kappa \ge 1$, the inequality $\alpha-\beta'\le 0$ holds for any $\beta'\in[\alpha,\alpha-1+\kappa]$. Then, because $\mat{U}^{(\alpha-\beta')}$ and $\mat{W}^{(\alpha-\beta')}$ are zero from Remark~\ref{remark:non_positive_index}, the third summation vanishes.
\end{itemize}

Therefore, we have
\begin{eqnarray}
  R
  &=& \sum_{\beta'=\beta^*}^{\alpha-1} \sum_{\kappa=0}^{d+1} w^{\kappa,d} \left(\mat{U}^{(\alpha-\beta')} \mat{q}^{\beta'-\kappa} - \mat{W}^{(\alpha-\beta')} \mat{u}^{\beta'-\kappa} \right)
  - \underbrace{\sum_{\kappa=0}^{d+1} \sum_{\beta'=\beta^*}^{\beta^*+\kappa-1} w^{\kappa,d} \left(\mat{U}^{(\alpha-\beta')} \mat{q}^{\beta'-\kappa} - \mat{W}^{(\alpha-\beta')} \mat{u}^{\beta'-\kappa} \right)}_{\displaystyle F}\nonumber\\
  &=& \left(\underbrace{\sum_{\beta'=\beta^*}^{\beta^*+d}}_{\displaystyle G}+\underbrace{\sum_{\beta'=\beta^*+d+1}^{\alpha-1}}_{\displaystyle 
H}\right) \sum_{\kappa=0}^{d+1} w^{\kappa,d} \left(\mat{U}^{(\alpha-\beta')} \mat{q}^{\beta'-\kappa} - \mat{W}^{(\alpha-\beta')} \mat{u}^{\beta'-\kappa} \right) - F,
\end{eqnarray}
where the summation over $\beta'$ was split into two parts, i.e. $G$ and $H$. Because $\sum_{\kappa=0}^{d+1}\sum_{\beta'=\beta^*}^{\beta^*+\kappa-1}=\sum_{\beta'=\beta^*}^{\beta^*+d} \sum_{\kappa=\beta'-\beta^*+1}^{d+1}$ holds in $F$, we have
\begin{eqnarray*}
  G-F
  &=&\sum_{\beta'=\beta^*}^{\beta^*+d}\sum_{\kappa=0}^{d+1} w^{\kappa,d} \left(\mat{U}^{(\alpha-\beta')} \mat{q}^{\beta'-\kappa} - \mat{W}^{(\alpha-\beta')} \mat{u}^{\beta'-\kappa} \right)
  -\sum_{\beta'=\beta^*}^{\beta^*+d} \sum_{\kappa=\beta'-\beta^*+1}^{d+1} w^{\kappa,d} \left(\mat{U}^{(\alpha-\beta')} \mat{q}^{\beta'-\kappa} - \mat{W}^{(\alpha-\beta')} \mat{u}^{\beta'-\kappa} \right)\nonumber\\
  &=&\sum_{\beta'=\beta^*}^{\beta^*+d} \left[\mat{U}^{(\alpha-\beta')}\left(\sum_{\kappa=0}^{d+1}w^{\kappa,d}\mat{q}^{\beta'-\kappa}-\sum_{\kappa=\beta'-\beta^*+1}^{d+1}w^{\kappa,d}\mat{q}^{\beta'-\kappa}\right)
  -\mat{W}^{(\alpha-\beta')}\left(\sum_{\kappa=0}^{d+1}w^{\kappa,d}\mat{u}^{\beta'-\kappa}-\sum_{\kappa=\beta'-\beta^*+1}^{d+1}w^{\kappa,d}\mat{u}^{\beta'-\kappa}\right)\right]\nonumber\\
  &=&\sum_{\beta'=\beta^*}^{\beta^*+d} \left (\mat{U}^{(\alpha-\beta')}\sum_{\kappa=0}^{\beta'-\beta^*} w^{\kappa,d}\mat{q}^{\beta'-\kappa} - \mat{W}^{(\alpha-\beta')}\sum_{\kappa=0}^{\beta'-\beta^*} w^{\kappa,d}\mat{u}^{\beta'-\kappa}\right).
\end{eqnarray*}
Meanwhile, we can rewrite $H$ as 
\begin{eqnarray*}
  H&=&\sum_{\beta'=\beta^*+d+1}^{\alpha-1}\sum_{\kappa=0}^{d+1} w^{\kappa,d} \left(\mat{U}^{(\alpha-\beta')} \mat{q}^{\beta'-\kappa} - \mat{W}^{(\alpha-\beta')} \mat{u}^{\beta'-\kappa} \right)\nonumber\\
   &=&\sum_{\beta'=\beta^*+d+1}^{\alpha-1}\left(\mat{U}^{(\alpha-\beta')} \sum_{\kappa=0}^{d+1} w^{\kappa,d} \mat{q}^{\beta'-\kappa} - \mat{W}^{(\alpha-\beta')} \sum_{\kappa=0}^{d+1} w^{\kappa,d} \mat{u}^{\beta'-\kappa} \right).
\end{eqnarray*}
Combining these, we have
\begin{eqnarray}
  R=G-F+H=\sum_{\beta=\beta^*}^{\alpha-1}\left(\mat{U}^{(\alpha-\beta)} \bm{\uptau}^{\beta} - \mat{W}^{(\alpha-\beta)} \bm{\upsigma}^{\beta}\right),
  \label{eq:RHS}
\end{eqnarray}
where we define the following boundary variable $\bm{\uptau}^\beta$ in terms of $\mat{q}^{\beta-\kappa}$:
\begin{eqnarray}
  \bm{\uptau}^{\beta}
  &:=&\begin{cases}
    \displaystyle\sum_{\kappa=0}^{\beta-\beta^*}w^{\kappa,d}\mat{q}^{\beta-\kappa} & \text{for $\beta^*\le\beta\le\beta^*+d$}\\
    \displaystyle\sum_{\kappa=0}^{d+1}w^{\kappa,d}\mat{q}^{\beta-\kappa} & \text{for $\beta^*+d+1\le\beta\le\alpha-1$}
  \end{cases}\nonumber\\
  &=&\sum_{\kappa=0}^{\min(d+1,\beta-\beta^*)}w^{\kappa,d}\mat{q}^{\beta-\kappa}\quad\text{for $\beta^*\le\beta\le\alpha-1$}.
     \label{eq:tau}
\end{eqnarray}
Similarly, we define the boundary variable $\bm{\sigma}^{\beta}$ as follows:
\begin{eqnarray}
  \bm{\upsigma}^{\beta}:=\sum_{\kappa=0}^{\min(d+1,\beta-\beta^*)}w^{\kappa,d}\mat{u}^{\beta-\kappa}\quad\text{for $\beta^*\le\beta\le\alpha-1$}.
  \label{eq:sigma}
\end{eqnarray}

Finally, in accordance with the fact that the index $\alpha$ starts from one (not zero), we increase the index $\beta$ by one in (\ref{eq:RHS}) to yield the simplified OBIE in (\ref{eq:bie_disc3}).

\section{Evaluation of the spatial integrals in $U_{ij}^{(\gamma)}$ and $W_{ij}^{(\gamma)}$}\label{s:element_integral}

First, we consider a coefficient
\begin{eqnarray}
  U_{ij}^{(\gamma)}:=\int_{E_j}\frac{(ct_{\gamma}-\abs{\bm{x}_i-\bm{y}})_+^d}{4\pi(c\Dt)^d\abs{\bm{x}_i-\bm{y}}} \diff S_y,
\end{eqnarray}
that is, the single-layer potential with respect to a collocation point $\bm{x}_i$ (simply denoted by P) and a triangular element $E_j$ (denoted by ABC). We let O be the foot of the perpendicular from P to the plane including ABC. Then, the integral over ABC can be evaluated as the summation of three integrals over sub-triangles OAB (denoted by $E_j^1$), OBC ($E_j^2$), and OCA ($E_j^3$).

In a sub-triangle $E_j^n$, we denote its (local) vertexes by O, V$_1$, and V$_2$. Next, we introduce a local Cartesian coordinate system $xyz$, where the positive $z$-direction is chosen as the direction of the normal vector to $E_j$ and the positive $x$-direction is from V$_1$ to V$_2$. In this system, we denote the $x$-coordinates of the vectors $\overrightarrow{\textrm{OV}_1}$ and $\overrightarrow{\textrm{OV}_2}$ by $x_1$ and $x_2$, respectively. Further, $y$ and $z$ denote the $y$- and $z$-coordinates common to both vectors. In addition, we introduce the polar coordinates $\rho$ and $\theta$ with O as the centre. We define the angles $\theta:=\tan^{-1}(y/x)$ and  $\theta_i:=\tan^{-1}(y/x_i)$ ($i=1,2$).

\begin{figure}[hbt]
  \begin{center}
    \includegraphics[width=.4\textwidth]{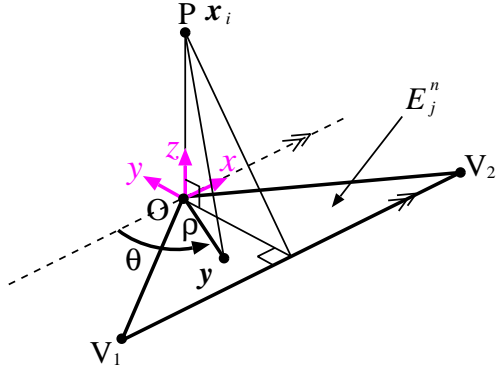}
    \caption{Symbols associated with the integration over a sub-triangle $E_j^n$.}
    \label{fig:zahyou}
  \end{center}
\end{figure}

Then, using the notations $s:=ct_{\gamma}$ (constant) and $r:=\abs{\bm{x}_i-\bm{y}}$, we can express $U_{ij}^{(\gamma)}$ as
\begin{eqnarray*}
  U_{ij}^{(\gamma)}
       =\frac{1}{4\pi(c\Dt)^d}\sum_{n=1}^3\int_{E_j^n}\frac{(s-r)_+^d}{r} \diff S_y
       =\frac{1}{4\pi(c\Dt)^d}\sum_{n=1}^3\int_{\theta_1}^{\theta_2}\int_0^\rho\frac{(s-\sqrt{\rho^2+z^2})_+^d}{r}\rho \diff\rho \diff\theta.
\end{eqnarray*}
By performing the integral with respect to $\rho$, we can obtain
\begin{eqnarray}
  U_{ij}^{(\gamma)}
  = \frac{1}{4\pi(c\Dt)^d}\left[\sum_{n=1}^3 I^d_n(\bm{x}_i,s) + I^d_0(\bm{x}_i,s)\right],
  \label{eq:Uij_tmp}
\end{eqnarray}
where the second term $I^d_0$  can be evaluated as follows:
\begin{eqnarray*}
  I^d_0
  := \sum_{n=1}^3\int_{\theta_1}^{\theta_2} \frac{(s-|z|)_+^{d+1}}{d+1}\diff\theta
  =\begin{cases}
    \displaystyle\frac{2\pi}{d+1}(s-|z|)_+^{d+1} & \text{if $\bm{x}_i$ is above/beneath $E_j$}\\
    \displaystyle 0 & \text{otherwise}
  \end{cases}.
\end{eqnarray*}

On the other hand, after the integration with respect to $\rho$, the first term $I^d_n$ in (\ref{eq:Uij_tmp}) reduces to
\begin{eqnarray*}
  I_n^d
  =\int_{\theta_1}^{\theta_2} \frac{-(s-r)_+^{d+1}}{d+1}\diff\theta.
\end{eqnarray*}
This term can be evaluated for $d=1$, $2$, and $3$ as follows:
\begin{eqnarray*}
  I^1_n
  &=& \left[\frac{xy}{2}+\frac{s^2+z^2}{2}\tan^{-1}\frac{x}{y}-sz\tan^{-1}\frac{xz}{yr}-sy\log\frac{x+r}{R}\right]_{x=x_1}^{x=x_2},\\
  I^2_n
  &=&\frac{1}{6}\left[ 2(s^3+3sz^2)\tan^{-1}\frac{x}{y}-2z(3s^2+z^2)\tan^{-1}\frac{xz}{yr}
     -xy(-6s+r)-y(6s^2+y^2+3z^2)\log\frac{x+r}{R}\right]_{x=x_1}^{x=x_2},\nonumber\\
  I^3_n
  &=&\frac{1}{12}\Biggl[
     xy(18s^2+x^2+3y^2+6z^2-6sr)+3(s^4+6s^2z^2+z^4)\tan^{-1}\frac{x}{y}\nonumber\\
    &&-12sz(s^2+z^2)\tan^{-1}\frac{xz}{yr}-6ys(2s^2+y^2+3z^2)\log\frac{x+r}{R}
      \Biggr]_{x=x_1}^{x=x_2}.
\end{eqnarray*}

By considering the relationship $W_{ij}^{(\gamma)} \equiv - \bm{n}\cdot\nabla^x U_{ij}^{(\gamma)}$, which follows from the property $\nabla\Gamma(\bm{x}-\bm{y})=-\nabla^x\Gamma(\bm{x}-\bm{y})$, we may differentiate $U_{ij}^{(\gamma)}$ with respect to the local coordinate $z$ to yield the double-layer potential $W_{ij}^{(\gamma)}$.

\section{Reduction of the discretised ordinary BIE in (\ref{eq:bie_disc3}) to the linear equations in (\ref{eq:linear})}\label{s:linear}

At the current time step $t_\alpha$ ($\alpha=1,2,\ldots$), we solve the discretised BIE in (\ref{eq:bie_disc3}) for the unknown components in the vectors $\mat{u}^{\alpha-1}$ and $\mat{q}^{\alpha-1}$. To this end, we rewrite the BIE in (\ref{eq:bie_disc3}) as follows:
\begin{eqnarray*}
  &&\mat{0}=\left(\sum_{\beta=\beta^*+1}^{\alpha-1}+\sum_{\beta=\alpha}^{\alpha}\right)\left(\mat{U}^{(\alpha-\beta+1)} \bm{\uptau}^{\beta-1} - \mat{W}^{(\alpha-\beta+1)}\bm{\upsigma}^{\beta-1}\right)\qquad\text{\hfill $\because$ The current time step $\alpha$ is separated}\nonumber\\
  &\quad\Leftrightarrow\quad&
  \left(\mat{W}^{(1)} \bm{\upsigma}^{\alpha-1} - \mat{U}^{(1)}\bm{\uptau}^{\alpha-1}\right)=-\sum_{\beta=\beta^*+1}^{\alpha-1}\left(\mat{W}^{(\alpha-\beta+1)}\bm{\upsigma}^{\beta-1} - \mat{U}^{(\alpha-\beta+1)} \bm{\uptau}^{\beta-1}\right)\\
  &\quad\Leftrightarrow\quad&
  w^{0,d} \left(\mat{W}^{(1)}\mat{u}^{\alpha-1} - \mat{U}^{(1)}\mat{q}^{\alpha-1}\right)
  +\left(\mat{W}^{(1)} \left(\sum_{\kappa=1}^{\min(d+1,\alpha-1-\beta^*)}w^{\kappa,d}\mat{u}^{\alpha-1-\kappa}\right) - \mat{U}^{(1)} \left(\sum_{\kappa=1}^{\min(d+1,\alpha-1-\beta^*)}w^{\kappa,d}\mat{q}^{\alpha-1-\kappa}\right)\right)\nonumber\\
  &&=-\sum_{\beta=\beta^*+1}^{\alpha-1}\left(\mat{W}^{(\alpha-\beta+1)}\bm{\upsigma}^{\beta-1} - \mat{U}^{(\alpha-\beta+1)} \bm{\uptau}^{\beta-1}\right)
  \qquad\text{$\because$ Eqs.(\ref{eq:tau}) and (\ref{eq:sigma})}\\
  &\quad\Leftrightarrow\quad&
  \mat{A}\mat{x}^{\alpha-1}
    =-\sum_{\beta=\beta^*+1}^{\alpha}\left(\mat{W}^{(\alpha-\beta+1)}\widetilde{\bm{\upsigma}}^{\beta-1} - \mat{U}^{(\alpha-\beta+1)} \widetilde{\bm{\uptau}}^{\beta-1}\right).
\end{eqnarray*}

\def\since#1{\qquad\text{#1}}
\def\dotL#1{\mat{L}^{(1)\tc{#1}}} 
\def\ddotL#1{\mat{L}^{(2)\tc{#1}}} 
\def\dddotL#1{\mat{L}^{(3)\tc{#1}}} 

\section{Heuristic justification of (\ref{eq:m2l_rec})}\label{s:m2l_rec}

In regard to $d=3$, we will see that the recurrence-type M2L in (\ref{eq:m2l_rec}) is derived from the near- and distant-future M2Ls in (\ref{eq:m2l_near}) and (\ref{eq:m2l_expand}), respectively. The case of $d=3$ differs from $d=2$ only in the existence of the third-order derivatives, i.e. $\mat{L}_{i,k}^{(3)}$ and $\mat{L}_i^{(3)}$. Hence, if we ignore the third-order derivatives from the result for $d=3$, we will be able to obtain the result for $d=2$. Similarly, ignoring the second-order derivatives leads to the result for $d=1$.

In what follows, we let $\mu=2$ (instead of $\mu=8$) for ease of explanation and denote the index of the current time interval by $k$. 

\begin{itemize}
\item $k=0$: We compute $\mat{M}_0$ at the end of the current time interval $I_0$ and cast it to $I_1,\ldots,I_{k+\mu+1}(=I_3)$ by the near-future M2L in (\ref{eq:m2l_near}) to yield $\mat{L}_1,\ldots,\mat{L}_{k+\mu+1}$, respectively, as follows:
  \begin{eqnarray}
    \mat{L}_1^{\tc{0}}&=&\mat{L}_1^{\tc{-1}}+\mat{U}_{1,0}\mat{M}_0^{\tc{0}}=\mat{0}+\mat{L}_{1,0}^{\tc{0}} \since{$\because$ (\ref{eq:m2l_near}) and (\ref{eq:L_l,k}) with $(k,l)=(0,1)$}\nonumber\\
    \mat{L}_2^{\tc{0}}&=&\mat{L}_2^{\tc{-1}}+\mat{U}_{2,0}\mat{M}_0^{\tc{0}}=\mat{0}+\mat{L}_{2,0}^{\tc{0}} \since{$\because$ (\ref{eq:m2l_near}) and (\ref{eq:L_l,k}) with $(k,l)=(0,2)$}\nonumber\\
    \mat{L}_3^{\tc{0}}&=&\mat{L}_3^{\tc{-1}}+\mat{U}_{3,0}\mat{M}_0^{\tc{0}}=\mat{0}+\mat{L}_{3,0}^{\tc{0}} \since{$\because$ (\ref{eq:m2l_near}) and (\ref{eq:L_l,k}) with $(k,l)=(0,3)$}\label{eq:k0L3}
  \end{eqnarray}
  Next, we compute $\mat{L}_4$ by (\ref{eq:m2l_expand}) with $(k,l)=(0,4)$, i.e. the distant-future M2L from $I_0$ to $I_4$ via $I_3$, as follows:
  \begin{eqnarray}
    \mat{L}_4^{\tc{0}}&=&\sum_{p=0}^d\mat{L}_{3,0}^{(p)\tc{0}}1^pT^p \since{$\because$ (\ref{eq:m2l_expand}) with $(k,l)=(0,4)$}\nonumber\\
    &=&\mat{L}_{3,0}^{\tc{0}}+\dotL{0}_{3,0}1^1T^1+\ddotL{0}_{3,0}1^2T^2+\dddotL{0}_{3,0}1^3T^3\label{eq:k0L4}\\
   &=&\mat{L}_3^{\tc{0}}+\dotL{0}_{3,0}(1^1-0^1)T^1+\ddotL{0}_{3,0}(1^2-0^2)T^2+\ddotL{0}_{3,0}(1^3-0^3)T^3 \since{$\because$ (\ref{eq:k0L3})}\nonumber\\
    &=&\mat{L}_3^{\tc{0}}+\mat{L}_{3}^{(1)\tc{0}}T^1+\mat{L}_{3}^{(2)\tc{0}}T^2+\mat{L}_{3}^{(3)\tc{0}}T^3 \since{$\because$ (\ref{eq:L_derivative})} \nonumber\\
    &=&\mat{L}_3^{\tc{0}}+\sum_{p=1}^{3}\mat{L}_{3}^{(p)\tc{0}}T^p\nonumber
  \end{eqnarray}
  The last equation corresponds to (\ref{eq:m2l_rec}) with $k=0$ and $d=3$.

\item $k=1$: We compute $\mat{M}_1$ to update $\mat{L}_2$, $\mat{L}_3$, and $\mat{L}_4$ with (\ref{eq:m2l_near}) as follows::
  \begin{eqnarray}
    \mat{L}_2^{\tc{1}}
    &=&\mat{L}_2^{\tc{0}}+\mat{L}_{2,1}^{\tc{1}} \since{$\because$ (\ref{eq:m2l_near}) and (\ref{eq:L_l,k}) with $(k,l)=(1,2)$}\nonumber\\
    \mat{L}_3^{\tc{1}}
    &=&\mat{L}_3^{\tc{0}}+\mat{L}_{3,1}^{\tc{1}} \since{$\because$ (\ref{eq:m2l_near}) and (\ref{eq:L_l,k}) with $(k,l)=(1,3)$}\nonumber\\
    \mat{L}_4^{\tc{1}}
    &=&\mat{L}_4^{\tc{0}}+\mat{L}_{4,1}^{\tc{1}} \since{$\because$ (\ref{eq:m2l_near}) and (\ref{eq:L_l,k}) with $(k,l)=(1,4)$} \label{eq:k1L4}
  \end{eqnarray}
  Finally, we compute $\mat{L}_5$ by the distant-future M2Ls from $I_0$ to $I_5$ and from $I_1$ to $I_5$ via $I_4$ as follows:
  \begin{eqnarray}
    \mat{L}_5^{\tc{1}}
    &=&\sum_{p=0}^3\mat{L}_{3,0}^{(p)\tc{0}}2^pT^p+\sum_{p=0}^3\mat{L}_{4,1}^{(p)\tc{1}}1^pT^p \since{$\because$ (\ref{eq:m2l_expand}) with $(k,l)=(0,5)$ and $(1,5)$}\nonumber\\
    &=&\left(\mat{L}^{\tc{0}}_{3,0}+\dotL{0}_{3,0}2^1T^1+\ddotL{0}_{3,0}2^2T^2+\dddotL{0}_{3,0}2^3T^3\right)+\left(\mat{L}^{\tc{1}}_{4,1}+\dotL{1}_{4,1}1^1T^1+\ddotL{1}_{4,1}1^2T^2+\dddotL{1}_{4,1}1^3T^3\right) \nonumber\\\label{eq:k1L5}\\
    &=&\mat{L}_4^{\tc{1}}+\left(\dotL{0}_{3,0}(2^1-1^1)+\dotL{1}_{4,1}(1^1-0^1)\right)T^1+\left(\ddotL{0}_{3,0}(2^2-1^2)+\ddotL{1}_{4,1}(1^2-0^2)\right)T^2\nonumber\\
    &&+ \left(\dddotL{0}_{3,0}(2^3-1^3)+\dddotL{1}_{4,1}(1^3-0^3)\right)T^3\since{$\because$ (\ref{eq:k0L4}) and (\ref{eq:k1L4})}\nonumber\\
    &=&\mat{L}_4^{\tc{1}}+\mat{L}_{4}^{(1)\tc{1}}T^1+\mat{L}_{4}^{(2)\tc{1}}T^2 + \mat{L}_{4}^{(3)\tc{1}}T^3\since{$\because$ (\ref{eq:L_derivative})}\nonumber\\
    &=&\mat{L}_4^{\tc{1}}+\sum_{p=1}^3\mat{L}_4^{(p)\tc{1}}T^p\nonumber
  \end{eqnarray}
  The last equation corresponds to (\ref{eq:m2l_rec}) with $k=1$ and $d=3$.

\item $k=2$: Similarly, we compute $\mat{M}_2$, update $\mat{L}_2,\ldots,\mat{L}_6$, and create $\mat{L}_7$ as follows:
  \begin{eqnarray}
    \mat{L}_3^{\tc{2}}
    &=&\mat{L}_3^{\tc{1}}+\mat{L}_{3,2}^{\tc{2}} \since{$\because$ (\ref{eq:m2l_near}) and (\ref{eq:L_l,k}) with $(k,l)=(2,3)$}\nonumber\\
    \mat{L}_4^{\tc{2}}
    &=&\mat{L}_4^{\tc{1}}+\mat{L}_{4,2}^{\tc{2}} \since{$\because$ (\ref{eq:m2l_near}) and (\ref{eq:L_l,k}) with $(k,l)=(2,4)$}\nonumber\\
    \mat{L}_5^{\tc{2}}
    &=&\mat{L}_5^{\tc{1}}+\mat{L}_{5,2}^{\tc{2}} \since{$\because$ (\ref{eq:m2l_near}) and (\ref{eq:L_l,k}) with $(k,l)=(2,5)$}\label{eq:k2L5}\\
    \mat{L}_6^{\tc{2}}&=&\left(\mat{L}_{3,0}^{\tc{0}}+\dotL{0}_{3,0}3^1T^1+\ddotL{0}_{3,0}3^2T^2+\dddotL{0}_{3,0}3^3T^3\right)\nonumber\\
                        &&+\left(\mat{L}^{\tc{1}}_{4,1}+\dotL{1}_{4,1}2^1T^1+\ddotL{1}_{4,1}2^2T^2+\dddotL{1}_{4,1}2^3T^3\right)\nonumber\\
                        &&+\left(\mat{L}^{\tc{2}}_{5,2}+\dotL{2}_{5,2}1^1T^1+\ddotL{2}_{5,2}1^2T^2+\dddotL{2}_{5,2}1^3T^3\right)\since{$\because$ (\ref{eq:m2l_expand}) with $(k,l)=(0,6)$, $(1,6)$, and $(2,6)$}\nonumber\\
    \label{eq:k2L6}\\
    &=&\mat{L}_5^{\tc{2}}+\left(\dotL{0}_{3,0}(3^1-2^1)+\dotL{1}_{4,1}(2^1-1^1)+\dotL{2}_{5,2}(1^1-0^1)\right)T^1\nonumber\\
    &&+\left(\ddotL{0}_{3,0}(3^2-2^2)+\ddotL{1}_{4,1}(2^2-1^2)+\ddotL{2}_{5,2}(1^2-0^2)\right)T^2\nonumber\\
    &&+\left(\dddotL{0}_{3,0}(3^3-2^3)+\dddotL{1}_{4,1}(2^3-1^3)+\dddotL{2}_{5,2}(1^3-0^3)\right)T^3\since{$\because$ (\ref{eq:k1L5}) and (\ref{eq:k2L5})}\nonumber\\
    &=&\mat{L}_5^{\tc{2}}+\mat{L}_5^{(1)\tc{2}}T^1+\mat{L}_5^{(2)\tc{2}}T^2+\mat{L}_5^{(3)\tc{2}}T^3 \since{$\because$ (\ref{eq:L_derivative})}\nonumber\\
    &=&\mat{L}_5^{\tc{2}}+\sum_{p=1}^3\mat{L}_5^{(p)\tc{2}}T^p\nonumber
  \end{eqnarray}
  The last equation corresponds to (\ref{eq:m2l_rec}) with $k=2$ and $d=3$.

\item $k=3$: Similarly, we compute $\mat{M}_3$, update $\mat{L}_3,\ldots,\mat{L}_7$, and create $\mat{L}_8$ as follows:
  \begin{eqnarray}
    \mat{L}_4^{\tc{3}}
    &=&\mat{L}_4^{\tc{2}}+\mat{L}_{4,3}^{\tc{3}} \since{$\because$ (\ref{eq:m2l_near}) and (\ref{eq:L_l,k}) with $(k,l)=(3,4)$}\nonumber\\
    \mat{L}_5^{\tc{3}}
    &=&\mat{L}_5^{\tc{2}}+\mat{L}_{5,3}^{\tc{3}} \since{$\because$ (\ref{eq:m2l_near}) and (\ref{eq:L_l,k}) with $(k,l)=(3,5)$}\nonumber\\
    \mat{L}_6^{\tc{3}}
    &=&\mat{L}_6^{\tc{2}}+\mat{L}_{6,3}^{\tc{3}} \since{$\because$ (\ref{eq:m2l_near}) and (\ref{eq:L_l,k}) with $(k,l)=(3,6)$}\label{eq:k3L6}\\
    \mat{L}_7^{\tc{3}}&=&\left(\mat{L}_{3,0}^{\tc{0}}+\dotL{0}_{3,0}4^1T^1+\ddotL{0}_{3,0}4^2T^2+\dddotL{0}_{3,0}4^3T^3\right)+\left(\mat{L}_{4,1}^{\tc{1}}+\dotL{1}_{4,1}3T^1+\ddotL{1}_{4,1}3^2T^2+\dddotL{1}_{4,1}3^3T^3\right)\nonumber\\
    &&+\left(\mat{L}_{5,2}^{\tc{2}}+\dotL{2}_{5,2}2^1T^1+\ddotL{2}_{5,2}2^2T^2+\dddotL{2}_{5,2}2^3T^3\right)+\left(\mat{L}_{6,3}^{\tc{3}}+\dotL{3}_{6,3}1^1T^1+\ddotL{3}_{6,3}1^2T^2+\dddotL{3}_{6,3}1^3T^3\right)\nonumber\\
    &&\qquad\qquad\since{$\because$ (\ref{eq:m2l_expand}) with $(k,l)=(0,7)$, $(1,7)$, $(2,7)$, and $(3,7)$}\nonumber\\    
    &=&\mat{L}_{6}^{\tc{3}}+\left(\dotL{0}_{3,0}(4^1-3^1)+\dotL{1}_{4,1}(3^1-2^1)+\dotL{2}_{5,2}(2^1-1^1)+\dotL{3}_{6,3}(1^1-0^1)\right)T^1\nonumber\\
    &&+\left(\ddotL{0}_{3,0}(4^2-3^2)+\ddotL{1}_{4,1}(3^2-2^2)+\ddotL{1}_{5,2}(2^2-1^2)+\ddotL{3}_{6,3}(1^2-0^2)\right)T^2\nonumber\\
    &&+\left(\dddotL{0}_{3,0}(4^3-3^3)+\dddotL{1}_{4,1}(3^3-2^3)+\dddotL{2}_{5,2}(2^3-1^3)+\dddotL{3}_{6,3}(1^3-0^3)\right)T^3 \since{$\because$ (\ref{eq:k2L6}) and (\ref{eq:k3L6})}\nonumber\\
    &=&\mat{L}_{6}^{\tc{3}}+ \mat{L}_6^{(1)\tc{3}} T^1
    +\mat{L}_6^{(2)\tc{3}} T^2
    +\mat{L}_6^{(3)\tc{3}} T^3 \since{$\because$ (\ref{eq:L_derivative})}\nonumber\\
    &=&\mat{L}_{6}^{\tc{3}} + \sum_{p=1}^3\mat{L}_6^{(p)\tc{3}} T^p
  \end{eqnarray}
  The last equation corresponds to (\ref{eq:m2l_rec}) with $k=3$ and $d=3$.
  
\end{itemize}

From the above, we can presumably conclude that (\ref{eq:m2l_rec}) holds in general.

\section*{Acknowledgements}

The authors would like to thank Naoshi Nishimura at Kyoto University for useful discussions on the stabilisation by using the Burton--Miller type BIE and Fumito Takase at Nagoya University for his assistance in building the models used in Section~\ref{s:num}. This study was partially supported by the KAKENHI (Grant numbers 18H03251 and 21H03454) and Shimizu Corporation (Project code 2720GZ046c).

\iffalse
\biboptions{numbers,sort&compress} 
\bibliographystyle{elsarticle-num-cpc} 
\bibliography{wave3d_paper}   

\begin{thebibliography}{10}
\expandafter\ifx\csname url\endcsname\relax
  \def\url#1{\texttt{#1}}\fi
\expandafter\ifx\csname urlprefix\endcsname\relax\def\urlprefix{URL }\fi
\expandafter\ifx\csname href\endcsname\relax
  \def\href#1#2{#2} \def\path#1{#1}\fi

\bibitem{takahashi2014}
T.~Takahashi, Journal of Computational Physics 258 (2014) 809--832.
\newline\url{https://www.sciencedirect.com/science/article/pii/S0021999113007584}

\bibitem{greengard1987}
L.~Greengard, V.~Rokhlin, Journal of Computational Physics 73~(2) (1987)
  325--348.

\bibitem{nishimura2002}
N.~Nishimura, Applied Mechanics Reviews 55~(4) (2002) 299--324.

\bibitem{liu2009book}
Y.~Liu, Fast Multipole Boundary Element Method: Theory and Applications in
  Engineering, Cambridge University Press, Cambridge, 2009.

\bibitem{liu2011review}
Y.~J. Liu, S.~Mukherjee, N.~Nishimura, M.~Schanz, W.~Ye, A.~Sutradhar, E.~Pan,
  N.~A. Dumont, A.~Frangi, A.~Saez, Applied Mechanics Reviews 64~(3) (2011)
  030802.

\bibitem{chew2001book}
W.~Chew, E.~Michielssen, J.~M. Song, J.~M. Jin (Eds.), Fast and Efficient
  Algorithms in Computational Electromagnetics, Artech House, Inc., Norwood,
  MA, USA, 2001.

\bibitem{ergin1998}
A.~Ergin, B.~Shanker, E.~Michielssen, Journal of Computational Physics 146~(1)
  (1998) 157--180.

\bibitem{ergin1999a}
A.~Ergin, B.~Shanker, E.~Michielssen, IEEE Antennas and Propagation Magazine
  41~(4) (1999) 39--52.

\bibitem{ergin1999c}
A.~A. Ergin, B.~Shanker, E.~Michielssen, Journal of the Acoustical Society of
  America 106 (1999) 2405--2416.

\bibitem{ergin2000}
A.~A. Ergin, B.~Shanker, E.~Michielssen, Journal of the Acoustical Society of
  America 107 (2000) 1168--1178.

\bibitem{lu2004}
M.~Lu, K.~Yegin, E.~Michielssen, B.~Shanker, Electromagnetics 24~(6) (2004)
  425--449.
\newline\url{http://www.tandfonline.com/doi/abs/10.1080/02726340490479977}

\bibitem{lu2004b}
M.~Lu, E.~Michielssen, B.~Shanker, Electromagnetics 24~(6) (2004) 451--470.
\newline\url{http://www.tandfonline.com/doi/abs/10.1080/02726340490467529}

\bibitem{shanker2003}
B.~Shanker, A.~Ergin, M.~Lu, E.~Michielssen, IEEE Transactions on Antennas and
  Propagation 51~(3) (2003) 628--641.

\bibitem{aygun2004}
K.~Aygun, B.~Fischer, J.~Meng, B.~Shanker, E.~Michielssen, IEEE Transactions on
  Microwave Theory and Techniques 52~(2) (2004) 573--583.

\bibitem{takahashi2001}
T.~Takahashi, N.~Nishimura, S.~Kobayashi, Transaction of the JSME, Series A
  67~(661) (2001) 1409--1416, (written in Japanese).

\bibitem{takahashi2003}
T.~Takahashi, N.~Nishimura, S.~Kobayashi, Engineering Analysis with Boundary
  Elements 27~(5) (2003) 491--506.

\bibitem{liu2014parallel_wavelet_PWTD}
Y.~{Liu}, A.~C. {Yücel}, H.~{Bağcı}, E.~{Michielssen}, A parallel
  wavelet-enhanced pwtd algorithm for analyzing transient scattering from
  electrically very large pec targets, in: 2014 USNC-URSI Radio Science Meeting
  (Joint with AP-S Symposium), 2014, pp. 177--177.
\newblock

\bibitem{hargreaves2007}
J.~Hargreaves, \href{http://usir.salford.ac.uk/id/eprint/16604/}{Time domain
  boundary element method for room acoustics}, Ph.D. thesis, University of
  Salford (April 2007).
\newline\url{http://usir.salford.ac.uk/id/eprint/16604/}

\bibitem{soares2007}
D.~Soares, W.~J. Mansur, {Computational Mechanics} 40~(2).

\bibitem{okamura2016}
R.~Okamura, H.~Yoshikawa, T.~Takahashi, T.~Takagi, K.~Kashiyama, Journal of
  Japan Society of Civil Engineers, Ser. A2 (Applied Mechanics) 72~(2) (2016)
  I\_257--I\_264, (written in Japanese).

\bibitem{ergin1999b}
A.~A. Ergin, B.~Shanker, E.~Michielssen, The Journal of the Acoustical Society
  of America 106~(5) (1999) 2396--2404.
\newline\url{https://doi.org/10.1121/1.428076}

\bibitem{burton_miller1971}
A.~J. Burton, G.~F. Miller, Proceedings of the Royal Society of London A:
  Mathematical, Physical and Engineering Sciences 323~(1553) (1971) 201--210.
\newline\url{http://rspa.royalsocietypublishing.org/content/323/1553/201}

\bibitem{fukuhara2019}
M.~Fukuhara, R.~Misawa, K.~Niino, N.~Nishimura, Engineering Analysis with
  Boundary Elements 108 (2019) 321--338.
\newline\url{https://www.sciencedirect.com/science/article/pii/S0955799719305600}

\bibitem{asakura2009}
J.~Asakura, T.~Sakurai, H.~Tadano, T.~Ikegami, K.~Kimura, JSIAM Letters 1
  (2009) 52--55.

\bibitem{chiyoda2019}
S.~Chiyoda, K.~Niino, N.~Nishimura, Proceedings of the Conference on
  Computational Engineering and Science 24 (2019) 3p, (written in Japanese).
\newline\url{https://ci.nii.ac.jp/naid/40021915153/en/}

\bibitem{aimi2009}
A.~Aimi, M.~Diligenti, C.~Guardasoni, I.~Mazzieri, S.~Panizzi, International
  Journal for Numerical Methods in Engineering 80~(9) (2009) 1196--1240.
\newline\url{https://onlinelibrary.wiley.com/doi/abs/10.1002/nme.2660}

\bibitem{joly2017}
P.~Joly, J.~Rodr$\acute{\i}$guez, Journal of Integral Equations and
  Applications 29~(1) (2017) 137 -- 187.
\newline\url{https://doi.org/10.1216/JIE-2017-29-1-137}

\bibitem{gimperlein2018}
H.~Gimperlein, C.~\"Ozdemir, E.~P. Stephan, Journal of Computational
  Mathematics 36~(1) (2018) 70--89.
\newline\url{http://global-sci.org/intro/article_detail/jcm/10583.html}

\bibitem{yoshikawa2003phd}
H.~Yoshikawa, Study on the application of the time-domain boundary integral
  equation method to the non-destructive evaluation using laser ultrasonic
  measurement, Ph.D. thesis, Kyoto university, (written in Japanese) (September
  2003).

\bibitem{Bowman}
J.~J. Bowman, T.~B.~A. Senior, P.~L.~E. Uslenghi, {Electromagnetic and Acoustic
  Scattering by Simple Shapes (Revised edition)}, Hemisphere Publishing Corp.,
  New York, 1987.

\bibitem{geuzaine2009Gmsh}
C.~Geuzaine, J.-F. Remacle, International Journal for Numerical Methods in
  Engineering 79~(11) (2009) 1309--1331.
\newline\url{http://dx.doi.org/10.1002/nme.2579}

\bibitem{powell2007}
M.~J.~D. Powell, Cambridge Uinversity Technical Report (2007) 10--12.
\newline\url{http://www.damtp.cam.ac.uk/user/na/NA_papers/NA2007_03.pdf}

\bibitem{jung2003}
B.~H. Jung, Y.-S. Chung, T.~K. Sarkar, Journal of Electromagnetic Waves and
  Applications 17~(5) (2003) 737--739.
\newline\url{https://doi.org/10.1163/156939303322226383}

\end{thebibliography}
\else

\fi

\end{document}